\input amstexl


\catcode`\@=11
\ifx\amstexloaded@\relax\else
 \errmessage{AmS-TeX must be loaded before LamS-TeX}\fi
\ifx\laxread@\undefined\else\catcode`\@=\active \fi
\def\err@#1{\errmessage{LamS-TeX error: #1}}
\def^^L{\par}
\let\+\tabalign
\def\newcount{\alloc@0\count\countdef\insc@unt}
\def\newdimen{\alloc@1\dimen\dimendef\insc@unt}
\def\newskip{\alloc@2\skip\skipdef\insc@unt}
\def\newmuskip{\alloc@3\muskip\muskipdef\@cclvi}
\def\newbox{\alloc@4\box\chardef\insc@unt}
\let\newtoks\relax
\def\newhelp#1#2{\newtoks#1#1\expandafter{\csname#2\endcsname}}
\def\newtoks{\alloc@5\toks\toksdef\@cclvi}
\def\newread{\alloc@6\read\chardef\sixt@@n}
\def\newwrite{\alloc@7\write\chardef\sixt@@n}
\def\newfam{\alloc@8\fam\chardef\sixt@@n}
\def\newlanguage{\alloc@9\language\chardef\@cclvi}
\def\newinsert#1{\global\advance\insc@unt by\m@ne
  \ch@ck0\insc@unt\count
  \ch@ck1\insc@unt\dimen
  \ch@ck2\insc@unt\skip
  \ch@ck4\insc@unt\box
  \allocationnumber=\insc@unt
  \global\chardef#1=\allocationnumber
  \wlog{\string#1=\string\insert\the\allocationnumber}}
\def\newif#1{\count@\escapechar \escapechar\m@ne
  \expandafter\expandafter\expandafter
   \edef\@if#1{true}{\let\noexpand#1=\noexpand\iftrue}%
  \expandafter\expandafter\expandafter
   \edef\@if#1{false}{\let\noexpand#1=\noexpand\iffalse}%
  \@if#1{false}\escapechar\count@}

\def\Err@#1{\errhelp\defaulthelp@\err@{#1}}
{\catcode`\@=\active
 \edef\next{\gdef\noexpand@{\futurelet\noexpand\next
  \csname at\string@\endcsname}}
 \next
}
\def\at@{\ifcat\noexpand\next a\let\next@\at@@\else
 \ifcat\noexpand\next0\let\next@\at@@\else
 \ifcat\noexpand\next\relax\let\next@\at@@\else
 \let\next@\at@@@\fi\fi\fi\next@}
\def\at@@@{\errhelp\athelp@\err@{Invalid use of @}}
\def\at@@#1{\expandafter
 \ifx\csname\string#1@at\endcsname\relax\let\next@\at@@@\else
 \DN@{\csname\string#1@at\endcsname}\fi\next@}
\def\atdef@#1{\expandafter\def\csname\string#1@at\endcsname}
\newif\iftest@
\def\tagin@#1{\tagin@false
 \DN@##1\tag##2##3\next@{\test@true\ifx\tagin@##2\test@false\fi}%
 \next@#1\tag\tagin@\next@\tagin@false\iftest@\tagin@true\fi}
\let\lkerns@\relax
\def\nolinebreak{\RIfM@\mathmodeerr@\nolinebreak\else
 \ifhmode\saveskip@\lastskip\unskip
 \nobreak\ifdim\saveskip@>\z@\hskip\saveskip@\fi\lkerns@
 \else\vmodeerr@\nolinebreak\fi\fi}
\def\allowlinebreak{\RIfM@\mathmodeerr@\allowlinebreak\else
 \ifhmode\saveskip@\lastskip\unskip
 \allowbreak\ifdim\saveskip@>\z@\hskip\saveskip@\fi\lkerns@
 \else\vmodeerr@\allowlinebreak\fi\fi}
\def\linebreak{\RIfM@\mathmodeerr@\linebreak\else
 \ifhmode\unskip\unkern\break\lkerns@
 \else\vmodeerr@\linebreak\fi\fi}
\let\nkerns@\relax
\def\newline{\RIfM@\mathmodeerr@\newline\else
 \ifhmode\unskip\unkern\null\hfill\break\nkerns@
 \else\vmodeerr@\newline\fi\fi}%
\def\newbox@{\alloc@@4\box\chardef\insc@unt}
\def\newcount@{\alloc@@0\count\countdef\insc@unt}
\def\accentedsymbol#1#2{\expandafter\newbox@\csname\exstring@#1@box\endcsname
 \setbox\csname\exstring@#1@box\endcsname\hbox{$\m@th#2$}%
 \define#1{\copy\csname\exstring@#1@box\endcsname{}}}
\def\rightadd@#1\to#2{\toks@{\\#1}\toks@@\expandafter{#2}\xdef#2{\the\toks@@
 \the\toks@}\toks@{}\toks@@{}}
\def\fontlist@{\\\tenrm\\\sevenrm\\\fiverm\\\teni\\\seveni\\\fivei
 \\\tensy\\\sevensy\\\fivesy\\\tenex\\\tenbf\\\sevenbf\\\fivebf
 \\\tensl\\\tenit}
\def\font@#1=#2 {\rightadd@#1\to\fontlist@\font#1=#2 }
\def\ismember@#1#2{\global\let\Next@ F\let\next@= #2%
 {\def\\##1{\let\nextii@##1\ifx\nextii@\next@\global\let\Next@ T\fi}#1}%
 \test@false\ifx\Next@ T\test@true\fi\let\next@\relax}
\def\FNSS@#1{\let\FNSS@@#1\FN@\FNSS@@@}
\def\FNSS@@@{\ifx\next\space@\def\FNSS@@@@. {\FN@\FNSS@@@}\else
 \def\FNSS@@@@.{\FNSS@@}\fi\FNSS@@@@.}
\atdef@"{\unskip
 \DN@{\ifx\next`\DN@`{\FN@\nextii@}%
  \else\ifx\next\lq\DN@\lq{\FN@\nextii@}%
  \else\DN@####1{\FN@\nextiii@}\fi\fi
  \next@}%
 \DNii@{\ifx\next`\DN@`{\sldl@``}%
  \else\ifx\next\lq\DN@\lq{\sldl@``}%
  \else\DN@{\dlsl@`}\fi\fi\next@}%
 \def\nextiii@{\ifx\next'\DN@'{\srdr@''}%
  \else\ifx\next\rq\DN@\rq{\srdr@''}%
  \else\DN@{\drsr@'}\fi\fi\next@}%
 \FNSS@\next@}
\def\root{%
  \DN@{\ifx\next\uproot\let\next@\nextii@\else
   \ifx\next\leftroot\let\next@\nextiii@\else
   \let\next@\plainroot@\fi\fi\next@}%
  \DNii@\uproot##1{\uproot@##1\relax\FNSS@\nextiv@}%
  \def\nextiv@{\ifx\next\leftroot\let\next@\nextv@\else
   \let\next@\plainroot@\fi\next@}%
  \def\nextv@\leftroot##1{\leftroot@##1\relax\plainroot@}%
  \def\nextiii@\leftroot##1{\leftroot@##1\relax\FNSS@\nextvi@}%
  \def\nextvi@{\ifx\next\uproot\let\next@\nextvii@\else
   \let\next@\plainroot@\fi\next@}%
  \def\nextvii@\uproot##1{\uproot@##1\relax\plainroot@}%
  \bgroup\uproot@\z@\leftroot@\z@
 \FNSS@\next@}
\def\loop#1\repeat{\def\iterate{#1\relax\expandafter\iterate\fi}%
 \iterate\let\iterate\relax}
\def\gloop@#1\repeat{\gdef\iterate@{#1\relax\expandafter\iterate@\fi}%
 \iterate@\global\let\iterate@\relax}
\def\printoptions{\W@{Do you want S(yntax check),
  G(alleys) or P(ages)?^^JType S, G or P, follow by <return>: }\loop
 \read\m@ne to\ans@
 \edef\next@{\def\noexpand\Ans@{\ans@}}\uppercase\expandafter{\next@}%
 \ifx\Ans@\S@\test@true\syntax\else
 \ifx\Ans@\G@\test@true\galleys\else
 \ifx\Ans@\P@\test@true\else
 \test@false\fi\fi\fi
 \iftest@\else\W@{Type S, G or P, follow by <return>: }%
 \repeat}
\expandafter\let\csname A@;\endcsname;
\expandafter\let\csname A@:\endcsname:
\expandafter\let\csname A@?\endcsname?
\expandafter\let\csname A@!\endcsname!
\def\APdef#1{\def\next@{\expandafter\let\csname A@\string#1\endcsname#1}%
 \afterassignment\next@\def#1}
\let\fextra@\,
\def\tdots@{\unskip
 \DN@{$\m@th\mathinner{\ldotp\ldotp\ldotp}\,
   \ifx\next,\,$\else\ifx\next.\,$\else
   \ifx\next;\,$\else
   \expandafter\ifx\csname A@\string;\endcsname\next\fextra@$\else
   \ifx\next:\,$\else
   \expandafter\ifx\csname A@\string:\endcsname\next\fextra@$\else
   \ifx\next?\,$\else
   \expandafter\ifx\csname A@\string?\endcsname\next\fextra@$\else
   \ifx\next!\,$\else
   \expandafter\ifx\csname A@\string!\endcsname\next\fextra@$\else
   $ \fi\fi\fi\fi\fi\fi\fi\fi\fi\fi}%
 \ \FN@\next@}
\def\extrap@#1{%
 \ifx\next,\DN@{#1\,}\else
 \ifx\next;\DN@{#1\,}\else
 \expandafter\ifx\csname A@\string;\endcsname\next\DN@{#1\fextra@}\else
 \ifx\next.\DN@{#1\,}\else\extra@
 \ifextra@\DN@{#1\,}\else
 \let\next@#1\fi\fi\fi\fi\fi\next@}
\def\dotsc{\DN@{\ifx\next;\plainldots@\,\else
 \expandafter\ifx\csname A@\string;\endcsname\next\plainldots@\fextra@\else
 \ifx\next.\plainldots@\,\else\extra@\plainldots@
 \ifextra@\,\fi\fi\fi\fi}%
 \FN@\next@}
\def\keybin@{\keybin@true
 \ifx\next+\else\ifx\next=\else\ifx\next<\else\ifx\next>\else\ifx\next-\else
 \ifx\next*\else\ifx\next:\else
 \expandafter\ifx\csname A@\string;\endcsname\next\else
 \keybin@false\fi\fi\fi\fi\fi\fi\fi\fi}
\def\boldkey#1{\ifcat\noexpand#1A%
  \ifcmmibloaded@{\fam\cmmibfam#1}\else
   \Err@{First bold symbol font not loaded}\fi
 \else
 \let\next=#1%
 \ifx#1!\mathchar"5\bffam@21 \else
 \expandafter\ifx\csname A@\string!\endcsname\next\mathchar"5\bffam@21 \else
 \ifx#1(\mathchar"4\bffam@28 \else\ifx#1)\mathchar"5\bffam@29 \else
 \ifx#1+\mathchar"2\bffam@2B \else\ifx#1:\mathchar"3\bffam@3A \else
 \expandafter\ifx\csname A@\string:\endcsname\next\mathchar"3\bffam@3A \else
 \ifx#1;\mathchar"6\bffam@3B \else
 \expandafter\ifx\csname A@\string;\endcsname\next\mathchar"6\bffam@3B \else
 \ifx#1=\mathchar"3\bffam@3D \else
 \ifx#1?\mathchar"5\bffam@3F \else
 \expandafter\ifx\csname A@\string?\endcsname\next\mathchar"5\bffam@3F \else
 \ifx#1[\mathchar"4\bffam@5B \else
 \ifx#1]\mathchar"5\bffam@5D \else
 \ifx#1,\mathchari@63B \else
 \ifx#1-\mathcharii@200 \else
 \ifx#1.\mathchari@03A \else
 \ifx#1/\mathchari@03D \else
 \ifx#1<\mathchari@33C \else
 \ifx#1>\mathchari@33E \else
 \ifx#1*\mathcharii@203 \else
 \ifx#1|\mathcharii@06A \else
 \ifx#10\bold0\else\ifx#11\bold1\else\ifx#12\bold2\else\ifx#13\bold3\else
 \ifx#14\bold4\else\ifx#15\bold5\else\ifx#16\bold6\else\ifx#17\bold7\else
 \ifx#18\bold8\else\ifx#19\bold9\else
  \Err@{\noexpand\boldkey can't be used with #1}%
 \fi\fi\fi\fi\fi\fi\fi\fi\fi\fi\fi\fi\fi\fi\fi
 \fi\fi\fi\fi\fi\fi\fi\fi\fi\fi\fi\fi\fi\fi\fi\fi\fi\fi}
\def\arabic#1{#1}
\def\alph#1{\count@#1\relax\advance\count@96 \ifnum\count@>122
 \Err@{\noexpand\alph invalid for numbers > 26}\else\char\count@\fi}
\def\Alph#1{\count@#1\relax\advance\count@64 \ifnum\count@>90
 \Err@{\noexpand\Alph invalid for numbers > 26}\else\char\count@\fi}

\def\Roman#1{\uppercase\expandafter{\romannumeral#1}}
\def\fnsymbol#1{\count@#1\relax
 \count@@\count@
 \advance\count@\m@ne\divide\count@7
 \count@@@\count@\advance\count@@@\@ne
 \multiply\count@7 \advance\count@@-\count@
 \count@\count@@@
 {\loop
  \ifcase\count@@\or*\or\dag\or\ddag\or\P\or\S\or\text{$\|$}\or\#\fi
  \advance\count@\m@ne\ifnum\count@>\z@\repeat}}
\def\cardnine@#1{\ifcase#1\or one\or two\or three\or four\or five\or
 six\or seven\or eight\or nine\fi}
\let\alloc@\alloc@@
\newcount\ten@
\ten@10
\def\cardinal#1{\count@#1\relax
 \ifnum\count@>99 \number\count@
 \else
  \ifnum\count@=\z@ zero%
  \else
   \ifnum\count@<\ten@\cardnine@\count@
   \else
    \ifnum\count@<20
     \advance\count@-\ten@
     \ifcase\count@ ten\or eleven\or twelve\or thirteen\or fourteen\or
      fifteen\or sixteen\or seventeen\or eighteen\or nineteen\fi
    \else
     \count@@\count@\count@@@\count@@
     \divide\count@\ten@\multiply\count@\ten@
     \advance\count@@@-\count@\divide\count@\ten@
     \ifcase\count@\or\or twenty\or thirty\or forty\or fifty\or sixty\or
      seventy\or eighty\or ninety\fi
     \ifnum\count@@@=\z@\else-\cardnine@\count@@@\fi
    \fi
   \fi
  \fi
 \fi}
\def\ordnine@#1{\ifcase#1\or first\or second\or third\or fourth\or fifth\or
 sixth\or seventh\or eighth\or ninth\fi}
\newcount\count@@@@
\def\ordsuffix@{\count@@@@\count@
 \divide\count@\ten@
 \count@@@\count@\count@@\count@
 \divide\count@@\ten@\multiply\count@@\ten@
 \advance\count@@@-\count@@
 \ifnum\count@@@=\@ne th%
 \else
  \count@@@\count@@@@
  \count@@\count@@@@
  \divide\count@@\ten@\multiply\count@@\ten@
  \advance\count@@@-\count@@
  \ifcase\count@@@ th\or st\or nd\or rd\else th\fi
 \fi}
\def\nordinal#1{\count@#1\relax\number\count@\ordsuffix@}
\def\spordinal#1{\count@#1\relax\number\count@$^{\text{\ordsuffix@}}$}
\def\ordinal#1{\count@#1\relax
 \ifnum\count@>99 \number\count@\ordsuffix@
 \else
   \ifnum\count@=\z@ zeroth%
  \else
    \ifnum\count@<\ten@\ordnine@\count@
    \else
     \ifnum\count@<20 \advance\count@-\ten@
      \ifcase\count@ tenth\or eleventh\or twelfth\or thirteenth\or
       fourteenth\or fifteenth\or sixteenth\or seventeenth\or eighteenth\or
       nineteenth\fi
     \else
      \count@@\count@
      \divide\count@\ten@\multiply\count@\ten@
      \count@@@\count@@\advance\count@@@-\count@
      \divide\count@\ten@
      \ifcase\count@\or\or twent\or thirt\or fort\or fift\or sixt\or sevent\or
       eight\or ninet\fi
      \ifnum\count@@@=\z@ ieth\else y-\ordnine@\count@@@\fi
     \fi
    \fi
  \fi
 \fi}
\font@\tensmc=cmcsc10
\textonlyfont@\smc\tensmc
\newtoks\noexpandtoks@
\noexpandtoks@{\let\arabic\relax\let\alph\relax\let\Alph\relax
 \let\Roman\relax\let\fnsymbol\relax\let\rm\relax
 \let\it\relax\let\bf\relax\let\sl\relax\let\smc\relax
 \let\/\relax\let\null\relax}
\def\noexpands@{\the\noexpandtoks@}
\def\Nonexpanding#1{\global\noexpandtoks@
 \expandafter{\the\noexpandtoks@\let#1\relax}}
\def\prevanish@{\saveskip@\z@\ifhmode\saveskip@\lastskip\unskip\fi}
\def\postvanish@{\ifdim\saveskip@>\z@\hskip\saveskip@\fi\FN@\postvanish@@}
\def\postvanish@@{\DN@.{}%
 \ifx\next\space@\ifdim\saveskip@>\z@\DN@. {}\fi\fi\next@.}
\def\invisible#1{\prevanish@\ignorespaces#1\unskip\postvanish@}
\def\vanishlist@{\\\invisible}
\let\noindent@\noindent
\def\noindent{\par\noindent@\FN@\pretendspace@}
\def\pretendspace@{\ismember@\vanishlist@\next
 \iftest@\nobreak\hskip-\p@\hskip\p@\fi}

\newtoks\everypartoks@
\def\noindent@@{\par\everypartoks@\expandafter{\the\everypar}\everypar{}%
 \noindent@\everypar\expandafter{\the\everypartoks@}}
\def\page{\Err@{\noexpand\page has no meaning by itself}}
\let\page@C\pageno
\let\page@P\empty
\let\page@Q\empty
\def\page@S#1{#1\/}
\def\page@F{\rm}
\def\page@N{\arabic}   
\newif\ifindexing@
\def\indexfile{\ifindexing@\else
 \alloc@@7\write\chardef\sixt@@n\ndx@
 \immediate\openout\ndx@=\jobname.ndx
 \global\indexing@true\fi}
\global\advance\insc@unt\m@ne
\ch@ck0\insc@unt\count
\ch@ck1\insc@unt\dimen
\ch@ck2\insc@unt\skip
\ch@ck4\insc@unt\box
\allocationnumber\insc@unt
\global\chardef\margin@\allocationnumber
\dimen\margin@\maxdimen
\count\margin@\z@
\skip\margin@\z@
\newif\ifindexproofing@
\def\indexproofing{\indexproofing@true}
\def\noindexproofing{\indexproofing@false}
\def\unmacro@#1:#2->#3\unmacro@{\def\macpar@{#2}\def\macdef@{#3}}
\def\starparts@#1{\def\stari@{#1}\def\starii@{#1}\let\stariii@\empty
 \test@false
 \DN@##1*##2##3\next@{\ifx\starparts@##2\test@false\else\test@true\fi}%
 \next@#1*\starparts@\next@
 \iftest@\DN@{\starparts@@#1\starparts@@}\else\let\next@\relax\fi\next@}
\def\starparts@@#1*#2\starparts@@{\def\starii@{#1}\def\stariii@{*#2}}
\def\windex@{\ifindexing@
 \expandafter\unmacro@\meaning\stari@\unmacro@
 \edef\macdef@{\string"\macdef@\string"}%
 \edef\next@{\write\ndx@{\macdef@}}\next@
 \write\ndx@{{\number\pageno}{\page@N}{\page@P}{\page@Q}}%
 \fi
 \ifindexproofing@
  \ifx\stariii@\empty\else
   \expandafter\unmacro@\meaning\stariii@\unmacro@\fi
  \insert\margin@{\hbox{\rm\vrule\height9\p@\depth2\p@\width\z@\starii@
  \ifx\stariii@\empty\else\tt\macdef@\fi}}\fi}
\catcode`\"=\active
\def"{\FN@\quote@}
\def\quote@{\ifx\next"\expandafter\quote@@\else\expandafter\quote@@@\fi}
\def\quote@@@#1"{\starparts@{#1}\starii@\windex@}
\def\quote@@"#1"{\prevanish@\starparts@{#1}\windex@\FN@\quote@@@@}
\def\quote@@@@{\ifx\next"\DN@"{\postvanish@}\else
 \let\next@\postvanish@\fi\next@}
\rightadd@"\to\vanishlist@
\def\idefine#1{\DN@{#1}\DNii@{\noexpand#1}%
 \afterassignment\idefine@\def\nextiii@}
\def\idefine@{\ifindexing@
 \expandafter\let\next@\nextiii@
 \expandafter\unmacro@\meaning\nextiii@\unmacro@
 \immediate\write\ndx@{\noexpand\define\nextii@\macpar@{\macdef@}}\fi}
\def\iabbrev*#1#2{\ifindexing@\toks@{#2}%
 \immediate\write\ndx@{\noexpand\abbrev*\noexpand#1{\the\toks@}}\fi}
\newread\laxread@
\newwrite\laxwrite@
\let\fnpages@\empty
\def\Finit@#1#2\Finit@{\let\nextii@#1\def\nextiii@{#2}}
\catcode`\~=11
\def\getparts@ @#1~#2~#3~#4~#5~#6{\def\nextiv@{#1}%
 \def\nextiii@{#2~#3~#4~#5~}\count@#6\relax}
\newif\ifdocument@
\def\document{\ifdocument@\else\global\document@true
 \let\fontlist@\empty
 \immediate\openin\laxread@=\jobname.lax\relax
 {\endlinechar\m@ne\noexpands@\catcode`\@=11 \catcode`\~=11
  \loop\ifeof\laxread@\else
   \read\laxread@ to\next@
   \ifx\next@\empty
   \else
    \expandafter\Finit@\next@\Finit@
    \if\nextii@ F%
     \expandafter\rightadd@\nextiii@\to\fnpages@
    \else
     \expandafter\getparts@\next@
     \edef\next@{\gdef\csname\nextiv@ @L\endcsname{\nextiii@\number\count@}}%
     \next@
    \fi
   \fi
  \repeat}%
 \immediate\closein\laxread@
 \immediate\openout\laxwrite@=\jobname.lax\relax\fi}
\let\thelabel@\relax
\def\thelabels@{\thelabel@ ~\thelabel@@ ~\thelabel@@@ ~\thelabel@@@@ ~}
\def\label#1{\prevanish@
 \ifx\thelabel@\relax
  \Err@{There's nothing here to be labelled}%
 \else
  {\noexpands@
  \expandafter\ifx\csname#1@L\endcsname\relax
   \expandafter\xdef\csname#1@L\endcsname{\thelabels@0}%
   \immediate\write\laxwrite@{@#1~\thelabels@1}%
  \else
   \edef\next@{@~\csname#1@L\endcsname}%
    \expandafter\getparts@\next@
    \ifodd\count@
    \expandafter\xdef\csname#1@L\endcsname{\thelabels@0}%
    \immediate\write\laxwrite@{@#1~\thelabels@1}%
   \else
    \Err@{Label #1 already used}%
   \fi
  \fi
  }%
 \fi
 \postvanish@}
\rightadd@\label\to\vanishlist@
\def\thepages@{\page@N{\number\page@C}~%
 \page@S{\page@P\page@N{\number\page@C}\page@Q}~%
 \number\page@C ~\page@P\page@N{\number\page@C}\page@Q ~}
\def\pagelabel#1{\prevanish@
 \expandafter\ifx\csname#1@L\endcsname\relax
  {\noexpands@
  \expandafter\xdef\csname#1@L\endcsname{\thepages@2}}%
  \write\laxwrite@{@#1~\thepages@3}%
 \else
  {\noexpands@
  \edef\next@{@~\csname#1@L\endcsname}%
  \expandafter\getparts@\next@
  \ifodd\count@
   \ifnum\count@=\@ne
    \expandafter\xdef\csname#1@L\endcsname{\thelabels@2}%
   \fi
   \write\laxwrite@{@#1~\thepages@3}%
  \else
   \Err@{Label #1 already used}%
  \fi
  }%
 \fi
 \postvanish@}
\rightadd@\pagelabel\to\vanishlist@
\newif\ifreferr@
\referr@true
\def\RefErrors{\global\referr@true}
\def\RefWarnings{\global\referr@false}
\setbox\z@\hbox{\global\count@=`^^30}
\ifnum\count@=48 \let\versionthree@\relax\fi
\def\nolabel@#1#2#3{\expandafter\ifx\csname#2@L\endcsname\relax
 \ifreferr@\Err@{No \noexpand\label found for #2}\else
 \W@{Warning: No \noexpand\label found for #2.}%
 \ifx\versionthree@\relax\W@{l.\number\inputlineno\space ... \string#1{#2}}\fi
 \fi#3\else}
\def\csL@#1{{\noexpands@\xdef\Next@{\csname#1@L\endcsname}}}
\def\ref#1{\nolabel@\ref{#1}\relax
 \DNii@##1~##2\nextii@{##1}%
 \csL@{#1}\expandafter\nextii@\Next@\nextii@\fi}
\def\Ref#1{\nolabel@\Ref{#1}\relax
 \DNii@##1~##2~##3\nextii@{##2}%
 \csL@{#1}\expandafter\nextii@\Next@\nextii@\fi}
\def\nref#1{\nolabel@\nref{#1}\relax
 \DNii@##1~##2~##3~##4\nextii@{##3}%
 \csL@{#1}\expandafter\nextii@\Next@\nextii@\fi}
\def\pref#1{\nolabel@\pref{#1}\relax
 \DNii@##1~##2~##3~##4~##5\nextii@{##4}%
 \csL@{#1}\expandafter\nextii@\Next@\nextii@\fi}
\let\pref@\pref
\def\Evaluatenref#1{\nolabel@\Evaluatenref{#1}{\gdef\Nref{-10000 }}%
 \DNii@##1~##2~##3~##4\nextii@{\DNii@{##3}}%
 \csL@{#1}\expandafter\nextii@\Next@\nextii@
 \xdef\Nref{\nextii@}\fi}
\def\Evaluatepref#1{\nolabel@\Evaluatepref{#1}{\global\let\Pref\empty}%
 \DNii@##1~##2~##3~##4~##5\nextii@{\DNii@{##4}}%
 \csL@{#1}\expandafter\nextii@\Next@\nextii@
 \xdef\Pref{\nextii@}\fi}
\def\readlax#1{\immediate\openin\laxread@=#1.lax\relax
 \ifeof\laxread@\W@{}\W@{File #1.lax not found.}\W@{}\fi
 {\endlinechar\m@ne\noexpands@\catcode`\@=11 \catcode`\~=11
  \loop\ifeof\laxread@\else
   \read\laxread@ to\nextv@
   \ifx\nextv@\empty
   \else
    \expandafter\Finit@\nextv@\Finit@
    \ifx\nextii@ F%
    \else
     \expandafter\getparts@\nextv@
     \expandafter\ifx\csname\nextiv@ @L\endcsname\relax
      \edef\next@{\gdef\csname\nextiv@ @L\endcsname
       {\nextiii@\ifnum\count@=\@ne0\else2\fi}}%
      \next@
     \else
      \Err@{Label \nextiv@\space in #1.lax already used}%
     \fi
    \fi
   \fi
  \repeat}%
 \immediate\closein\laxread@}
\catcode`\~=\active

\def\FNSSP@{\FNSS@\pretendspace@}
\everydisplay{\csname displaymath \endcsname}
\expandafter\def\csname displaymath \endcsname#1$${#1$$\FNSSP@}
\def\locallabel@{\let\thelabel@\Thelabel@\let\thelabel@@\Thelabel@@
 \let\thelabel@@@\Thelabel@@@\let\thelabel@@@@\Thelabel@@@@}
\newcount\tag@C
\tag@C\z@
\let\tag@P\empty
\let\tag@Q\empty
\def\tag@S#1{{\rm(}{#1\/}{\rm)}}
\let\tag@N\arabic
\def\tag@F{\rm}
\def\maketag@{\FN@\maketag@@}
\def\maketag@@{\ifx\next\relax\DN@\relax{\FN@\maketag@@}\else
 \ifx\next"\let\next@\maketag@@@\else
 \let\next@\maketag@@@@\fi\fi\next@}
\def\xdefThelabel@#1{\xdef\Thelabel@{#1{\Thelabel@@@}}}
\def\xdefThelabel@@#1{\xdef\Thelabel@@{#1{\Thelabel@@@@}}}
\def\maketag@@@@#1\maketag@{\global\advance\tag@C\@ne
 {\noexpands@
  \xdef\Thelabel@@@{\number\tag@C}%
  \xdefThelabel@\tag@N
  \xdef\Thelabel@@@@{\ifmathtags@$\tag@P\Thelabel@\tag@Q$\else
   \tag@P\Thelabel@\tag@Q\fi}%
  \xdefThelabel@@\tag@S
  }%
 \locallabel@
 \hbox{\tag@F\thelabel@@}%
 #1}
\def\Qlabel@#1{{\noexpands@\xdef\Thelabel@@{#1}%
 \let\style\empty\xdef\Thelabel@@@@{#1}%
 \let\pre\empty\let\post\empty\xdef\Thelabel@{#1}%
 \let\numstyle\empty\xdef\Thelabel@@@{#1}}}
\def\maketag@@@"#1"#2\maketag@{%
 {\let\pre\tag@P\let\post\tag@Q\let\style\tag@S\let\numstyle\tag@N
  \hbox{\tag@F#1}%
  \noexpands@
  \Qlabel@{#1}%
  }%
 \locallabel@
 #2}
\def\align@{\inalign@true\inany@true
 \vspace@\allowdisplaybreak@\displaybreak@\intertext@
 \def\tag{\global\tag@true\ifnum\and@=\z@
  \DN@{&\omit\global\rwidth@\z@&\relax}\else
  \DN@{&\relax}\fi\next@}%
 \iftagsleft@\DN@{\csname align \endcsname}\else
  \DN@{\csname align \space\endcsname}\fi\next@}
\def\noset@{\def\Offset##1##2{\prevanish@\postvanish@}%
 \def\Reset##1##2{\prevanish@\postvanish@}}
\def\measure@#1\endalign{\global\lwidth@\z@\global\rwidth@\z@
 \global\maxlwidth@\z@\global\maxrwidth@\z@
 \global\and@\z@
 \setbox\z@\vbox
  {\noset@\everycr{\noalign{\global\tag@false\global\and@\z@}}\Let@
  \halign{\setboxz@h{$\m@th\displaystyle{\@lign##}$}%
   \global\lwidth@\wdz@
   \ifdim\lwidth@>\maxlwidth@\global\maxlwidth@\lwidth@\fi
   \global\advance\and@\@ne
   &\setboxz@h{$\m@th\displaystyle{{}\@lign##}$}\global\rwidth@\wdz@
   \ifdim\rwidth@>\maxrwidth@\global\maxrwidth@\rwidth@\fi
   \global\advance\and@\@ne
   &\Tag@\eat@{##}\crcr#1\crcr}}%
 \totwidth@\maxlwidth@\advance\totwidth@\maxrwidth@}
\def\prepost@{\global\let\tag@P@\tag@P\global\let\tag@Q@\tag@Q}
\def\reprepost@{\let\tag@P\tag@P@\let\tag@Q\tag@Q@}
\expandafter\def\csname align \space\endcsname#1\endalign
 {\measure@#1\endalign\global\and@\z@
 \ifingather@\everycr{\noalign{\global\and@\z@}}\else\displ@y@\fi
 \Let@\tabskip\centering@
 \halign to\displaywidth
  {\hfil\strut@\setboxz@h{$\m@th\displaystyle{\@lign##\prepost@}$}%
  \boxz@\global\advance\and@\@ne
  \tabskip\z@skip
  &\setboxz@h{$\m@th\displaystyle{{}\@lign##\prepost@}$}%
  \global\rwidth@\wdz@\boxz@\hfil\global\advance\and@\@ne
  \tabskip\centering@
  &\setboxz@h{\@lign\strut@\reprepost@\maketag@##\maketag@}%
  \dimen@\displaywidth\advance\dimen@-\totwidth@
  \divide\dimen@\tw@\advance\dimen@\maxrwidth@\advance\dimen@-\rwidth@
  \ifdim\dimen@<\tw@\wdz@\llap{\vtop{\normalbaselines\null\boxz@}}%
  \else\llap{\boxz@}\fi
  \tabskip\z@skip
  \crcr#1\crcr
  \black@\totwidth@}}
\expandafter\def\csname align \endcsname#1\endalign{\measure@#1\endalign
 \global\and@\z@
 \ifdim\totwidth@>\displaywidth\let\displaywidth@\totwidth@\else
  \let\displaywidth@\displaywidth\fi
 \ifingather@\everycr{\noalign{\global\and@\z@}}\else\displ@y@\fi
 \Let@\tabskip\centering@\halign to\displaywidth
  {\hfil\strut@\setboxz@h{$\m@th\displaystyle{\@lign##\prepost@}$}%
  \global\lwidth@\wdz@\global\lineht@\ht\z@
  \boxz@\global\advance\and@\@ne
  \tabskip\z@skip&\setboxz@h{$\m@th\displaystyle{{}\@lign##\prepost@}$}%
  \ifdim\ht\z@>\lineht@\global\lineht@\ht\z@\fi
  \boxz@\hfil\global\advance\and@\@ne
  \tabskip\centering@&\kern-\displaywidth@
  \setboxz@h{\@lign\strut@\reprepost@\maketag@##\maketag@}%
  \dimen@\displaywidth\advance\dimen@-\totwidth@
  \divide\dimen@\tw@\advance\dimen@\maxlwidth@\advance\dimen@-\lwidth@
  \ifdim\dimen@<\tw@\wdz@
   \rlap{\vbox{\normalbaselines\boxz@\vbox to\lineht@{}}}\else
   \rlap{\boxz@}\fi
  \tabskip\displaywidth@\crcr#1\crcr\black@\totwidth@}}
\def\attag@#1{\let\Maketag@\maketag@\let\TAG@\Tag@
 \let\Prepost@\prepost@\let\Reprepost@\reprepost@
 \let\Tag@\relax\let\maketag@\relax
 \let\prepost@\relax\let\reprepost@\relax
 \ifmeasuring@
  \def\llap@##1{\setboxz@h{##1}\hbox to\tw@\wdz@{}}%
  \def\rlap@##1{\setboxz@h{##1}\hbox to\tw@\wdz@{}}%
 \else\let\llap@\llap\let\rlap@\rlap\fi
 \toks@{\hfil\strut@
  $\m@th\displaystyle{\@lign\the\hashtoks@\prepost@}$%
  \tabskip\z@skip\global\advance\and@\@ne&
  $\m@th\displaystyle{{}\@lign\the\hashtoks@\prepost@}$\hfil
  \ifxat@\tabskip\centering@\fi\global\advance\and@\@ne}%
 \iftagsleft@
  \toks@@{\tabskip\centering@&\Tag@\kern-\displaywidth
   \rlap@{\@lign\reprepost@\maketag@\the\hashtoks@\maketag@}%
   \global\advance\and@\@ne\tabskip\displaywidth}\else
  \toks@@{\tabskip\centering@&\Tag@\llap@{\@lign\reprepost@\maketag@
   \the\hashtoks@\maketag@}\global\advance\and@\@ne\tabskip\z@skip}\fi
 \atcount@#1\relax\advance\atcount@\m@ne
 \loop\ifnum\atcount@>\z@
  \toks@\expandafter{\the\toks@&\hfil$\m@th\displaystyle{\@lign
  \the\hashtoks@\prepost@}$\global\advance\and@\@ne
  \tabskip\z@skip
  &$\m@th\displaystyle{{}\@lign\the\hashtoks@\prepost@}$\hfil\ifxat@
  \tabskip\centering@\fi\global\advance\and@\@ne}\advance\atcount@\m@ne
 \repeat
 \edef\preamble@{\the\toks@\the\toks@@}%
 \edef\preamble@@{\preamble@}%
 \let\maketag@\Maketag@\let\Tag@\TAG@
 \let\prepost@\Prepost@\let\reprepost@\Reprepost@}
\def\unlabel@{\def\label##1{\prevanish@\postvanish@}%
 \def\pagelabel##1{\prevanish@\postvanish@}}
\newcount\tag@CC
\expandafter\def\csname alignat \endcsname#1#2\endalignat
 {\inany@true\xat@false
 \def\tag{\global\tag@true
  \count@#1\relax\multiply\count@\tw@\advance\count@\m@ne
  \gdef\tag@{&}%
  \loop\ifnum\count@>\and@\xdef\tag@{&\omit\tag@}%
  \advance\count@\m@ne\repeat
  \tag@\relax}%
 \vspace@\allowdisplaybreak@\displaybreak@\intertext@
 \displ@y@\measuring@true\tag@CC\tag@C
 \setbox\savealignat@\hbox{\noset@\unlabel@$\m@th\displaystyle\Let@
  \attag@{#1}\vbox{\halign{\span\preamble@@\crcr#2\crcr}}$}%
 \measuring@false
 \Let@\attag@{#1}\tag@C\tag@CC
 \tabskip\centering@\halign to\displaywidth
  {\span\preamble@@\crcr#2\crcr\black@{\wd\savealignat@}}}
\expandafter\def\csname xalignat \endcsname#1#2\endxalignat
 {\inany@true\xat@true
 \def\tag{\global\tag@true
  \count@#1\relax\multiply\count@\tw@\advance\count@\m@ne
  \gdef\tag@{&}%
  \loop\ifnum\count@>\and@\xdef\tag@{&\omit\tag@}%
  \advance\count@\m@ne\repeat
  \tag@\relax}%
 \vspace@\allowdisplaybreak@\displaybreak@\intertext@
 \displ@y@\measuring@true\tag@CC\tag@C
 \setbox\savealignat@\hbox{\noset@\unlabel@$\m@th\displaystyle\Let@
  \attag@{#1}\vbox{\halign{\span\preamble@@\crcr#2\crcr}}$}%
 \measuring@false\Let@\attag@{#1}\tag@C\tag@CC
 \tabskip\centering@\halign to\displaywidth
 {\span\preamble@@\crcr#2\crcr\black@{\wd\savealignat@}}}
\def\gather{\RIfMIfI@\DN@{\onlydmatherr@\gather}\else
 \ingather@true\inany@true\def\tag{&\relax}%
 \vspace@\allowdisplaybreak@\displaybreak@\intertext@
 \displ@y\Let@
 \iftagsleft@\DN@{\csname gather \endcsname}\else
  \DN@{\csname gather \space\endcsname}\fi\fi
 \else\DN@{\onlydmatherr@\gather}\fi\next@}
\def\exstring@{\expandafter\eat@\string}
\def\newcounter#1{\define#1{}%
 \edef\next@{\def\noexpand#1{\futurelet\noexpand\next
  \csname\exstring@#1@Z\endcsname}}\next@
 \edef\next@{\def\csname\exstring@#1@Z\endcsname
  {\global\advance\csname\exstring@#1@C\endcsname\@ne
  {\csname\exstring@#1@F\endcsname\csname\exstring@#1@S\endcsname
   {\csname\exstring@#1@P\endcsname\csname\exstring@#1@N\endcsname
   {\noexpand\number\csname\exstring@#1@C\endcsname}%
   \csname\exstring@#1@Q\endcsname}}%
  \noexpand\ifx\noexpand\next\noexpand\label
   \def\noexpand\next@\noexpand\label########1{{\noexpand\noexpands@
    \xdef\noexpand\Thelabel@{\csname\exstring@#1@N\endcsname
     {\noexpand\number\csname\exstring@#1@C\endcsname}}%
    \xdef\noexpand\Thelabel@@@{\noexpand\number
     \csname\exstring@#1@C\endcsname}%
    \xdef\noexpand\Thelabel@@{\csname\exstring@#1@S\endcsname
     {\csname\exstring@#1@P\endcsname
     \csname\exstring@#1@N\endcsname
     {\noexpand\number\csname\exstring@#1@C\endcsname}%
     \csname\exstring@#1@Q\endcsname}}%
    \xdef\noexpand\Thelabel@@@@{\csname\exstring@#1@P\endcsname
     \csname\exstring@#1@N\endcsname
     {\noexpand\number\csname\exstring@#1@C\endcsname}%
     \csname\exstring@#1@Q\endcsname}}%
    {\noexpand\locallabel@\noexpand\label{########1}}}%
   \noexpand\else\let\noexpand\next@\relax\noexpand\fi\noexpand\next@}}\next@
 \expandafter\newcount@\csname\exstring@#1@C\endcsname
 \expandafter\let\csname\exstring@#1@N\endcsname\arabic
 \expandafter\def\csname\exstring@#1@S\endcsname##1{##1\/}%
 \expandafter\let\csname\exstring@#1@P\endcsname\empty
 \expandafter\let\csname\exstring@#1@Q\endcsname\empty
 \expandafter\def\csname\exstring@#1@F\endcsname{\rm}%
 }
\def\HASH@#1#2{\ifnum#2=\z@\else
 \edef\next@{\toks@{\the\toks@\the\hashtoks@#2}%
 \toks@@{\the\toks@@{\the\hashtoks@#2}}}\next@\expandafter\HASH@\fi}
\def\HASH@@{\toks@{}\toks@@{}\expandafter\HASH@\macpar@00}
\def\usecounter#1#2{\expandafter\ifx\csname\exstring@#1@Z\endcsname
 \relax\Err@{\noexpand#1not created with \string\newcounter}\fi
 \expandafter\let\csname\exstring@#1@@Z\endcsname\relax
 \expandafter\let\csname\exstring@#1@@Z@\endcsname\relax
 \expandafter\let\csname\exstring@#1@@Z@@\endcsname\relax
 \edef\next@{\def\noexpand#2{\futurelet\noexpand\next
  \csname\exstring@#1@@Z\endcsname}}\next@
 \edef\next@{\def\csname\exstring@#1@@Z\endcsname{\noexpand\ifx
  \noexpand\next\noexpand\label\def\noexpand\next@\noexpand\label
   ########1{\csname\exstring@#1@@Z@\endcsname
   {\noexpand#1\noexpand\label{########1}}}%
   \noexpand\else\noexpand\ifx\noexpand\next
   \noexpand"\def\noexpand\next@\noexpand"########1\noexpand"%
   {\csname\exstring@#1@@Z@\endcsname{{\expandafter\noexpand
   \csname\exstring@#1@F\endcsname
   \let\noexpand\pre\expandafter\noexpand\csname\exstring@#1@P\endcsname
   \let\noexpand\post\expandafter\noexpand\csname\exstring@#1@Q\endcsname
   \let\noexpand\style\expandafter\noexpand\csname\exstring@#1@S\endcsname
   \let\noexpand\numstyle\expandafter\noexpand\csname\exstring@#1@N\endcsname
   ########1}}}\noexpand\else
   \def\noexpand\next@{\csname\exstring@#1@@Z@\endcsname{\noexpand#1}}%
   \noexpand\fi\noexpand\fi\noexpand\next@}}\next@
 \def\next@{\expandafter\expandafter\expandafter\unmacro@\expandafter
  \meaning\csname\exstring@#1@@Z@@\endcsname\unmacro@
  \HASH@@
  \edef\next@{\def\csname\exstring@#1@@Z@\endcsname\the\toks@{%
   \expandafter\noexpand\csname\exstring@#1@@Z@@\endcsname\the\toks@@
   \noexpand\FNSSP@}}\next@}%
 \afterassignment\next@
 \expandafter\def\csname\exstring@#1@@Z@@\endcsname}
\def\listbi@{\penalty50 \medskip}
\def\listbii@{\penalty100 \smallskip}
\let\listbiii@\relax
\let\listbiv@\relax
\let\listbv@\relax
\def\listmi@{\advance\leftskip30\p@\relax}
\let\listmii@\listmi@
\let\listmiii@\listmi@
\let\listmiv@\listmi@
\let\listmv@\listmi@
\def\itemi@#1{\noindent@@\llap{#1\hskip5\p@}}
\let\itemii@\itemi@
\let\itemiii@\itemi@
\let\itemiv@\itemi@
\let\itemv@\itemi@
\def\liste@{\penalty-50 \medskip}
\def\listei@{\penalty-100 \smallskip}
\let\listeii@\relax
\let\listeiii@\relax
\let\listeiv@\relax
\expandafter\newcount\csname list@C1\endcsname
\csname list@C1\endcsname\z@
\expandafter\newcount\csname list@C2\endcsname
\csname list@C2\endcsname\z@
\expandafter\newcount\csname list@C3\endcsname
\csname list@C3\endcsname\z@
\expandafter\newcount\csname list@C4\endcsname
\csname list@C4\endcsname\z@
\expandafter\newcount\csname list@C5\endcsname
\csname list@C5\endcsname\z@
\expandafter\let\csname list@P1\endcsname\empty
\expandafter\let\csname list@P2\endcsname\empty
\expandafter\let\csname list@P3\endcsname\empty
\expandafter\let\csname list@P4\endcsname\empty
\expandafter\let\csname list@P5\endcsname\empty
\expandafter\let\csname list@Q1\endcsname\empty
\expandafter\let\csname list@Q2\endcsname\empty
\expandafter\let\csname list@Q3\endcsname\empty
\expandafter\let\csname list@Q4\endcsname\empty
\expandafter\let\csname list@Q5\endcsname\empty
\expandafter\def\csname list@S1\endcsname#1{{\rm(}{#1\/}{\rm)}}
\expandafter\def\csname list@S2\endcsname#1{{\rm(}{#1\/}{\rm)}}
\expandafter\def\csname list@S3\endcsname#1{{\rm(}{#1\/}{\rm)}}
\expandafter\def\csname list@S4\endcsname#1{{\rm(}{#1\/}{\rm)}}
\expandafter\def\csname list@S5\endcsname#1{{\rm(}{#1\/}{\rm)}}
\expandafter\let\csname list@N1\endcsname\arabic
\expandafter\let\csname list@N2\endcsname\arabic
\expandafter\let\csname list@N3\endcsname\arabic
\expandafter\let\csname list@N4\endcsname\arabic
\expandafter\let\csname list@N5\endcsname\arabic
\expandafter\def\csname list@F1\endcsname{\rm}
\expandafter\def\csname list@F2\endcsname{\rm}
\expandafter\def\csname list@F3\endcsname{\rm}
\expandafter\def\csname list@F4\endcsname{\rm}
\expandafter\def\csname list@F5\endcsname{\rm}
\newcount\listlevel@
\listlevel@\z@
\def\list@@C{\csname list@C\number\listlevel@\endcsname}
\def\list@@P{\csname list@P\number\listlevel@\endcsname}
\def\list@@Q{\csname list@Q\number\listlevel@\endcsname}
\def\list@@S{\csname list@S\number\listlevel@\endcsname}
\def\list@@N{\csname list@N\number\listlevel@\endcsname}
\def\list@@F{\csname list@F\number\listlevel@\endcsname}
\newif\iffirstitemi@
\newif\iffirstitemii@
\newif\iffirstitemiii@
\newif\iffirstitemiv@
\newif\iffirstitemv@
\def\Firstitem@true{\csname firstitem\romannumeral\listlevel@
 @true\endcsname}
\def\Firstitem@false{\csname firstitem\romannumeral\listlevel@
 @false\endcsname}
\def\Listm@{\csname listm\romannumeral\listlevel@ @\endcsname}
\def\Item@{\csname item\romannumeral\listlevel@ @\endcsname}
\def\Liste@{\csname liste\romannumeral\listlevel@ @\endcsname}
\newif\iflistcontinue@
\def\keepitem{\listcontinue@true}
\newcount\list@C@
\def\list{%
 \iflistcontinue@\csname list@C1\endcsname\csname list@C@\endcsname\fi
 \global\csname list@C2\endcsname\z@
 \global\csname list@C3\endcsname\z@
 \global\csname list@C4\endcsname\z@
 \global\csname list@C5\endcsname\z@
 \begingroup
 \firstitemi@true
 \listlevel@\@ne
 \def\item{\FN@\item@}%
 \FN@\list@}
\Invalid@\runinitem
\def\list@{\ifx\next\par
 \DN@\par{\FN@\list@}\else
 \ifx\next\runinitem
  \DN@\runinitem{\FN@\runinitem@}\else
  \DN@{\par\dimen@\parskip\parskip\dimen@}\fi\fi\next@}
\newif\ifoutlevel@
\newif\ifrunin@
\def\item@{%
 \ifoutlevel@\Liste@\outlevel@false\fi
 \ifrunin@\runin@false\par
  \dimen@\parskip\parskip\dimen@
  \Listm@\fi
 \iffirstitemi@\listbi@\listmi@\firstitemi@false\else\par\fi
 \iffirstitemii@\listbii@\listmii@\firstitemii@false\else\par\fi
 \iffirstitemiii@\listbiii@\listmiii@\firstitemiii@false\else\par\fi
 \iffirstitemiv@\listbiv@\listmiv@\firstitemiv@false\else\par\fi
 \iffirstitemv@\listbv@\listmv@\firstitemv@false\else\par\fi
 \DN@"##1"{{\let\pre\list@@P\let\post\list@@Q
  \let\style\list@@S\let\numstyle\list@@N
  \vskip-\parskip
  \Item@{\list@@F##1}%
  \noexpands@
  \Qlabel@{##1}}%
  \locallabel@
  \FNSSP@}%
 \DNii@{\global\advance\list@@C\@ne
  {\noexpands@
   \xdef\Thelabel@@@{\number\list@@C}%
   \xdefThelabel@\list@@N
   \xdef\Thelabel@@@@{\list@@P\Thelabel@\list@@Q}%
   \xdefThelabel@@\list@@S
  }%
  \locallabel@
  \vskip-\parskip
  \Item@{\list@@F\thelabel@@}%
  \FN@\pretendspace@}%
 \ifx\next"\expandafter\next@\else\expandafter\nextii@\fi}
\def\runinitem@{%
  \runin@true
  \Firstitem@false
  \DN@"##1"{{\let\pre\list@@P\let\post\list@@Q
   \let\style\list@@S\let\numstyle\list@@N
   \unskip\space{\list@@F##1} %
   \noexpands@
   \Qlabel@{##1}}%
   \locallabel@
   \ignorespaces}%
  \DNii@{\global\advance\list@@C\@ne
   {\noexpands@
    \xdef\Thelabel@@@{\number\list@@C}%
    \xdefThelabel@\list@@N
    \xdef\Thelabel@@@@{\list@@P\Thelabel@\list@@Q}%
    \xdefThelabel@@\list@@S
   }%
   \locallabel@
   \unskip\space{\list@@F\thelabel@@} }%
  \ifx\next"\expandafter\next@\else\expandafter\nextii@\fi}
\def\inlevel{\ifnum\listlevel@=5
 \DN@{\Err@{Already 5 levels down}}\else
 \DN@{\begingroup\advance\listlevel@\@ne
 \Firstitem@true\FN@\inlevel@}\fi\next@}
\def\inlevel@{\ifx\next\par
 \DN@\par{\FN@\inlevel@}\else
 \ifx\next\runinitem
  \DN@\runinitem{\FN@\runinitem@}\else
  \let\next@\relax\fi\fi\next@}
\def\outlevel{\ifnum\listlevel@=\@ne
 \Err@{At top level}\else
 \par\global\list@@C\z@\endgroup\outlevel@true\fi}
\def\endlist{%
 \expandafter\global\csname list@C@\endcsname\csname list@C1\endcsname
 \par
 \global\toks\@ne{}\count@\listlevel@
 {\loop
  \ifnum\count@>\z@\global\toks\@ne\expandafter{\the\toks\@ne\endgroup}%
  \advance\count@\m@ne
  \repeat}%
 \the\toks\@ne
 \liste@
 \listcontinue@false\global\csname list@C1\endcsname\z@
 \vskip-\parskip
 \noindent@@
 \FN@\pretendspace@}
\newif\iffirstdescribe@
\def\describe{\par
 \begingroup\firstdescribe@true
 \def\item##1{%
  \iffirstdescribe@\penalty50 \medskip\vskip-\parskip
  \firstdescribe@false\else\par\fi
  \noindent@@\hangindent2pc\hangafter\@ne
  {\bf##1}\hskip.5em}}

\Invalid@\pullin
\Invalid@\pullinmore
\newif\iffirstpull@
\def\margins{\par\begingroup\firstpull@true
 \def\pullin##1##2{\par
  \iffirstpull@\firstpull@false\else\endgroup\fi
  \begingroup\DN@{##1}%
  \ifx\next@\empty\leftskip\z@\else\ifx\next@\space\leftskip\z@
  \else\leftskip##1\fi\fi
  \DN@{##2}\ifx\next@\empty\rightskip\z@\else\ifx\next@\space
  \rightskip\z@\else\rightskip##2\fi\fi\ignorespaces}%
 \def\pullinmore##1##2{\par
  \xdef\Next@{\leftskip\the\leftskip\relax\rightskip\the\rightskip\relax}%
  \iffirstpull@\firstpull@false\else\endgroup\fi
  \begingroup\Next@
  \DN@{##1}%
  \ifx\next@\empty\else\ifx\next@\space\else\advance\leftskip##1\fi\fi
  \DN@{##2}\ifx\next@\empty\else\ifx\next@\space\else
  \advance\rightskip##2\fi\fi\ignorespaces}}

\newif\ifnopunct@
\newif\ifnospace@
\newif\ifoverlong@
\let\nofrillslist@\empty
\let\overlonglist@\empty
\def\nopunct{\nopunct@true\FN@\nopunct@}
\def\nospace{\nospace@true\FN@\nospace@}
\def\overlong{\overlong@true\FN@\overlong@}
\def\nopunct@{\ifx\next\nospace
 \DN@\nospace{\nospace@true\FN@\nopnos@}\else\ifx\next\overlong
 \DN@\overlong{\overlong@true\FN@\nopol@}\else
 \let\next@\nopunct@@\fi\fi\next@}
\def\nopunct@@#1{\ismember@\nofrillslist@#1%
 \iftest@\let\next@#1\else
 \DN@{\nopunct@false\Err@{\noexpand\nopunct can't be used with
 \string#1}#1}\fi\next@}
\def\nospace@{\ifx\next\nopunct
 \DN@\nopunct{\nopunct@true\FN@\nopnos@}\else\ifx\next\overlong
 \DN@\overlong{\overlong@true\FN@\nosol@}\else
 \let\next@\nospace@@\fi\fi\next@}
\def\nospace@@#1{\ismember@\nofrillslist@#1%
 \iftest@\let\next@#1\else
 \DN@{\nospace@false\Err@{\noexpand\nospace can't be used with
 \string#1}#1}\fi\next@}
\def\overlong@{\ifx\next\nopunct
 \DN@\nopunct{\nopunct@true\FN@\nopol@}\else\ifx\next\nospace
 \DN@\nospace{\nospace@true\FN@\nosol@}\else
 \let\next@\overlong@@\fi\fi\next@}
\def\overlong@@#1{\ismember@\overlonglist@#1%
 \iftest@\let\next@#1\else
 \DN@{\overlong@false\Err@{\noexpand\overlong can't be used with
 \string#1}#1}\fi\next@}
\def\nopnos@{\ifx\next\overlong
 \DN@\overlong{\overlong@true\nopnosol@}\else
 \let\next@\nopnos@@\fi\next@}
\def\nopol@{\ifx\next\nospace
 \DN@\nospace{\nospace@true\nopnosol@}\else
 \let\next@\nopol@@\fi\next@}
\def\nosol@{\ifx\next\nopunct
 \DN@\nopunct{\nopunct@true\nopnosol@}\else
 \let\next@\nosol@@\fi\next@}
\def\nopnos@@#1{\ismember@\nofrillslist@#1%
 \iftest@\let\next@#1\else
 \DN@{\nopunct@false\nospace@false
  \Err@{\noexpand\nopunct\noexpand\nospace
   can't be used with \string#1}#1}\fi\next@}
\def\testii@#1{\ismember@\nofrillslist@#1%
 \iftest@\let\nextiii@ T\else\let\nextiii@ F\fi
 \ismember@\overlonglist@#1%
 \iftest@\let\nextiv@ T\else\let\nextiv@ F\fi
 \test@false\if\nextiii@ T\if\nextiv@ T\test@true\fi\fi}
\def\nopol@@#1{\testii@{#1}%
 \iftest@\let\next@#1%
 \else\DN@{\if\nextiii@ T\else\nopunct@false\fi
  \if\nextiv@ T\else\overlong@false\fi
  \Err@{\if\nextiii@ T\else\noexpand\nopunct\fi
  \if\nextiv@ T\else\noexpand\overlong\fi can't be used
  with \string#1}#1}\fi\next@}
\def\nosol@@#1{\testii@{#1}%
 \iftest@\let\next@#1%
 \else\DN@{\if\nextiii@ T\else\nospace@false\fi
  \if\nextiv@ T\else\overlong@false\fi
  \Err@{\if\nextiii@ T\else\noexpand\nospace\fi
  \if\nextiv@ T\else\noexpand\overlong\fi can't be used
  with \string#1}#1}\fi\next@}
\def\nopnosol@#1{\testii@{#1}%
 \iftest@\let\next@#1%
 \else\DN@{\if\nextiii@ T\else\nopunct@false\nospace@false\fi
  \if\nextiv@ T\else\overlong@false\fi
  \Err@{\if\nextiii@ T\else\noexpand\nopunct\noexpand\nospace\fi
  \if\nextiv@ T\else\noexpand\overlong\fi can't be used
  with \string#1}#1}\fi\next@}
\def\punct@#1{\ifnopunct@\else#1\fi}
\def\addspace@#1{\ifnospace@\else#1\fi}
\def\hss@{\ifoverlong@\z@ plus\@m\p@ minus\@m\p@
 \else \z@ plus\@m\p@\fi}
\rightadd@\demo\to\nofrillslist@
\newif\ifclaim@
\def\exxx@{\expandafter\expandafter\expandafter\eat@\expandafter\string}
\let\colon@:
\def\demo#1{\ifclaim@
 \Err@{Previous \expandafter\noexpand\claimtype@ has
  no matching \string\end\exxx@\claimtype@}%
 \let\next@\relax
 \else
  \par
  \ifdim\lastskip<\smallskipamount\removelastskip\smallskip\fi
  \begingroup
  \noindent@@{\smc\ignorespaces#1\unskip
   \punct@{\null\colon@}\addspace@\enspace}%
  \nopunct@false\nospace@false
  \rm
  \DN@{\FNSSP@}%
 \fi
 \next@}
\def\enddemo{\par\endgroup\nopunct@false\nospace@false\smallskip}
\rightadd@\claim\to\nofrillslist@
\def\claim@F{\smc}
\def\claim@@@F{\csname\exxx@\claimtype@ @F\endcsname}
\def\claimformat@#1#2#3{%
 \medbreak\noindent@@{\smc#1 {\claim@@@F#2} #3%
 \punct@{\null.}\addspace@\enspace}\sl}
\def\claimformat@@#1#2{\claimformat@{\ignorespaces#1\unskip}%
 {\ifx\thelabel@@\empty\unskip\else\thelabel@@\fi}%
 {\ignorespaces#2\unskip}%
 \let\Claimformat@@\claimformat@@\FNSSP@}
\let\Claimformat@@\claimformat@@
\def\claim@@@P{\csname\exxx@\claimtype@ @P\endcsname}
\def\claim@@@Q{\csname\exxx@\claimtype@ @Q\endcsname}
\def\claim@@@S{\csname\exxx@\claimtype@ @S\endcsname}
\def\claim@@@N{\csname\exxx@\claimtype@ @N\endcsname}
\def\claim@@@C{\csname claim@C\claimclass@\endcsname}
\newcount\claim@C
\claim@C\z@
\let\claim@P\empty
\let\claim@Q\empty
\def\claim@S#1{#1\/}
\let\claim@N\arabic
\def\claim{\claim@true\let\claimclass@\empty
 \def\claimtype@{\claim}\FN@\claim@}
\def\claim@{%
 \ifx\next\c
  \let\next@\claim@c
 \else
  \ifx\next"%
   \let\next@\claim@q
  \else
   \begingroup\global\advance\claim@C\@ne
   {\noexpands@
    \xdef\Thelabel@@@{\number\claim@C}%
    \xdefThelabel@\claim@N
    \xdef\Thelabel@@@@{\claim@P\Thelabel@\claim@Q}%
    \xdefThelabel@@\claim@S
   }%
   \locallabel@
   \let\next@\Claimformat@@
  \fi
 \fi
 \next@}
\def\claim@c\c#1{\claim@true\begingroup
 \expandafter
 \ifx\csname claim@C#1\endcsname\relax
  \expandafter\newcount@\csname claim@C#1\endcsname
  \global\csname claim@C#1\endcsname\@ne
 \else
  \global\advance\csname claim@C#1\endcsname\@ne
 \fi
 \def\claimclass@{#1}%
 {\noexpands@
  \xdef\Thelabel@@@{\number\claim@@@C}%
  \xdefThelabel@\claim@@@N
  \xdef\Thelabel@@@@{\claim@@@P\Thelabel@\claim@@@Q}%
  \xdefThelabel@@\claim@@@S
 }%
 \locallabel@
 \FNSS@\claim@c@}
\def\claim@q"#1"{\begingroup
 {\let\pre\claim@@@P\let\post\claim@@@Q
  \let\style\claim@@@S\let\numstyle\claim@@@N
  \noexpands@
  \Qlabel@{#1}}%
 \locallabel@
 \FNSS@\claim@q@}
\def\claim@c@{\ifx\next"%
 \global\advance\claim@@@C\m@ne\let\next@\claim@cq
 \else\let\next@\Claimformat@@\fi\next@}
\def\claim@cq"#1"{{\let\pre\claim@@@P\let\post\claim@@@Q
 \let\style\claim@@@S\let\numstyle\claim@@@N
 \noexpands@
 \Qlabel@{#1}}%
 \locallabel@
 \FNSS@\Claimformat@@}
\def\claim@q@{\ifx\next\c\expandafter\claim@qc
 \else\expandafter\Claimformat@@\fi}
\def\claim@qc\c#1{\expandafter\ifx\csname claim@C#1\endcsname\relax
 \expandafter\newcount@\csname claim@C#1\endcsname
 \global\csname claim@C#1\endcsname\z@\fi
 \FNSS@\Claimformat@@}
\def\endclaim{\endgroup\claim@false\nopunct@false\nospace@false
 \let\Claimformat@@\claimformat@@\medbreak}
\Invalid@\claimclause
\def\newclaim{\FN@\newclaim@}
\def\newclaim@{\ifx\next\claimclause
 \DN@\claimclause##1{\newclaim@@{##1}}\else
 \DN@{\newclaim@@\relax}\fi\next@}
\def\claimlist@{\\\claim}
\newtoks\claim@i
\newtoks\claim@v
\let\noclaimclause@=F
\def\newclaim@@#1#2#3\c#4#5{\define#2{}%
 \rightadd@#2\to\claimlist@\rightadd@#2\to\nofrillslist@%
 \expandafter\def\csname\exstring@#2@P\endcsname{\claim@P}%
 \expandafter\def\csname\exstring@#2@Q\endcsname{\claim@Q}%
 \expandafter\def\csname\exstring@#2@S\endcsname{\claim@S}%
 \expandafter\def\csname\exstring@#2@N\endcsname{\claim@N}%
 \expandafter\def\csname\exstring@#2@F\endcsname{\claim@F}%
 \expandafter\def\csname end\exstring@#2\endcsname{\endclaim}%
 \expandafter\ifx\csname claim@C#4\endcsname\relax
  \expandafter\newcount@\csname claim@C#4\endcsname
  \global\csname claim@C#4\endcsname\z@\fi
 \edef\next@{\let\csname\exstring@#2@C\endcsname
   \csname claim@C#4\endcsname}\next@
 \def#2{\ifx\noclaimclause@ T\else#1\fi
  \global\claim@i{#1}\gdef\claim@iv{#4}\global\claim@v{#5}%
  \def\claimtype@{#2}\def\Claimformat@@{\claimformat@@{#5}}\claim@c\c{#4}}}
\def\shortenclaim#1#2{\define#2{}%
 \ismember@\claimlist@#1%
 \iftest@
  \rightadd@#2\to\nofrillslist@%
  \expandafter\def\csname\exstring@#2@P\endcsname
   {\csname\exstring@#1@P\endcsname}%
  \expandafter\def\csname\exstring@#2@Q\endcsname
   {\csname\exstring@#1@Q\endcsname}%
  \expandafter\def\csname\exstring@#2@S\endcsname
   {\csname\exstring@#1@S\endcsname}%
  \expandafter\def\csname\exstring@#2@N\endcsname
   {\csname\exstring@#1@N\endcsname}%
  \expandafter\def\csname\exstring@#2@F\endcsname
   {\csname\exstring@#1@F\endcsname}%
  \expandafter\def\csname end\exstring@#2\endcsname{\endclaim}%
  \edef\next@{\let\csname\exstring@#2@C\endcsname
    \csname claim\exstring@#1C\endcsname}\next@
  \setbox\z@\vbox{\let\noclaimclause@ T#1""\relax\endgroup}%
  \edef#2{\the\claim@i
   \def\noexpand\claimtype@{\noexpand#2}%
   \def\noexpand\Claimformat@@{\noexpand\claimformat@@{\the\claim@v}\relax}%
   \noexpand\claim@c\noexpand\c{\claim@iv}}%
 \else
  \Err@{\noexpand#1not yet created by \string\newclaim}%
 \fi}
\def\classtest@#1{\DN@{#1}\ifx\next@\claimclass@
 \test@true\else\test@false\fi}
\def\typetest@#1{\DN@{#1}\ifx\next@\claimtype@\test@true\else
  \test@false\fi}
\newif\iftoc@
\def\tocfile{\iftoc@\else\alloc@@7\write\chardef\sixt@@n\toc@
 \immediate\openout\toc@=\jobname.toc
 \alloc@@7\write\chardef\sixt@@n\tic@
 \immediate\openout\tic@=\jobname.tic
 \global\toc@true\fi}
\rightadd@\hl\to\nofrillslist@
\rightadd@\HL\to\overlonglist@
\def\HL@@C{\csname HL@C\HLlevel@\endcsname}
\def\HL@@P{\csname HL@P\HLlevel@\endcsname}
\def\HL@@Q{\csname HL@Q\HLlevel@\endcsname}
\def\HL@@S{\csname HL@S\HLlevel@\endcsname}
\def\HL@@N{\csname HL@N\HLlevel@\endcsname}
\def\HL@@F{\csname HL@F\HLlevel@\endcsname}
\def\HL@@@C{\csname\exxx@\HLtype@ @C\endcsname}
\def\HL@@@P{\csname\exxx@\HLtype@ @P\endcsname}
\def\HL@@@Q{\csname\exxx@\HLtype@ @Q\endcsname}
\def\HL@@@S{\csname\exxx@\HLtype@ @S\endcsname}
\def\HL@@@N{\csname\exxx@\HLtype@ @N\endcsname}
\def\HL#1{\expandafter
 \ifx\csname HL@C#1\endcsname\relax
  \DN@{\Err@{\string\HL#1 not defined in this style}}%
 \else
  \DN@{\gdef\HLlevel@{#1}\def\HLname@{\HL{#1}}\let\HLtype@\relax\FNSS@\HL@}%
 \fi
 \next@}%
\newif\ifquoted@
\let\aftertoc@\relax
\def\HL@{%
 \DN@"##1"##2\endHL{\def\entry@{##2}\quoted@true
  {\noexpands@
  \ifx\HLtype@\relax
   \let\pre\HL@@P\let\post\HL@@Q\let\style\HL@@S\let\numstyle\HL@@N
  \else
   \let\pre\HL@@@P\let\post\HL@@@Q\let\style\HL@@@S\let\numstyle\HL@@@N
  \fi
  \Qlabel@{##1}\let\style\relax\xdef\Qlabel@@@@{##1}%
  \xdef\Thepref@{\Thelabel@@@@}}%
  \csname HL@\HLlevel@\endcsname##2\endHL
  \let\pref\Thepref@
  \csname HL@I\HLlevel@\endcsname
  \csname HL@J\HLlevel@\endcsname
  \let\pref\pref@
  \HLtoc@
  \aftertoc@
  \let\aftertoc@\relax\overlong@false}%
 \DNii@##1\endHL{\def\entry@{##1}\quoted@false
  {\noexpands@
  \ifx\HLtype@\relax
   \global\advance\HL@@C\@ne
   \xdef\Thelabel@@@{\number\HL@@C}%
   \xdefThelabel@{\HL@@N}%
   \xdef\Thelabel@@@@{\HL@@P\Thelabel@\HL@@Q}%
   \xdefThelabel@@{\HL@@S}%
  \else
   \global\advance\HL@@@C\@ne
   \xdef\Thelabel@@@{\number\HL@@@C}%
   \xdefThelabel@{\HL@@@N}%
   \xdef\Thelabel@@@@{\HL@@@P\Thelabel@\HL@@@Q}%
   \xdefThelabel@@{\HL@@@S}%
  \fi
  \xdef\Thepref@{\Thelabel@@@@}}%
  \csname HL@\HLlevel@\endcsname##1\endHL
  \let\pref\Thepref@
  \csname HL@I\HLlevel@\endcsname
  \csname HL@J\HLlevel@\endcsname
  \let\pref\pref@
  \HLtoc@
  \aftertoc@
  \let\aftertoc@\relax\overlong@false}%
 \ifx\next"\expandafter\next@\else\expandafter\nextii@\fi}%
\Invalid@\endHL
\def\hl@@C{\csname hl@C\hllevel@\endcsname}
\def\hl@@P{\csname hl@P\hllevel@\endcsname}
\def\hl@@Q{\csname hl@Q\hllevel@\endcsname}
\def\hl@@S{\csname hl@S\hllevel@\endcsname}
\def\hl@@N{\csname hl@N\hllevel@\endcsname}
\def\hl@@F{\csname hl@F\hllevel@\endcsname}
\def\hl@@@C{\csname\exxx@\hltype@ @C\endcsname}
\def\hl@@@P{\csname\exxx@\hltype@ @P\endcsname}
\def\hl@@@Q{\csname\exxx@\hltype@ @Q\endcsname}
\def\hl@@@S{\csname\exxx@\hltype@ @S\endcsname}
\def\hl@@@N{\csname\exxx@\hltype@ @N\endcsname}
\def\hl#1{\expandafter
 \ifx\csname hl@C#1\endcsname\relax
  \DN@{\Err@{\string\hl#1 not defined in this style}}%
 \else
  \DN@{\gdef\hllevel@{#1}\def\hlname@{\hl{#1}}\let\hltype@\relax\FNSS@\hl@}%
 \fi
 \next@}
\def\hl@{%
 \DN@"##1"##2{\def\entry@{##2}\quoted@true
  {\noexpands@
  \ifx\hltype@\relax
   \let\pre\hl@@P\let\post\hl@@Q\let\style\hl@@S\let\numstyle\hl@@N
  \else
   \let\pre\hl@@@P\let\post\hl@@@Q\let\style\hl@@@S\let\numstyle\hl@@@N
  \fi
  \Qlabel@{##1}\let\style\relax\xdef\Qlabel@@@@{##1}%
  \xdef\Thepref@{\Thelabel@@@@}}%
  \csname hl@\hllevel@\endcsname{##2}%
  \let\pref\Thepref@
  \csname hl@I\hllevel@\endcsname
  \csname hl@J\hllevel@\endcsname
  \let\pref\pref@
  \hltoc@
  \aftertoc@
  \let\aftertoc@\relax\nopunct@false\nospace@false\FNSSP@}%
 \DNii@##1{\def\entry@{##1}\quoted@false
  {\noexpands@
  \ifx\hltype@\relax
   \global\advance\hl@@C\@ne
   \xdef\Thelabel@@@{\number\hl@@C}%
   \xdefThelabel@{\hl@@N}%
   \xdef\Thelabel@@@@{\hl@@P\Thelabel@\hl@@Q}%
   \xdefThelabel@@{\hl@@S}%
  \else
   \global\advance\hl@@@C\@ne
   \xdef\Thelabel@@@{\number\hl@@@C}%
   \xdefThelabel@{\hl@@@N}%
   \xdef\Thelabel@@@@{\hl@@@P\Thelabel@\hl@@@Q}%
   \xdefThelabel@@{\hl@@@S}%
  \fi
  \xdef\Thepref@{\Thelabel@@@@}}%
  \csname hl@\hllevel@\endcsname{##1}%
  \let\pref\Thepref@
  \csname hl@I\hllevel@\endcsname
  \csname hl@J\hllevel@\endcsname
  \let\pref\pref@
  \hltoc@
  \aftertoc@
  \let\aftertoc@\relax\nopunct@false\nospace@false\FNSSP@}%
 \ifx\next"\expandafter\next@\else\expandafter\nextii@\fi}%
\def\six@#1#2 #3 #4 #5 #6 #7 {\DN@{#2}\ifx\next@\empty
 \DN@##1\six@{}\else
 \write#1{ #2 #3 #4 #5 #6 #7}\DN@{\six@#1}\fi
 \next@}
\def\Sixtoc@{\ifx\macdef@\empty\else
 \DN@##1##2\next@{\def\macdef@{##1##2}}%
 \expandafter\next@\macdef@\next@
 \edef\next@
  {\noexpand\six@\toc@\macdef@
  \space\space\space\space\space\space\space\space\space\space\space\space
  \noexpand\six@}%
 \next@\let\macdef@\relax\fi}
\def\QorThelabel@@@@{\ifquoted@
 \noexpand\noexpand\noexpand"\Qlabel@@@@\noexpand\noexpand\noexpand"\else
 \Thelabel@@@@\fi}
\def\HLtoc@{%
 \iftoc@
 \expandafter\expandafter\expandafter\unmacro@
  \expandafter\meaning\csname HL@W\HLlevel@\endcsname\unmacro@
  {\noexpands@\let\style\relax
   \edef\next@{\write\toc@{\noexpand\noexpand\expandafter\noexpand\HLname@
   {\macdef@}{\QorThelabel@@@@}}}%
  \next@}%
  \expandafter\unmacro@\meaning\entry@\unmacro@
  \Sixtoc@
  \write\toc@{\noexpand\Page{\number\pageno}{\page@N}%
   {\page@P}{\page@Q}^^J}%
 \fi}
\def\hltoc@{%
 \iftoc@
 \expandafter\expandafter\expandafter\unmacro@
  \expandafter\meaning\csname hl@W\hllevel@\endcsname\unmacro@
  {\noexpands@\let\style\relax
  \edef\next@{\write\toc@{%
   \ifnopunct@\noexpand\noexpand\noexpand\nopunct\fi
   \ifnospace@\noexpand\noexpand\noexpand\nospace\fi
   \noexpand\noexpand\expandafter\noexpand\hlname@
   {\macdef@}{\QorThelabel@@@@}}}%
  \next@}%
  \expandafter\unmacro@\meaning\entry@\unmacro@
  \Sixtoc@
  \write\toc@{\noexpand\Page{\number\pageno}{\page@N}%
   {\page@P}{\page@Q}^^J}%
 \fi}
\def\mainfile#1{\def\mainfile@{#1}}
\def\checkmainfile@{\ifx\mainfile@\undefined
 \Err@{No \noexpand\mainfile specified}\fi}
\expandafter\newcount@\csname HL@C1\endcsname
\csname HL@C1\endcsname\z@
\expandafter\def\csname HL@S1\endcsname#1{#1\null.}
\expandafter\let\csname HL@N1\endcsname\arabic
\expandafter\let\csname HL@P1\endcsname\empty
\expandafter\let\csname HL@Q1\endcsname\empty
\expandafter\def\csname HL@F1\endcsname{\bf}
\expandafter\let\csname HL@W1\endcsname\empty
\expandafter\newcount@\csname hl@C1\endcsname
\csname hl@C1\endcsname\z@
\expandafter\def\csname hl@S1\endcsname#1{#1\/}
\expandafter\let\csname hl@N1\endcsname\arabic
\expandafter\let\csname hl@P1\endcsname\empty
\expandafter\let\csname hl@Q1\endcsname\empty
\expandafter\def\csname hl@F1\endcsname{\bf}
\expandafter\let\csname hl@W1\endcsname\empty
\expandafter\def\csname HL@1\endcsname#1\endHL{\bigbreak
 {\locallabel@
  \global\setbox\@ne\vbox{\Let@\tabskip\hss@
  \halign to\hsize{\bf\hfil\ignorespaces##\unskip\hfil\cr
  \expandafter\ifx\csname HL@W1\endcsname\empty\else
   \csname HL@W1\endcsname\space\fi
  {\HL@@F\ifx\thelabel@@\empty\else\thelabel@@\space\fi}%
  \ignorespaces#1\crcr}}%
  }%
 \unvbox\@ne\nobreak\medskip}
\expandafter\def\csname hl@1\endcsname#1{\medbreak\noindent@@
 {\locallabel@
 \bf{\hl@@F\ifx\thelabel@@\empty\else\thelabel@@\space\fi}%
 \ignorespaces#1\unskip\punct@{\null.}\addspace@\enspace}}
\expandafter\def\csname HL@I1\endcsname{\Reset\hl1{1}%
 \ifx\pref\empty\newpre\hl1{}\else\newpre\hl1{\pref.}\fi}
\def\NameHL#1#2{\define#2{}%
 \expandafter\ifx\csname HL@R#1\endcsname\relax
 \else
  \def\nextiv@{\let\nextiii@}%
  \expandafter\nextiv@\csname HL@R#1\endcsname
  \expandafter\let\nextiii@\undefined
  \expandafter\let\csname\exxx@\nextiii@ @C\endcsname\relax
  \expandafter\let\csname\exxx@\nextiii@ @P\endcsname\relax
  \expandafter\let\csname\exxx@\nextiii@ @Q\endcsname\relax
  \expandafter\let\csname\exxx@\nextiii@ @S\endcsname\relax
  \expandafter\let\csname\exxx@\nextiii@ @N\endcsname\relax
  \expandafter\let\csname\exxx@\nextiii@ @F\endcsname\relax
  \expandafter\let\csname\exxx@\nextiii@ @W\endcsname\relax
  \expandafter\let\csname end\exxx@\nextiii@\endcsname\undefined
 \fi
 \expandafter\gdef\csname HL@R#1\endcsname{#2}%
 \expandafter\gdef\csname\exstring@#2@R\endcsname{{HL}{#1}}%
 \iftoc@\write\toc@{\noexpand\NameHL#1\noexpand#2^^J}\fi
 \rightadd@#2\to\overlonglist@
 \edef\next@{\let\csname\exstring@#2@C\endcsname\expandafter\noexpand
  \csname HL@C#1\endcsname}\next@
 \edef\next@{\let\csname\exstring@#2@P\endcsname\expandafter\noexpand
  \csname HL@P#1\endcsname}\next@
 \edef\next@{\let\csname\exstring@#2@Q\endcsname\expandafter\noexpand
  \csname HL@Q#1\endcsname}\next@
 \edef\next@{\let\csname\exstring@#2@S\endcsname\expandafter\noexpand
  \csname HL@S#1\endcsname}\next@
 \edef\next@{\let\csname\exstring@#2@N\endcsname\expandafter\noexpand
  \csname HL@N#1\endcsname}\next@
 \edef\next@{\let\csname\exstring@#2@F\endcsname\expandafter\noexpand
  \csname HL@F#1\endcsname}\next@
 \edef\next@{\let\csname\exstring@#2@W\endcsname\expandafter\noexpand
  \csname HL@W#1\endcsname}\next@
 \edef\next@{\def\noexpand#2####1\expandafter\noexpand
  \csname end\exstring@#2\endcsname
  {\def\noexpand\HLtype@{\noexpand#2}%
   \def\noexpand\HLname@{\noexpand#2}%
   \gdef\noexpand\HLlevel@{#1}%
   \noexpand\FNSS@\noexpand\HL@####1\noexpand\endHL}}%
  \next@
 \edef\next@{\noexpand\Invalid@\expandafter\noexpand
  \csname end\exstring@#2\endcsname}%
 \next@}
\def\Namehl#1#2{\define#2{}%
 \expandafter\ifx\csname hl@R#1\endcsname\relax
 \else
  \def\nextiv@{\let\nextiii@}%
  \expandafter\nextiv@\csname hl@R#1\endcsname
  \expandafter\let\nextiii@\undefined
  \expandafter\let\csname\exxx@\nextiii@ @C\endcsname\relax
  \expandafter\let\csname\exxx@\nextiii@ @P\endcsname\relax
  \expandafter\let\csname\exxx@\nextiii@ @Q\endcsname\relax
  \expandafter\let\csname\exxx@\nextiii@ @S\endcsname\relax
  \expandafter\let\csname\exxx@\nextiii@ @N\endcsname\relax
  \expandafter\let\csname\exxx@\nextiii@ @F\endcsname\relax
  \expandafter\let\csname\exxx@\nextiii@ @W\endcsname\relax
 \fi
 \expandafter\gdef\csname hl@R#1\endcsname{#2}%
 \expandafter\gdef\csname\exstring@#2@R\endcsname{{hl}{#1}}%
 \iftoc@\write\toc@{\noexpand\Namehl#1\noexpand#2^^J}\fi
 \rightadd@#2\to\nofrillslist@%
 \edef\next@{\let\csname\exstring@#2@C\endcsname\expandafter\noexpand
  \csname hl@C#1\endcsname}\next@
 \edef\next@{\let\csname\exstring@#2@P\endcsname\expandafter\noexpand
  \csname hl@P#1\endcsname}\next@
 \edef\next@{\let\csname\exstring@#2@Q\endcsname\expandafter\noexpand
  \csname hl@Q#1\endcsname}\next@
 \edef\next@{\let\csname\exstring@#2@S\endcsname\expandafter\noexpand
  \csname hl@S#1\endcsname}\next@
 \edef\next@{\let\csname\exstring@#2@N\endcsname\expandafter\noexpand
  \csname hl@N#1\endcsname}\next@
 \edef\next@{\let\csname\exstring@#2@F\endcsname\expandafter\noexpand
  \csname hl@F#1\endcsname}\next@
 \edef\next@{\let\csname\exstring@#2@W\endcsname\expandafter\noexpand
  \csname hl@W#1\endcsname}\next@
 \edef\next@{\def\noexpand#2{%
  \def\noexpand\hltype@{\noexpand#2}%
  \def\noexpand\hlname@{\noexpand#2}%
  \gdef\noexpand\hllevel@{#1}%
  \noexpand\FNSS@\noexpand\hl@}}%
 \next@}%
\def\Initialize{\FN@\Init@}
\def\Init@{\ifx\next\HL\let\next@\InitH@\else\ifx\next\hl\let\next@\InitH@
  \else\let\next@\InitS@\fi\fi\next@}
\def\InitH@#1#2{\expandafter\ifx\csname\exstring@#1@C#2\endcsname\relax
 \DN@{\Err@{\noexpand#1level #2 not defined in this style}}\else
 \DN@{\expandafter\gdef\csname\exstring@#1@J#2\endcsname}\fi\next@}
\def\InitC@#1#2{\edef\nextii@{\expandafter\noexpand\csname#1\endcsname{#2}}}
\def\InitS@#1{\expandafter\ifx\csname\exstring@#1@R\endcsname\relax
 \Err@{\noexpand#1not defined in this style}\let\next@\relax\else
 \DN@{\let\next@}\expandafter\next@\csname\exstring@#1@R\endcsname
 \expandafter\InitC@\next@
 \DN@{\expandafter\InitH@\nextii@}\fi\next@}
\def\value#1{\expandafter
 \ifx\csname\exstring@#1@C\endcsname\relax
  \expandafter\ifx\csname\exstring@#1@C1\endcsname\relax
   \DN@{\Err@{\noexpand\value can't be used with \string#1}}%
  \else
   \DN@{\value@#1}%
  \fi
 \else
  \DN@{\number\csname\exstring@#1@C\endcsname\relax}%
 \fi
 \next@}
\def\value@#1#2{\expandafter
 \ifx\csname\exstring@#1@C#2\endcsname\relax
  \DN@{\Err@{\string\value\string#1 can't be followed by \string#2}}%
 \else
  \DN@{\number\csname\exstring@#1@C#2\endcsname\relax}%
 \fi
 \next@}
\newcount\Value
\def\Evaluate#1{\expandafter
 \ifx\csname\exstring@#1@C\endcsname\relax
  \expandafter\ifx\csname\exstring@#1@C1\endcsname\relax
   \DN@{\Err@{\noexpand\Evaluate can't be used with \string#1}}%
  \else
   \DN@{\Evaluate@#1}%
  \fi
 \else
  \DN@{\global\Value\csname\exstring@#1@C\endcsname}%
 \fi
 \next@}
\def\Evaluate@#1#2{\expandafter
 \ifx\csname\exstring@#1@C#2\endcsname\relax
  \DN@{\Err@{\string\Evaluate\string#1 can't be followed by \string#2}}%
 \else
  \DN@{\global\Value\csname\exstring@#1@C#2\endcsname}%
 \fi\next@}
\def\pre#1{\expandafter
 \ifx\csname\exstring@#1@P\endcsname\relax
  \expandafter\ifx\csname\exstring@#1@P1\endcsname\relax
   \DN@{\Err@{\noexpand\pre can't be used with \string#1}}%
  \else
   \DN@{\pre@#1}%
  \fi
 \else
  \DN@{{\csname\exstring@#1@P\endcsname}}%
 \fi
 \next@}
\def\pre@#1#2{\expandafter
 \ifx\csname\exstring@#1@P#2\endcsname\relax
  \DN@{\Err@{\string\pre\string#1 can't be followed by \string#2}}%
 \else
  \DN@{{\csname\exstring@#1@P#2\endcsname}}%
 \fi
 \next@}
\def\post#1{\expandafter
 \ifx\csname\exstring@#1@Q\endcsname\relax
  \expandafter\ifx\csname\exstring@#1@Q1\endcsname\relax
   \DN@{\Err@{\noexpand\post can't be used with \string#1}}%
  \else
   \DN@{\post@#1}%
  \fi
 \else
  \DN@{{\csname\exstring@#1@Q\endcsname}}%
 \fi
 \next@}
\def\post@#1#2{\expandafter
 \ifx\csname\exstring@#1@Q#2\endcsname\relax
  \DN@{\Err@{\string\post\string#1 can't be followed by \string#2}}%
 \else
  \DN@{{\csname\exstring@#1@Q#2\endcsname}}%
 \fi
 \next@}
\def\style#1{\expandafter
 \ifx\csname\exstring@#1@S\endcsname\relax
  \expandafter\ifx\csname\exstring@#1@S1\endcsname\relax
   \DN@{\Err@{\noexpand\style can't be used with \string#1}}%
  \else
   \DN@{\style@#1}%
  \fi
 \else
  \DN@{\csname\exstring@#1@S\endcsname}%
 \fi
 \next@}
\def\style@#1#2{\expandafter
 \ifx\csname\exstring@#1@S#2\endcsname\relax
  \DN@{\Err@{\string\style\string#1 can't be followed by \string#2}}%
 \else
  \DN@{\csname\exstring@#1@S#2\endcsname}%
 \fi
 \next@}
\def\fontstyle#1{\expandafter
 \ifx\csname\exstring@#1@F\endcsname\relax
  \expandafter\ifx\csname\exstring@#1@F1\endcsname\relax
   \DN@{\Err@{\noexpand\fontstyle can't be used with \string#1}}%
  \else
   \DN@{\fontstyle@#1}%
  \fi
 \else
  \DN@##1{{\csname\exstring@#1@F\endcsname##1}}%
 \fi
 \next@}
\def\fontstyle@#1#2{\expandafter
 \ifx\csname\exstring@#1@F#2\endcsname\relax
  \DN@{\Err@{\string\fontstyle\string#1 can't be followed by \string#2}}%
 \else
  \DN@##1{{\csname\exstring@#1@F#2\endcsname##1}}%
 \fi
 \next@}
\def\Reset#1{\expandafter
 \ifx\csname\exstring@#1@C\endcsname\relax
  \expandafter\ifx\csname\exstring@#1@C1\endcsname\relax
   \DN@{\Err@{\noexpand\Reset can't be used with \string#1}}%
  \else
   \DN@{\Reset@#1}%
  \fi
 \else
  \DN@##1{\count@##1\relax\ifx#1\page\else\advance\count@\m@ne\fi
   \global\csname\exstring@#1@C\endcsname\count@}%
 \fi
 \next@}
\def\Reset@#1#2{\expandafter
 \ifx\csname\exstring@#1@C#2\endcsname\relax
  \DN@{\Err@{\string\Reset\string#1 can't be followed by \string#2}}%
 \else
  \DN@##1{\count@##1\relax\advance\count@\m@ne
   \global\csname\exstring@#1@C#2\endcsname\count@}%
 \fi
 \next@}
\def\Offset#1{\expandafter
 \ifx\csname\exstring@#1@C\endcsname\relax
  \expandafter\ifx\csname\exstring@#1@C1\endcsname\relax
   \DN@{\Err@{\noexpand\Offset can't be used with \string#1}}%
  \else
   \DN@{\Offset@#1}%
  \fi
 \else
  \DN@##1{\count@##1\relax\advance\count@\m@ne\global\advance
   \csname\exstring@#1@C\endcsname\count@}%
 \fi
 \next@}
\def\Offset@#1#2{\expandafter
 \ifx\csname\exstring@#1@C#2\endcsname\relax
  \DN@{\Err@{\string\Offset\string#1 can't be followed by \string#2}}%
 \else
  \DN@##1{\count@##1\relax\advance\count@\m@ne
   \global\advance\csname\exstring@#1@C#2\endcsname\count@}%
 \fi
 \next@}
\def\getR@#1#2{\def\nextiv@{\let\nextiii@}\expandafter\nextiv@
 \csname\exstring@#1@R#2\endcsname}
\def\letR@#1#2#3{\expandafter\let\csname#1@#3#2\endcsname\Next@}
\def\letR@@#1#2{\expandafter\let\csname\exstring@#1@#2\endcsname\Next@}
\def\newpre#1{\expandafter
 \ifx\csname\exstring@#1@P\endcsname\relax
  \expandafter\ifx\csname\exstring@#1@P1\endcsname\relax
   \DN@{\Err@{\noexpand\newpre can't be used with \string#1}}%
  \else
   \DN@{\newpre@#1}%
  \fi
 \else
  \DN@{%
   \DNii@{%
    \endgroup
    \expandafter\let\csname\exstring@#1@P\endcsname\Next@
    \expandafter\ifx\csname\exstring@#1@R\endcsname\relax\else
    \getR@#1{}\expandafter\letR@\nextiii@ P\fi
    }%
   \begingroup\noexpands@\afterassignment\nextii@\xdef\Next@}%
 \fi
 \next@}
\def\newpre@#1#2{\expandafter
 \ifx\csname\exstring@#1@P#2\endcsname\relax
  \DN@{\Err@{\string\newpre\string#1 can't be followed by \string#2}}%
 \else
  \DN@{%
   \DNii@{%
    \endgroup
    \expandafter\let\csname\exstring@#1@P#2\endcsname\Next@
    \expandafter\ifx\csname\exstring@#1@R#2\endcsname\relax\else
    \getR@#1{#2}\expandafter\letR@@\nextiii@ P\fi
    }%
   \begingroup\noexpands@\afterassignment\nextii@\xdef\Next@}%
 \fi
 \next@}
\def\newpost#1{\expandafter
 \ifx\csname\exstring@#1@Q\endcsname\relax
  \expandafter\ifx\csname\exstring@#1@Q1\endcsname\relax
   \DN@{\Err@{\noexpand\newpost can't be used with \string#1}}%
  \else
   \DN@{\newpost@#1}%
  \fi
 \else
  \DN@{%
   \DNii@{%
    \endgroup
    \expandafter\let\csname\exstring@#1@Q\endcsname\Next@
    \expandafter\ifx\csname\exstring@#1@R\endcsname\relax\else
    \getR@#1{}\expandafter\letR@\nextiii@ Q\fi
    }%
   \begingroup\noexpands@\afterassignment\nextii@\xdef\Next@}%
 \fi
 \next@}
\def\newpost@#1#2{\expandafter
 \ifx\csname\exstring@#1@Q#2\endcsname\relax
  \DN@{\Err@{\string\newpost\string#1 can't be followed by \string#2}}%
 \else
  \DN@{%
   \DNii@{%
    \endgroup
    \expandafter\let\csname\exstring@#1@Q#2\endcsname\Next@
    \expandafter\ifx\csname\exstring@#1@R#2\endcsname\relax\else
    \getR@#1{#2}\expandafter\letR@@\nextiii@ Q\fi
    }%
   \begingroup\noexpands@\afterassignment\nextii@\xdef\Next@}%
 \fi
 \next@}
\def\newstyle#1{\expandafter
 \ifx\csname\exstring@#1@S\endcsname\relax
  \expandafter\ifx\csname\exstring@#1@S1\endcsname\relax
   \DN@{\Err@{\noexpand\newstyle can't be used
    with \string#1}}%
  \else
   \DN@{\newstyle@#1}%
  \fi
 \else
  \DN@{%
   \DNii@{%
    \expandafter\let\csname\exstring@#1@S\endcsname\Next@
    \expandafter\ifx\csname\exstring@#1@R\endcsname\relax\else
    \getR@#1{}\expandafter\letR@\nextiii@ S\fi
    }%
   \afterassignment\nextii@\gdef\Next@}%
 \fi
 \next@}
\def\newstyle@#1#2{\expandafter
 \ifx\csname\exstring@#1@S#2\endcsname\relax
  \DN@{\Err@{\string\newstyle\string#1 can't be followed by
   \string#2}}%
 \else
  \DN@{%
   \DNii@{%
    \expandafter\let\csname\exstring@#1@S#2\endcsname\Next@
    \expandafter\ifx\csname\exstring@#1@R#2\endcsname\relax\else
    \getR@#1{#2}\expandafter\letR@@\nextiii@ S\fi
    }%
   \afterassignment\nextii@\gdef\Next@}%
 \fi
 \next@}
\def\newnumstyle#1{\expandafter
 \ifx\csname\exstring@#1@N\endcsname\relax
  \expandafter\ifx\csname\exstring@#1@N1\endcsname\relax
   \DN@{\Err@{\noexpand\newnumstyle can't be used with
    \string#1}}%
  \else
   \DN@{\newnumstyle@#1}%
  \fi
 \else
  \DN@##1{%
   \gdef\Next@{##1}%
    \expandafter\let\csname\exstring@#1@N\endcsname\Next@
    \expandafter\ifx\csname\exstring@#1@R\endcsname\relax\else
    \getR@#1{}\expandafter\letR@\nextiii@ N\fi
    }%
 \fi
 \next@}
\def\newnumstyle@#1#2{\expandafter
 \ifx\csname\exstring@#1@N#2\endcsname\relax
  \DN@{\Err@{\string\newnumstyle\string#1 can't be followed by
   \string#2}}%
 \else
  \DN@##1{%
   \gdef\Next@{##1}%
    \expandafter\let\csname\exstring@#1@N#2\endcsname\Next@
    \expandafter\ifx\csname\exstring@#1@R#2\endcsname\relax\else
    \getR@#1{#2}\expandafter\letR@@\nextiii@ N\fi
    }%
  \fi
 \next@}
\def\newfontstyle#1{\expandafter
 \ifx\csname\exstring@#1@F\endcsname\relax
  \expandafter\ifx\csname\exstring@#1@F1\endcsname\relax
   \DN@{\Err@{\noexpand\newfontstyle can't be used with
    \string#1}}%
  \else
   \DN@{\newfontstyle@#1}%
  \fi
 \else
  \DN@##1{%
   \gdef\Next@{##1}%
    \expandafter\let\csname\exstring@#1@F\endcsname\Next@
    \expandafter\ifx\csname\exstring@#1@R\endcsname\relax\else
    \getR@#1{}\expandafter\letR@\nextiii@ F\fi
    }%
 \fi
 \next@}
\def\newfontstyle@#1#2{\expandafter
 \ifx\csname\exstring@#1@F#2\endcsname\relax
  \DN@{\Err@{\string\newfontstyle\string#1 can't be followed by
   \string#2}}%
 \else
  \DN@##1{%
   \gdef\Next@{##1}%
    \expandafter\let\csname\exstring@#1@F#2\endcsname\Next@
    \expandafter\ifx\csname\exstring@#1@R#2\endcsname\relax\else
    \getR@#1{#2}\expandafter\letR@@\nextiii@ F\fi
    }%
 \fi
 \next@}
\def\word#1{\expandafter
 \ifx\csname\exstring@#1@W\endcsname\relax
  \expandafter\ifx\csname\exstring@#1@W1\endcsname\relax
   \DN@{\Err@{\noexpand\word can't be used with \string#1}}%
  \else
   \DN@{\word@#1}%
  \fi
 \else
  \DN@{{\csname\exstring@#1@W\endcsname}}%
 \fi
 \next@}
\def\word@#1#2{\expandafter
 \ifx\csname\exstring@#1@W#2\endcsname\relax
  \DN@{\Err@{\string\word\noexpand#1can't be followed by \string#2}}%
 \else
  \DN@{{\csname\exstring@#1@W#2\endcsname}}%
 \fi
 \next@}
\def\newword#1{\expandafter
 \ifx\csname\exstring@#1@W\endcsname\relax
  \expandafter\ifx\csname\exstring@#1@W1\endcsname\relax
   \DN@{\Err@{\noexpand\newword can't be used  with \string#1}}%
  \else
   \DN@{\newword@#1}%
  \fi
 \else
  \DN@{%
   \DNii@{%
    \expandafter\let\csname\exstring@#1@W\endcsname\Next@
    \expandafter\ifx\csname\exstring@#1@R\endcsname\relax\else
     \getR@#1{}\expandafter\letR@\nextiii@ W\fi
    }%
   \afterassignment\nextii@\gdef\Next@}%
 \fi
 \next@}
\def\newword@#1#2{\expandafter
 \ifx\csname\exstring@#1@W#2\endcsname\relax
  \DN@{\Err@{\string\newword\noexpand#1can't be followed by \string#2}}%
 \else
  \DN@{%
   \DNii@{%
    \expandafter\let\csname\exstring@#1@W#2\endcsname\Next@
    \expandafter\ifx\csname\exstring@#1@R#2\endcsname\relax\else
     \getR@#1{#2}\expandafter\letR@@\nextiii@ W\fi
    }%
   \afterassignment\nextii@\gdef\Next@}%
 \fi
 \next@}
\newif\iffn@
\newcount\footmark@C
\footmark@C\z@
\def\footmark@S#1{$^{#1}$}
\let\footmark@N\arabic
\def\footmark@F{\rm}
\def\foottext@S#1{$^{#1}$}
\def\foottext@F{\rm}
\let\modifyfootnote@\relax
\def\modifyfootnote#1{\def\modifyfootnote@{#1}}
\def\vfootnote@#1{\insert\footins
 \bgroup
 \floatingpenalty\@MM\interlinepenalty\interfootnotelinepenalty
 \leftskip\z@\rightskip\z@\spaceskip\z@\xspaceskip\z@
 \rm\splittopskip\ht\strutbox\splitmaxdepth\dp\strutbox
 \locallabel@\noindent@@{\foottext@F#1}\modifyfootnote@
 \footstrut\FN@\fo@t}
\def\fo@t{\ifcat\bgroup\noexpand\next\expandafter\f@@t\else
 \expandafter\f@t\fi}
\def\f@t#1{#1\@foot}
\def\f@@t{\bgroup\aftergroup\@foot\afterassignment\FNSSP@\let\next@}
\def\@foot{\unskip\lower\dp\strutbox\vbox to\dp\strutbox{}\egroup
 \iffn@\expandafter\fn@false\else
 \expandafter\postvanish@\fi}
\newif\ifplainfn@
\plainfn@true
\def\fancyfootnotes{\plainfn@false}
\newcount\fancyfootmarkcount@
\fancyfootmarkcount@\z@
\newcount\lastfnpage@
\lastfnpage@-\@M
\let\justfootmarklist@\empty
\def\footmark{\let\@sf\empty
 \ifhmode\edef\@sf{\spacefactor\the\spacefactor}\/\fi
 \DN@{\ifx"\next\expandafter\nextii@\else\expandafter\footmark@\fi}%
 \DNii@"##1"{%
  \iffirstchoice@
   {\let\style\footmark@S\let\numstyle\footmark@N
   \footmark@F##1%
   \noexpands@
   \let\style\foottext@S
   \Qlabel@{##1}%
   }%
   \iffn@\else
    {\noexpands@
    \xdef\Next@{{\Thelabel@}{\Thelabel@@}{\Thelabel@@@}{\Thelabel@@@@}}%
    }%
    \expandafter\rightappend@\Next@\to\justfootmarklist@
   \fi
  \fi
  \@sf\relax}%
 \FN@\next@}
\def\footmark@{%
 \iffirstchoice@
  \global\advance\footmark@C\@ne
  \ifplainfn@
   \xdef\adjustedfootmark@{\number\footmark@C}%
  \else
   {\let\\\or\xdef\Next@{\ifcase\number\footmark@C\fnpages@\else
     -\@M\fi}}%
   \ifnum\Next@=-\@M
    \xdef\adjustedfootmark@{\number\footmark@C}%
   \else
    \ifnum\Next@=\lastfnpage@
     \global\advance\fancyfootmarkcount@\@ne
    \else
     \global\fancyfootmarkcount@\@ne
     \global\lastfnpage@\Next@
    \fi
    \xdef\adjustedfootmark@{\number\fancyfootmarkcount@}%
   \fi
  \fi
  {\noexpands@
  \xdef\Thelabel@@@{\adjustedfootmark@}%
  \xdefThelabel@\footmark@N
  \xdef\Thelabel@@@@{\Thelabel@}%
  \xdefThelabel@@\foottext@S
  }%
  \iffn@\else
   {\noexpands@
   \xdef\Next@{{\Thelabel@}{\Thelabel@@}{\Thelabel@@@}{\Thelabel@@@@}}%
   }%
   \expandafter\rightappend@\Next@\to\justfootmarklist@
  \fi
  \ifplainfn@
  \else
   \edef\next@{\write\laxwrite@{F\noexpand\the\pageno}}\next@
  \fi
 \fi
 \footmark@S{\footmark@N{\adjustedfootmark@}}%
 \@sf\relax}
\def\foottext{\prevanish@
 \ifx\justfootmarklist@\empty
  \Err@{There is no \noexpand\footmark for this \string\foottext}\fi
 \DN@\\##1##2\next@{\DN@{##1}\gdef\justfootmarklist@{##2}}%
 \expandafter\next@\justfootmarklist@\next@
 \expandafter\foottext@\next@}
\def\foottext@#1#2#3#4{{\noexpands@
  \xdef\Thelabel@{#1}\xdef\Thelabel@@{#2}%
  \xdef\Thelabel@@@{#3}\xdef\Thelabel@@@@{#4}}%
  \vfootnote@{\thelabel@@}}
\rightadd@\foottext\to\vanishlist@
\def\footnote{\fn@true
 \let\@sf\empty
 \ifhmode\edef\@sf{\spacefactor\the\spacefactor}\/\fi
 \DN@{\ifx"\next\expandafter\nextii@\else\expandafter\nextiii@\fi}%
 \DNii@"##1"{\footmark"##1"\vfootnote@{\let\style\foottext@S
  \let\numstyle\footmark@N##1}}%
 \def\nextiii@{\footmark\vfootnote@{\foottext@S{\footmark@N
  {\adjustedfootmark@}}}}%
 \FN@\next@}
\newdimen\litindent
\litindent20\p@
\newbox\litbox@
\newbox\Litbox@
\newcount\interlitpenalty@
\interlitpenalty@\@M
\newcount\litlines@
{\obeyspaces\gdef\defspace@{\def {\allowbreak\hskip.5emminus.15em}}}
{\obeylines\gdef\letM@{\let^^M\CtrlM@}}
\def\CtrlM@{\egroup
 \ifcase\litlines@\advance\litlines@\@ne\or
 \box\litbox@\advance\litlines@\@ne\else
 \penalty\interlitpenalty@\box\litbox@\fi
 \Lit@}
\def\Lit@{\setbox\litbox@\hbox\bgroup\litdefs@\hskip\litindent}
\newcount\littab@
\littab@8
\def\littab#1{\littab@#1\relax}
{\catcode`\^^I=\active\gdef\letTAB@{\let^^I\TAB@}}
\def\TAB@{\egroup
 \dimen@\wd\litbox@
 \advance\dimen@-\litindent
 \setboxz@h{\tt0}%
 \dimen@ii\littab@\wdz@
 \divide\dimen@\dimen@ii
 \multiply\dimen@\dimen@ii
 \advance\dimen@\littab@\wdz@
 \advance\dimen@\litindent
 \setbox\litbox@\hbox\bgroup\litdefs@\hbox to\dimen@{\unhbox\litbox@\hfil}}
{\catcode`\`=\active\gdef`{\relax\lq}}
\let\litbs@\relax
\let\litbs@@\relax
\def\litbackslash#1{%
 \edef\litbs@{\catcode`\string#1=\z@
 \def\noexpand\litbs@@{\def\expandafter\noexpand\csname\string#1\endcsname
  {\char`\string#1}}}}
\def\litcodes@{\catcode`\\=12
 \catcode`\{=12 \catcode`\}=12
 \catcode`\$=12 \catcode`\&=12
 \catcode`\#=12
 \catcode`\^=12 \catcode`\_=12
 \catcode`\@=12 \catcode`\~=12 \catcode`\"=12
 \catcode`\;=12 \catcode`\:=12 \catcode`\!=12 \catcode`\?=12
 \catcode`\%=12 \litbs@\catcode`\`=\active\obeyspaces\defspace@}
\def\activate@#1#2{{\lccode`\~=`#2%
 \lowercase{%
  \if0#1%
  \gdef\Next@{\def~{\egroup\endgroup\bigskip\vskip-\parskip
   \def\next@{\noindent@@\FN@\pretendspace@}\FNSS@\next@}}\else
  \gdef\Next@{\def~{\egroup\egroup\endgroup}}\fi
  }%
 }}
\def\litdefs@{\let\0\empty\let\1\litdelim@\def\ {\char32 }\litbs@@}%
\def\litdelimiter#1{%
 \edef\litdelim@{\char`#1}%
 \def\lit#1{\leavevmode\begingroup\litcodes@\litdefs@
  \tt\hyphenchar\tentt\m@ne\lit@}%
 \def\lit@##1#1{##1\endgroup\null}%
 \def\Lit#1{\ifhmode$$\abovedisplayskip\bigskipamount
  \abovedisplayshortskip\bigskipamount
  \belowdisplayskip\z@\belowdisplayshortskip\z@
  \postdisplaypenalty\@M
  $$\vskip-\baselineskip\else\bigskip\fi
  \begingroup\litlines@\z@
  \catcode`#1=\active\activate@0#1\Next@
  \def\displaybreak{\egroup\break\litlines@\z@\Lit@}%
  \def\allowdisplaybreak{\egroup\allowbreak\litlines@\z@\Lit@}%
  \def\allowdisplaybreaks{\egroup\allowbreak\interlitpenalty@\z@
   \litlines@\z@\Lit@}%
  \litcodes@\tt\catcode`\^^I=\active\letTAB@
  \obeylines\letM@\Lit@}%
 \def\Litbox##1=#1{\begingroup\ifodd##1\relax\aftergroup\global\fi
  \aftergroup\setbox\aftergroup##1\aftergroup\box\aftergroup\Litbox@
  \def\allowdisplaybreak{\egroup\allowbreak\litlines@\z@\Lit@}%
  \def\allowdisplaybreaks{\egroup\allowbreak\interlitpenalty@\z@
   \litlines@\z@\Lit@}%
  \catcode`#1=\active\activate@1#1\Next@
  \litcodes@\tt\catcode`\^^I=\active\letTAB@
  \obeylines\letM@\global\setbox\Litbox@\vbox\bgroup\litindent\z@%
  \litlines@\z@\Lit@}%
}
\newbox\titlebox@
\setbox\titlebox@\vbox{}
\rightadd@\title\to\overlonglist@
\def\title{\begingroup\Let@
 \global\setbox\titlebox@\vbox\bgroup\tabskip\hss@
 \halign to\hsize\bgroup\bf\hfil\ignorespaces##\unskip\hfil\cr}
\def\endtitle{\crcr\egroup\egroup\endgroup\overlong@false}
\newbox\authorbox@
\rightadd@\author\to\overlonglist@
\def\author{\begingroup\Let@
 \global\setbox\authorbox@\vbox\bgroup\tabskip\hss@
 \halign to\hsize\bgroup\rm\hfil\ignorespaces##\unskip\hfil\cr}
\def\endauthor{\crcr\egroup\egroup\endgroup\overlong@false}
\newbox\affilbox@
\def\affil{\begingroup\Let@
 \global\setbox\affilbox@\vbox\bgroup\tabskip\hss@
 \halign to\hsize\bgroup\rm\hfil\ignorespaces##\unskip\hfil\cr}%
\def\endaffil{\crcr\egroup\egroup\endgroup\overlong@false}
\let\date@\relax
\def\date#1{\gdef\date@{\ignorespaces#1\unskip}}
\def\today{\ifcase\month\or January\or February\or March\or April\or May\or
 June\or July\or August\or September\or October\or November\or December\fi
 \space\number\day, \number\year}
\def\maketitle{\hrule\height\z@\vskip-\topskip
 \vskip24\p@ plus12\p@ minus12\p@
 \unvbox\titlebox@
 \ifvoid\authorbox@\else\vskip12\p@ plus6\p@ minus3\p@\unvbox\authorbox@\fi
 \ifvoid\affilbox@\else\vskip10\p@ plus5\p@ minus2\p@\unvbox\affilbox@\fi
 \ifx\date@\relax\else\vskip6\p@ plus2\p@ minus\p@\centerline{\rm\date@}\fi
 \vskip18\p@ plus12\p@ minus6\p@}
\def\cite{%
 \DNii@(##1)##2{{\rm[}{##2}, {##1\/}{\rm]}}%
 \def\nextiii@##1{{\rm[}{##1\/}{\rm]}}%
 \DN@{\ifx\next(\expandafter\nextii@\else\expandafter\nextiii@\fi}%
 \FN@\next@}
\def\makebib@W{Bibliography}
\def\makebibp@W{References}
\def\makebib{\begingroup\rm\bigbreak\centerline{\smc\makebib@W}%
 \nobreak\medskip
 \sfcode`\.=\@m\everypar{}\parindent\z@
 \def\nopunct{\nopunct@true}\def\nospace{\nospace@true}%
 \nopunct@false\nospace@false
 \def\lkerns@{\null\kern\m@ne sp\kern\@ne sp}%
 \def\nkerns@{\null\kern-\tw@ sp\kern\tw@ sp}%
}

\newif\ifnoprepunct@
\newif\ifnoprespace@
\newif\ifnoquotes@
\def\noprepunct{\noprepunct@true}
\def\noprespace{\noprespace@true}
\def\noquotes{\noquotes@true}
\newbox\nobox@
\newbox\keybox@
\newbox\bybox@
\newbox\paperbox@
\newbox\paperinfobox@
\newbox\jourbox@
\newbox\volbox@
\newbox\issuebox@
\newbox\yrbox@
\newbox\pgbox@
\newbox\ppbox@
\newbox\bookbox@
\newbox\inbookbox@
\newbox\bookinfobox@
\newbox\publbox@
\newbox\publaddrbox@
\newbox\edbox@
\newbox\edsbox@
\newbox\langbox@
\newbox\translbox@
\newbox\finalinfobox@
\def\setbibinfo@#1{\edef\next@{\ifnopunct@1\else0\fi
 \ifnospace@1\else0\fi\ifnoprepunct@1\else0\fi\ifnoprespace@1\else0\fi
 \ifnoquotes@1\else0\fi}%
 \DNii@{00000}%
 \ifx\next@\nextii@\else\xdef\bibinfo@{\bibinfo@\the#1,\next@}%
 \fi}
\def\getbibinfo@#1{\ifx\bibinfo@\empty
 \let\next@0\let\nextii@0\let\nextiii@0\let\nextiv@0\let\nextv@0\else
 \edef\next@{\def
  \noexpand\next@####1\the#1,####2####3####4####5####6####7\noexpand\next@
  {\let\noexpand\next@####2\let\noexpand\nextii@####3%
  \let\noexpand\nextiii@####4\let\noexpand\nextiv@####5%
  \let\noexpand\nextv@####6}%
  \noexpand\next@\bibinfo@\the#1,00000\noexpand\next@}\next@
 \fi}
\newif\ifbookinquotes@
\def\bookinquotes{\bookinquotes@true}
\newif\ifpaperinquotes@
\def\paperinquotes{\paperinquotes@true}
\newif\ifininbook@
\def\ininbook{\ininbook@true}
\newif\ifopenquotes@
\def\closequotes@{\ifopenquotes@''\openquotes@false\fi}
\newif\ifbeginbib@
\newif\ifendbib@
\newif\ifprevjour@
\newif\ifprevbook@
\newdimen\bibindent@
\bibindent@20\p@
\def\bib{\global\let\bibinfo@\empty\global\let\translinfo@\relax\beginbib@true
 \begingroup\noindent@
 \hangindent\bibindent@\hangafter\@ne\bib@}
\def\v@id#1{\setbox#1\box\voidb@x}
\def\bib@{\v@id\nobox@\v@id\keybox@\v@id\bybox@\v@id\paperbox@
 \v@id\paperinfobox@\v@id\jourbox@\v@id\volbox@\v@id\issuebox@
 \v@id\yrbox@\v@id\pgbox@\v@id\ppbox@\v@id\bookbox@\v@id\inbookbox@
 \v@id\bookinfobox@\v@id\publbox@\v@id\publaddrbox@\v@id\edbox@
 \v@id\edsbox@\v@id\langbox@\v@id\translbox@\v@id\finalinfobox@
 \bgroup}
\def\Setnonemptybox@#1#2{\unskip\setbibinfo@#1\egroup#2%
 \def\aftergroup@{\ifdim\wd#1=\z@\setbox#1\box\voidb@x\fi}%
 \setbox#1\vbox\bgroup\aftergroup\aftergroup@\hsize\maxdimen\leftskip\z@
 \rightskip\z@\hbadness\@M\hfuzz\maxdimen\noindent}
\def\setnonemptybox@#1{\Setnonemptybox@#1\relax}
\def\no{\setnonemptybox@\nobox@}
\def\key{\setnonemptybox@\keybox@\bf}
\def\by{\setnonemptybox@\bybox@}
\def\bysame{\setnonemptybox@\bybox@\leaders\hrule\hskip3em\null}
\def\paper{\setnonemptybox@\paperbox@
 \ifpaperinquotes@\getbibinfo@\paperbox@
 \if\nextv@1\else``\fi\else\it\fi}
\def\paperinfo{\setnonemptybox@\paperinfobox@}
\def\jour{\Setnonemptybox@\jourbox@\prevjour@true}
\def\vol{\setnonemptybox@\volbox@\bf}
\def\issue{\setnonemptybox@\issuebox@}
\def\yr{\setnonemptybox@\yrbox@}
\def\toappear{\noprepunct\finalinfo(to appear)}
\def\pg{\setnonemptybox@\pgbox@}
\def\pp{\setnonemptybox@\ppbox@}
\def\book{\Setnonemptybox@\bookbox@\prevbook@true
 \ifbookinquotes@\getbibinfo@\bookbox@
 \if\nextv@1\else``\fi\else\it\fi}
\def\inbook{\Setnonemptybox@\inbookbox@\prevbook@true
 \ifininbook@ in \fi\ifbookinquotes@\getbibinfo@\inbookbox@
 \if\nextv@1\else``\fi\fi}
\def\bookinfo{\setnonemptybox@\bookinfobox@}
\def\publ{\setnonemptybox@\publbox@}
\def\publaddr{\setnonemptybox@\publaddrbox@}
\def\ed{\setnonemptybox@\edbox@}
\def\eds{\setnonemptybox@\edsbox@}
\def\lang{\setnonemptybox@\langbox@}
\def\finalinfo{\setnonemptybox@\finalinfobox@}
\def\setboxzl@{\setbox\z@\lastbox}
\def\getbox@#1{\setbox\z@\vbox{\vskip-\@M\p@
 \unvbox#1%
 \setboxzl@
 \global\setbox\@ne\hbox{\unhbox\z@\unskip\unskip\unpenalty}%
 \ifdim\lastskip=-\@M\p@\else
 \loop\ifdim\lastskip=-\@M\p@
 \else\unskip\unpenalty\setboxzl@
 \global\setbox\@ne\hbox{\unhbox\z@\unhbox\@ne}%
 \repeat\fi}%
 \unhbox\@ne}
\def\adjustpunct@#1{\count@\lastkern
 \ifnum\count@=\z@#1\closequotes@\else
 \ifnum\count@>\tw@#1\closequotes@\else
 \ifnum\count@<-\tw@#1\closequotes@\else
  \unkern\unkern\setboxzl@
  \skip@\lastskip\unskip
  \count@@\lastpenalty\unpenalty
  \ifnum\count@=\tw@\unskip\setboxzl@\fi
  \ifdim\skip@=\z@\else\hskip\skip@\fi
  #1\closequotes@
  \ifnum\count@=\tw@\null\hfill\fi
  \penalty\count@@
 \fi\fi\fi}
\def\prepunct@#1#2{\getbibinfo@#2%
 \ifnopunct@
 \else
  \if\nextiii@0\adjustpunct@#1\fi
 \fi
 \closequotes@
 \ifnospace@
 \else
  \if\nextiv@0\space\else\fi
 \fi
 \nopunct@false\nospace@false
 \if\next@1\nopunct@true\fi
 \if\nextii@1\nospace@true\fi}
\def\ppunbox@#1#2{\prepunct@{#1}#2%
 \getbox@#2}
\let\semicolon@;
\def\endbib@{%
 \ifbeginbib@
  \ifvoid\nobox@
   \ifvoid\keybox@\else\hbox to\bibindent@{[\getbox@\keybox@]\hss}\fi
  \else\hbox to\bibindent@{\hss\getbox@\nobox@. }\fi
  \ifvoid\bybox@\else\getbox@\bybox@\fi
 \else
  \nopunct@true
  \ifvoid\bybox@\else\ppunbox@\relax\bybox@\fi
 \fi
 \ifvoid\translbox@\else\ppunbox@,\translbox@\fi
 \ifvoid\paperbox@\else\ppunbox@,\paperbox@\ifpaperinquotes@
  \if\nextv@1\else\openquotes@true\fi\fi
 \fi
 \ifvoid\paperinfobox@\else\ppunbox@,\paperinfobox@\fi
 \test@false
 \ifvoid\jourbox@\else\test@true\ppunbox@,\jourbox@\fi
 \ifprevjour@\test@true\fi
 \iftest@
  \ifvoid\volbox@\else\ppunbox@\relax\volbox@\fi
  \ifvoid\issuebox@
   \else\prepunct@\relax\issuebox@ no.~\getbox@\issuebox@\fi
  \ifvoid\yrbox@\else\prepunct@\relax\yrbox@(\getbox@\yrbox@)\fi
  \ifvoid\ppbox@\else\ppunbox@,\ppbox@\fi
  \ifvoid\pgbox@\else\prepunct@,\pgbox@ p.~\getbox@\pgbox@\fi
 \fi
 \test@false
 \ifvoid\bookbox@\else\test@true\ppunbox@,\bookbox@\ifbookinquotes@
  \if\nextv@1\else\openquotes@true\fi\fi\fi
 \ifvoid\inbookbox@\else\test@true\ppunbox@,\inbookbox@\ifbookinquotes@
  \if\nextv@1\else\openquotes@true\fi\fi\fi
 \ifprevbook@\test@true\fi
 \iftest@
  \ifvoid\edbox@\else\prepunct@\relax\edbox@(\getbox@\edbox@, ed.)\fi
  \ifvoid\edsbox@\else\prepunct@\relax\edsbox@(\getbox@\edsbox@, eds.)\fi
  \ifvoid\bookinfobox@\else\ppunbox@,\bookinfobox@\fi
  \ifvoid\publbox@\else\ppunbox@,\publbox@\fi
  \ifvoid\publaddrbox@\else\ppunbox@,\publaddrbox@\fi
  \ifvoid\yrbox@\else\ppunbox@,\yrbox@\fi
  \ifvoid\ppbox@\else\prepunct@,\ppbox@ pp.~\getbox@\ppbox@\fi
  \ifvoid\pgbox@\else\prepunct@,\pgbox@ p.~\getbox@\pgbox@\fi
 \fi
 \ifvoid\finalinfobox@
  \ifendbib@
   \ifnopunct@\else.\closequotes@\fi
  \else
  \ifvoid\langbox@\else\space(\getbox@\langbox@)\fi
   \/\semicolon@\closequotes@
  \fi
 \else
  \ifendbib@
   \ppunbox@{.\spacefactor3000\relax}\finalinfobox@
    \ifnopunct@\else.\fi
  \else
   \ppunbox@,\finalinfobox@\/\semicolon@\fi
 \fi
 \ifvoid\langbox@\else\space(\getbox@\langbox@)\fi
}
\def\endbib{\unskip\egroup\endbib@true\endbib@\par\endgroup}
\def\morebib{\unskip\egroup
 \endbib@false\endbib@
 \global\let\bibinfo@\empty\beginbib@false
 \bib@}
\def\anotherbib{\unskip\egroup
 \endbib@false\endbib@
 \global\let\bibinfo@\empty\beginbib@false
 \prevjour@false\prevbook@false\bib@}
\def\transl{\unskip
 \xdef\translinfo@{\the\translbox@,\ifnopunct@1\else0\fi
 \ifnospace@1\else0\fi\ifnoprepunct@1\else0\fi\ifnoprespace@1\else0\fi0}%
 \egroup\endbib@false\endbib@
 \global\let\bibinfo@\translinfo@\beginbib@false
 \bib@
 \egroup
 \def\aftergroup@{\ifdim\wd\translbox@=\z@\setbox\translbox@\box\voidb@x\fi}%
 \setbox\translbox@\vbox\bgroup\aftergroup\aftergroup@
 \hsize\maxdimen\leftskip\z@\rightskip\z@\hbadness\@M\hfuzz\maxdimen
 \noindent}
\newwrite\auxwrite@
\newread\bbl@
\def\UseBibTeX{\immediate\openout\auxwrite@=\jobname.aux
 \let\cite\BTcite@
 \def\nocite##1{\immediate\write\auxwrite@{\string\citation{##1}}}%
 \def\bibliographystyle##1{\immediate\write\auxwrite@{\string
  \bibstyle{##1}}}%
 \def\bibliography@W{Bibliography}%
 \def\bibliography##1{\immediate\write\auxwrite@{\string\bibdata{##1}}%
  \immediate\openin\bbl@=\jobname.bbl
  \ifeof\bbl@
   \W@{No .bbl file}%
  \else
   \immediate\closein\bbl@
   \begingroup\input bibtex \input\jobname.bbl \endgroup
  \fi}%
 }
\def\BTcite@{%
 \DNii@(##1)##2{{\rm[}\BTcite@@##2,\BTcite@@{\rm, }{##1\/}{\rm]}%
  \immediate\write\auxwrite@{\string\citation{##2}}}%
 \def\nextiii@##1{{\rm[}\BTcite@@##1,\BTcite@@\/{\rm]}%
  \immediate\write\auxwrite@{\string\citation{##1}}}%
 \DN@{\ifx\next(\expandafter\nextii@\else\expandafter\nextiii@\fi}%
 \FN@\next@}%
\def\BTcite@@#1,{\BTcite@@@{#1}\FN@\BTcite@@@@}
\def\BTcite@@@@{\ifx\next\BTcite@@
 \expandafter\eat@\else{\rm, }\expandafter\BTcite@@\fi}
\catcode`\~=11
\def\BTcite@@@#1{\nolabel@\cite{#1}\relax
 \DNii@##1~##2\nextii@{##1}%
 \csL@{#1}\expandafter\nextii@\Next@\nextii@\fi}
\catcode`\~=\active

\def\beginthebibliography@#1{\rm\setboxz@h{#1\ }\bibindent@\wdz@
 \bigbreak\centerline{\smc\bibliography@W}\nobreak\medskip
 \sfcode`\.=\@m\everypar{}\parindent\z@}

\newif\iffigproofing@
\def\Figureproofing{\figproofing@true}
\def\noFigureproofing{\figproofing@false}
\newif\ifHby@
\def\Hbyw#1{\global\Hby@true\hbyw\vsize{#1}}
\def\hbyw#1#2{%
 \hbox{%
  \ifHby@
  \else
   \iffigproofing@
    \setbox\z@\vbox{\hrule\width5\p@}\ht\z@\z@
    \vbox to#1{\hrule\height5\p@\width.4\p@\vfil\hrule\height5\p@\width.4\p@}%
    \kern-.4\p@\rlap{\copy\z@}\raise#1\hbox{\rlap{\copy\z@}}%
   \fi
  \fi
  \vbox to#1{\hbox to#2{}\vfil}%
  \ifHby@
  \else
   \iffigproofing@
    \vbox to#1{\hrule\height5\p@\width.4\p@\vfil\hrule\height5\p@\width.4\p@}%
    \kern-.4\p@\llap{\copy\z@}\raise#1\hbox{\llap{\boxz@}}%
   \fi
  \fi}}
\newcount\island@C
\let\island@P\empty
\let\island@Q\empty
\def\island@S#1{#1\null.}
\let\island@N\arabic
\def\island@F{\rm}
\def\island@@@P{\csname\exxx@\islandtype@ @P\endcsname}
\def\island@@@Q{\csname\exxx@\islandtype@ @Q\endcsname}
\def\island@@@S{\csname\exxx@\islandtype@ @S\endcsname}
\def\island@@@N{\csname\exxx@\islandtype@ @N\endcsname}
\def\island@@@F{\csname\exxx@\islandtype@ @F\endcsname}
\def\island@@@C{\csname island@C\islandclass@\endcsname}
\newif\ifplace@
\newif\ifisland@
\def\island{%
 \ifplace@
  \DN@{\let\islandclass@\empty\def\islandtype@{\island}\FN@\island@}%
 \else
  \long\DN@##1\endisland{\Err@{\noexpand\island must be used after some
   type of \string\...place}}%
 \fi
 \next@}
\def\island@{\ifx\next\c\let\next@\island@c\else
 \DN@{\FN@\island@@}\fi\next@}
\def\island@@{\ifcat\bgroup\noexpand\next\let\next@\island@@@\else
 \DN@{\Err@{\noexpand\island must be followed by a {prefix} for
 \string\caption's}}\fi\next@}
\newbox\islandbox@
\newcount\captioncount@
\def\island@@@#1{\def\captionprefix@{#1}\captioncount@\z@
 \global\setbox\islandbox@\vbox\bgroup}
\def\island@c\c#1{%
 \ifplace@
 \DN@{\def\islandclass@{#1}%
  \expandafter\ifx\csname island@C#1\endcsname\relax
  \expandafter\newcount@\csname island@C#1\endcsname
   \global\csname island@C#1\endcsname\z@\fi
  \FNSS@\island@c@}%
 \else
 \DN@{\edef\next@{\long\def\noexpand\next@########1\expandafter\noexpand
  \csname end\exxx@\islandtype@\endcsname{\noexpand\Err@{\noexpand\noexpand
  \expandafter\noexpand
  \islandtype@ must be used after some type of \noexpand\string
   \noexpand\...place}}}\next@\next@}%
 \fi
 \next@}
\def\island@c@{%
 \ifcat\bgroup\noexpand\next
  \let\next@\island@c@@
 \else
  \DN@{\Err@{\noexpand\island\string\c{\expandafter\string\islandclass@} must
   be followed by a {prefix} for \string\caption's}}%
 \fi\next@}
\def\island@c@@#1{\def\captionprefix@{#1}%
 \captioncount@\z@\global\setbox\islandbox@\vbox\bgroup}
\rightadd@\caption\to\nofrillslist@
\newbox\captionbox@
\newbox\Captionbox@
\def\caption{%
 \ifnum\captioncount@=\z@
  \ifnopunct@
   \DN@{\egroup\nopunct@true}%
  \else
   \let\next@\egroup
  \fi
 \else
  \let\next@\relax
 \fi
 \next@
 \advance\captioncount@\@ne
 \FN@\caption@}
\def\caption@{\ifx\next"\expandafter\caption@q\else\expandafter\caption@@\fi}
\def\caption@q"#1"{\quoted@true
 {\noexpands@
 \let\pre\island@@@P\let\post\island@@@Q
 \let\style\island@@@S\let\numstyle\island@@@N
 \Qlabel@{#1}\let\style\relax\xdef\Qlabel@@@@{#1}}%
 \finishcaption@}
\def\caption@@{\quoted@false
 \global\advance\island@@@C\@ne
 {\noexpands@
 \xdef\Thelabel@@@{\number\island@@@C}%
 \xdefThelabel@\island@@@N
 \xdef\Thelabel@@@@{\island@@@P\Thelabel@\island@@@Q}%
 \xdefThelabel@@\island@@@S
 \xdef\Thepref@{\Thelabel@@@@}}%
 \finishcaption@}
\long\def\captionformat@#1#2#3{\rm\strut#1 {\island@@@F#2} #3%
 \punct@.\strut}
\long\def\widerthanisland@#1#2#3{\test@true\setbox\z@\vbox{\hsize\maxdimen
 \noindent@@\captionformat@{#1}{#2}{#3}\par\setboxzl@}%
 \ifdim\wdz@=\z@
  \global\setbox\captionbox@\hbox{\noset@\unlabel@
   \captionformat@{#1}{#2}{#3}}%
  \ifdim\wd\captionbox@>\wd\islandbox@\else\test@false\fi
 \fi}
\long\def\captionformat@@#1#2#3{\widerthanisland@{#1}{#2}{#3}%
 \iftest@
  \global\setbox\captionbox@\vbox{\hsize\wd\islandbox@
   \vskip-\parskip\noindent@@\noset@\unlabel@
   \captionformat@{#1}{#2}{#3}\par}%
 \else
  \global\setbox\captionbox@
   \hbox to\wd\islandbox@{\hfil\box\captionbox@\hfil}%
 \fi}
\long\def\finishcaption@#1{\def\entry@{#1}%
 {\locallabel@
 \captionformat@@
  {\expandafter\ignorespaces\captionprefix@\unskip}%
  {\ifx\thelabel@@\empty\unskip\else\thelabel@@\fi}%
  {\ignorespaces#1\unskip}%
 \ifnum\captioncount@=\@ne
  \global\setbox\islandbox@\vbox{\ticwrite@\vbox{\box\islandbox@}}%
  \global\setbox\Captionbox@\vbox{\box\captionbox@}%
 \else
  \global\setbox\islandbox@\vbox{\unvbox\islandbox@\setboxzl@
   \ticwrite@\boxz@}%
  \global\setbox\Captionbox@\vbox{\unvbox\Captionbox@
   \smallskip\box\captionbox@}%
 \fi}%
 \nopunct@false\nospace@false\ignorespaces}
\def\Sixtic@{\ifx\macdef@\empty\else
 \DN@##1##2\next@{\def\macdef@{##1##2}}%
 \expandafter\next@\macdef@\next@
 \edef\next@
  {\noexpand\six@\tic@\macdef@
  \space\space\space\space\space\space\space\space\space\space\space\space
  \noexpand\six@}%
 \next@\let\macdef@\relax\fi}
\def\ticwrite@{%
 \iftoc@
  {\noexpands@\let\style\relax
  \DN@{\island}%
  \edef\next@{\write\tic@{%
   \ifnopunct@\noexpand\noexpand\noexpand\nopunct\fi
   \ifx\islandtype@\next@\noexpand\noexpand\noexpand\island
    \noexpand\string\noexpand\c{\islandclass@}{\captionprefix@}%
     {\QorThelabel@@@@}\else\noexpand\noexpand\expandafter\noexpand
     \islandtype@{\QorThelabel@@@@}}\fi}%
  \next@}%
  \expandafter\unmacro@\meaning\entry@\unmacro@
  \Sixtic@
  \write\tic@{\noexpand\Page{\number\pageno}{\page@N}{\page@P}{\page@Q}^^J}%
 \fi}
\def\Htrim@#1{%
 \ifHby@
  \dimen@\vsize
  \ifnum\captioncount@=\z@
  \else
   \advance\dimen@-\ht\Captionbox@
   \advance\dimen@-#1%
  \fi
  \global\Hby@false
  \dimen@ii\wd\islandbox@
  \global\setbox\islandbox@\vbox
   {\unvbox\islandbox@\setboxzl@
   \vbox to\z@{\vss\boxz@}\nointerlineskip\hbyw\dimen@\dimen@ii}%
  \global\Hby@true
 \fi}
\newif\ifdata@
\def\iclasstest@#1{\DN@{#1}\ifx\next@\islandclass@
 \test@true\else\test@false\fi}
\skipdef\skipi@=1
\def\endisland{\ifnum\captioncount@=\z@\expandafter\egroup\fi
 \ifdata@
 \else
  \iclasstest@{T}%
  \iftest@
   {\rm\global\skipi@-\dp\strutbox}\global\advance\skipi@\bigskipamount
   \Htrim@\skipi@
   \global\setbox\islandbox@\vbox
    {\ifnum\captioncount@=\z@\else
     \box\Captionbox@
     \nointerlineskip
     \vskip\skipi@\fi
     \box\islandbox@}%
  \else
   {\rm\global\skipi@\dp\strutbox}\global\advance\skipi@\medskipamount
   \Htrim@\skipi@
   \global\setbox\islandbox@\vbox
    {\box\islandbox@
     \ifnum\captioncount@=\z@\else
     \nointerlineskip
     \vskip\skipi@
     \box\Captionbox@
     \fi}%
  \fi
  \ifHby@
  \else
   \dimen@\ht\islandbox@\advance\dimen@\dp\islandbox@
   \ifdim\dimen@>\vsize
    \DN@{\island}%
    \Err@{%
     \ifx\islandtype@\next@\noexpand\island\else
      \expandafter\noexpand\islandtype@\fi
     \ifnum\captioncount@=\z@\else
       with \noexpand\caption\fi
      is larger than page}%
     \ht\islandbox@=\vsize
   \fi
  \fi
 \fi
 \global\Hby@false\island@true}
\def\newisland#1\c#2#3{\define#1{}%
 \iftoc@\immediate\write\tic@{\noexpand\newisland\noexpand#1%
  \string\c{#2}{#3}^^J}\fi
 \expandafter\def\csname\exstring@#1@S\endcsname{\island@S}%
 \expandafter\def\csname\exstring@#1@N\endcsname{\island@N}%
 \expandafter\def\csname\exstring@#1@P\endcsname{\island@P}%
 \expandafter\def\csname\exstring@#1@Q\endcsname{\island@Q}%
 \expandafter\def\csname\exstring@#1@F\endcsname{\island@F}%
 \expandafter\def\csname end\exstring@#1\endcsname{\endisland}%
 \expandafter
 \ifx\csname island@C#2\endcsname\relax
  \expandafter\newcount@\csname island@C#2\endcsname
  \global\csname island@C#2\endcsname\z@
 \fi
 \edef\next@{\noexpand\expandafter\noexpand\let\noexpand
  \csname\exstring@#1@C\noexpand\endcsname
  \csname island@C#2\endcsname}%
 \next@
 \def#1{\def\islandtype@{#1}\island@c\c{#2}{#3}}}
\newisland\Figure\c{F}{Figure}
\newisland\Table\c{T}{Table}
\newbox\islandboxi
\newbox\islandboxii
\newbox\islandboxiii
\newbox\captionboxi
\newbox\captionboxii
\newbox\captionboxiii
\long\def\islandpairdata#1#2{{\data@true
 \place@true
 #1%
 \global\setbox\islandboxi\box\islandbox@
 \global\setbox\captionboxi\box\Captionbox@
 #2%
 \global\setbox\islandboxii\box\islandbox@
 \global\setbox\captionboxii\box\Captionbox@
 }}
\long\def\islandpairbox#1#2{\islandpairdata{#1}{#2}%
 \dimen@\ht\captionboxi
 \ifdim\ht\captionboxii>\dimen@\dimen@\ht\captionboxii\fi
 \ifdim\dimen@>\z@
  \ifdim\ht\captionboxi<\dimen@
   \global\setbox\captionboxi\vbox to\dimen@{\unvbox\captionboxi\vfil}\fi
  \ifdim\ht\captionboxii<\dimen@
   \global\setbox\captionboxii\vbox to\dimen@{\unvbox\captionboxii\vfil}\fi
 \fi
 \global\setbox\islandbox@\vbox
 {\hbox to\hsize{\hfil\box\islandboxi\hfil\box\islandboxii\hfil}%
 \ifdim\dimen@>\z@\nointerlineskip
 {\rm\global\skipi@\dp\strutbox}\global\advance\skipi@\medskipamount
  \vskip\skipi@
  \hbox to\hsize{\hfil\box\captionboxi\hfil\box\captionboxii\hfil}\fi}}
\long\def\islandpairboxa#1#2{\islandpairdata{#1}{#2}%
 \dimen@\ht\captionboxi
 \ifdim\ht\captionboxii>\dimen@\dimen@\ht\captionboxii\fi
 \ifdim\dimen@>\z@
  \ifdim\ht\captionboxi<\dimen@
   \global\setbox\captionboxi\vbox to\dimen@{\vfil\unvbox\captionboxi}\fi
  \ifdim\ht\captionboxii<\dimen@
   \global\setbox\captionboxii\vbox to\dimen@{\vfil\unvbox\captionboxii}\fi
 \fi
 \dimen@ii\ht\islandboxi
 \ifdim\ht\islandboxii>\dimen@ii \dimen@ii\ht\islandboxii\fi
 \ifdim\dimen@ii>\z@
  \ifdim\ht\islandboxi<\dimen@ii
   \global\setbox\islandboxi\vbox to\dimen@ii{\box\islandboxi\vfil}\fi
  \ifdim\ht\islandboxii<\dimen@ii
   \global\setbox\islandboxii\vbox to\dimen@ii{\box\islandboxii\vfil}\fi
 \fi
 \global\setbox\islandbox@\vbox{\ifdim\dimen@>\z@
  \hbox to\hsize{\hfil\box\captionboxi\hfil\box\captionboxii\hfil}%
  \nointerlineskip{\rm\global\skipi@-\dp\strutbox}%
  \global\advance\skipi@\bigskipamount\vskip\skipi@\fi
  \hbox to\hsize{\hfil\box\islandboxi\hfil\box\islandboxii\hfil}}}
\long\def\islandtripledata#1#2#3{{\data@true\place@true
 #1%
 \global\setbox\islandboxi\box\islandbox@
 \global\setbox\captionboxi\box\Captionbox@
 #2%
 \global\setbox\islandboxii\box\islandbox@
 \global\setbox\captionboxii\box\Captionbox@
 #3%
 \global\setbox\islandboxiii\box\islandbox@
 \global\setbox\captionboxiii\box\Captionbox@
 }}
\long\def\islandtriplebox#1#2#3{\islandtripledata{#1}{#2}{#3}%
 \dimen@\ht\captionboxi
 \ifdim\ht\captionboxii>\dimen@ \dimen@\ht\captionboxii\fi
 \ifdim\ht\captionboxiii>\dimen@ \dimen@\ht\captionboxiii\fi
 \ifdim\dimen@>\z@
  \ifdim\ht\captionboxi<\dimen@
   \global\setbox\captionboxi\vbox to\dimen@{\unvbox\captionboxi\vfil}\fi
  \ifdim\ht\captionboxii<\dimen@
   \global\setbox\captionboxii\vbox to\dimen@{\unvbox\captionboxii\vfil}\fi
  \ifdim\ht\captionboxiii<\dimen@
   \global\setbox\captionboxiii\vbox to\dimen@{\unvbox\captionboxiii\vfil}\fi
 \fi
 \global\setbox\islandbox@\vbox
  {\hbox to\hsize{\hfil\box\islandboxi\hfil\box\islandboxii\hfil
   \box\islandboxiii\hfil}%
 \ifdim\dimen@>\z@\nointerlineskip
  {\rm\global\skipi@\dp\strutbox}\global\advance\skipi@\medskipamount
  \vskip\skipi@
  \hbox to\hsize{\hfil\box\captionboxi\hfil\box\captionboxii\hfil
   \box\captionboxiii\hfil}\fi}}
\def\islandtripleboxa#1#2#3{\islandtripledata{#1}{#2}{#3}%
 \dimen@\ht\captionboxi
 \ifdim\ht\captionboxii>\dimen@ \dimen@\ht\captionboxii\fi
 \ifdim\ht\captionboxiii>\dimen@ \dimen@\ht\captionboxiii\fi
 \ifdim\dimen@>\z@
  \ifdim\ht\captionboxi<\dimen@
   \global\setbox\captionboxi\vbox to\dimen@{\vfil\unvbox\captionboxi}\fi
  \ifdim\ht\captionboxii<\dimen@
   \global\setbox\captionboxii\vbox to\dimen@{\vfil\unvbox\captionboxii}\fi
  \ifdim\ht\captionboxiii<\dimen@
   \global\setbox\captionboxiii\vbox to\dimen@{\vfil\unvbox\captionboxiii}\fi
 \fi
 \dimen@ii\ht\islandboxi
 \ifdim\ht\islandboxii>\dimen@ii \dimen@ii\ht\islandboxii\fi
 \ifdim\ht\islandboxiii>\dimen@ii \dimen@ii\ht\islandboxiii\fi
 \ifdim\dimen@ii>\z@
  \ifdim\ht\islandboxi<\dimen@ii
   \global\setbox\islandboxi\vbox to\dimen@ii{\box\islandboxi\vfil}\fi
  \ifdim\ht\islandboxii<\dimen@ii
   \global\setbox\islandboxii\vbox to\dimen@ii{\box\islandboxii\vfil}\fi
  \ifdim\ht\islandboxiii<\dimen@ii
   \global\setbox\islandboxiii\vbox to\dimen@ii{\box\islandboxiii\vfil}\fi
 \fi
 \global\setbox\islandbox@\vbox
  {\ifdim\dimen@>\z@
  \hbox to\hsize{\hfil\box\captionboxi\hfil\box\captionboxii\hfil
   \box\captionboxiii\hfil}%
  \nointerlineskip{\rm\global\skipi@-\dp\strutbox}%
  \global\advance\skipi@\bigskipamount\vskip\skipi@\fi
  \hbox to\hsize{\hfil\box\islandboxi\hfil\box\islandboxii\hfil
   \box\islandboxiii\hfil}}}
\def\Figurepair#1\and#2\endFigurepair{\island@true
 \islandpairbox{\Figure#1\endFigure}{\Figure#2\endFigure}}
\def\Figuretriple#1\and#2\and#3\endFiguretriple{\island@true
 \islandtriplebox{\Figure#1\endFigure}{\Figure#2\endFigure}%
  {\Figure#3\endFigure}}
\def\Tablepair#1\and#2\endTablepair{\island@true
 \islandpairboxa{\Table#1\endTable}{\Table#2\endTable}}
\def\Tabletriple#1\and#2\and#3\endTabletriple{\island@true
 \islandtripleboxa{\Table#1\endTable}{\Table#2\endTable}%
 {\Table#3\endTable}}
\def\place#1{\place@true\island@false
 #1%
 \ifisland@
  \box\islandbox@
 \else
  \Err@{Whoa ... there's no \string\Figure, \string\Table,
   etc., here}%
 \fi
 \place@false}
\newskip\belowtopfigskip
\belowtopfigskip 15\p@ plus 5\p@ minus5\p@
\newskip\abovebotfigskip
\abovebotfigskip 18\p@ plus 6\p@ minus6\p@
\newdimen\minpagesize
\minpagesize 5pc
\dimen@\belowtopfigskip
\advance\dimen@-\abovebotfigskip
\skip\topins\dimen@
\dimen\topins\z@
\newcount\topinscount@
\newbox\topinsdims@
\def\storedim@{\global\setbox\topinsdims@
 \vbox{\hbox to\dimen@{}\unvbox\topinsdims@}}
\def\advancedimtopins@{%
 \ifnum\pageno=\@ne
 \else
   \advance\dimen@\dimen\topins
   \global\dimen\topins\dimen@
 \fi}
\newcount\flipcount@
\def\fliptopins@{%
 \global\flipcount@\z@
 \ifvoid\topins\else
 \setbox\z@\vbox
  {\vskip\p@
   \unvbox\topins
   \global\setbox\topins\vbox{}%
   \loop
    \test@false
    \ifdim\lastskip=\z@\unskip
     \ifdim\lastskip=\z@
      \test@true\fi\fi
    \iftest@
    \global\advance\flipcount@\@ne
    \setboxzl@
    \global\setbox\topins\vbox{\unvbox\topins\boxz@}%
    \unpenalty
   \repeat}\fi}
\newif\ifPar@
\newcount\Parcount@
\newbox\Parbox@
\expandafter\newbox\csname Parfigbox1\endcsname
\expandafter\newbox\csname Parfigbox2\endcsname
\expandafter\newbox\csname Parfigbox3\endcsname
\expandafter\newbox\csname Parfigbox4\endcsname
\expandafter\newbox\csname Parfigbox5\endcsname
\expandafter\newdimen\csname Parprev1\endcsname
\expandafter\newdimen\csname Parprev2\endcsname
\expandafter\newdimen\csname Parprev3\endcsname
\expandafter\newdimen\csname Parprev4\endcsname
\expandafter\newdimen\csname Parprev5\endcsname
\expandafter\newdimen\csname Parprev6\endcsname
\def\Par{\par\global\csname Parprev1\endcsname\prevdepth
 \global\Parcount@\@ne
 \global\Par@true\global\let\Parlist@\empty
 \global\setbox\Parbox@\vbox\bgroup\break}
\def\place@#1#2{%
 \ifisland@
  \ifhmode
   \ifPar@
    \ifnum\Parcount@>5
     \Err@{Only 5 \string\place's allowed per
      \string\Par...\noexpand\endPar paragraph}%
    \else
     \expandafter\expandafter\expandafter
      \global\expandafter\setbox
       \csname Parfigbox\number\Parcount@\endcsname\box\islandbox@
     \global\advance\Parcount@\@ne
     \xdef\Parlist@{\Parlist@#1}%
    \fi
   \else
    \vadjust{#2}%
   \fi
  \else
   #2%
  \fi
 \else
  \Err@{Whoa ... there's no \string\Figure,
   \string\Table, etc., here}%
 \fi
 \place@false}
\long\def\Aplace#1{\prevanish@
 \place@true\island@false
 #1%
 \place@ a\Aplace@
 \postvanish@}
\long\def\AAplace#1{\prevanish@\place@true\island@false
 #1%
 \place@ A\AAplace@
 \postvanish@}
\newif\ifAA@
\def\AAplace@{\AA@true\Aplace@\AA@false}
\let\AAlist@\empty
\def\Aplace@{\allowbreak
 \dimen@=\ht\islandbox@
 \advance\dimen@\abovebotfigskip
 \ht\islandbox@\dimen@
 \advance\dimen@\dp\islandbox@
 \storedim@
 \ifAA@
  \xdef\AAlist@{\AAlist@1}%
  \advancedimtopins@
 \else
  \xdef\AAlist@{\AAlist@0}%
  \ifnum\topinscount@>\@ne\else\advancedimtopins@\fi
 \fi
 \insert\topins{\penalty\z@\splittopskip\z@\floatingpenalty\z@
  \box\islandbox@}%
 \global\advance\topinscount@\@ne}
\long\def\Bplace#1{\prevanish@\place@true\island@false
 #1%
 \place@ b\Bplace@
 \postvanish@}
\def\Bplace@{\allowbreak
 \ifnum\topinscount@=\z@
  \setbox\z@\vbox{\vbox to-\belowtopfigskip{}}%
  \dimen@-\skip\topins
  \ht\z@\dimen@
  \storedim@
  \advancedimtopins@
  \insert\topins{\boxz@}%
  \global\advance\topinscount@\@ne
  \xdef\AAlist@{\AAlist@0}%
 \fi
 \dimen@\ht\islandbox@
 \advance\dimen@\abovebotfigskip
 \ht\islandbox@\dimen@
 \advance\dimen@\dp\islandbox@
 \storedim@
 \xdef\AAlist@{\AAlist@0}%
 \ifnum\topinscount@>\@ne\else\advancedimtopins@\fi
 \insert\topins{\penalty\z@\splittopskip\z@
  \floatingpenalty\z@
  \box\islandbox@}%
 \global\advance\topinscount@\@ne}
\def\breakisland@{\global\setbox\@ne\lastbox\global\skipi@\lastskip\unskip
 \global\setbox\thr@@\lastbox}%
\def\printisland@{\centerline{\box\thr@@}\nobreak\nointerlineskip
 \vskip\skipi@
 \ifdim\ht\@ne<\z@\box\@ne\else\centerline{\box\@ne}\fi}
\def\bottomfigs@{%
 \count@\@ne
 \loop
  \ifnum\count@<\flipcount@
  \nointerlineskip
  \vskip\abovebotfigskip
  \global\setbox\topins\vbox{\unvbox\topins\setboxzl@
   \unvbox\z@
   \breakisland@}%
  \printisland@
  \advance\count@\@ne
  \repeat}
\def\resetdimtopins@{%
 \global\advance\topinscount@-\flipcount@
 \global\setbox\topinsdims@\vbox
  {\unvbox\topinsdims@
   \count@\z@
   \DN@##1##2\next@{\gdef\AAlist@{##2}}%
   \loop
    \ifnum\count@<\flipcount@\setboxzl@
    \expandafter\next@\AAlist@\next@
    \advance\count@\@ne
    \repeat
   \dimen@\z@
   \count@\z@
   \setbox\tw@\vbox{}%
   \edef\nextiii@{\AAlist@}%
   \DN@##1##2\next@{\DNii@{##1}\def\nextiii@{##2}}%
   \loop
    \test@false
    \ifnum\count@<\topinscount@
    \expandafter\next@\nextiii@\next@
     \ifnum\count@<\tw@
      \test@true
     \else
      \if\nextii@ 1\test@true\fi
     \fi
    \fi
    \iftest@
     \setboxzl@
     \advance\dimen@\wdz@
     \setbox\tw@\vbox{\boxz@\unvbox\tw@}%
     \advance\count@\@ne
    \repeat
    \unvbox\tw@
    \global\dimen\topins\dimen@}}
\def\Place@#1#2{%
 \ifisland@
  \ifhmode
   \ifPar@
    \ifnum\Parcount@>5
     \Err@{Only 5 \string\place's allowed per
       \string\Par...\noexpand\endPar paragraph}%
    \else
     \expandafter\expandafter\expandafter\global\expandafter\setbox
      \csname Parfigbox\number\Parcount@\endcsname\box\islandbox@
     \global\advance\Parcount@\@ne
     \xdef\Parlist@{\Parlist@#1}%
     \vadjust{\break}%
    \fi
   \else
    \Err@{\noexpand#2allowed only in a \string\Par...\noexpand\endPar
     paragraph}%
   \fi
  \else
   #2%
  \fi
 \else
  \Err@{Who ... there's no \string\Figure, \string\Table,
   etc., here}%
 \fi
 \place@false}
\newif\ifC@
\newdimen\Cdim@
\long\def\Cplace#1{\prevanish@\place@true\island@false
 #1%
 \Place@ c\Cplace@
 \postvanish@}
\def\Cplace@{\allowbreak
 \ifnum\topinscount@>\z@\else
  \global\C@true\global\Cdim@\pagetotal\fi
 \Aplace@}
\long\def\Mplace#1{\prevanish@\place@true\island@false
 #1%
 \Place@ m\Mplace@
 \postvanish@}
\long\def\MXplace#1{\prevanish@\place@true\island@false
 #1%
 \Place@ M\MXplace@
 \postvanish@}
\newif\ifMX@
\def\MXplace@{\MX@true\Mplace@\MX@false}
\def\Mplace@{\allowbreak
 \dimen@\ht\islandbox@\advance\dimen@\dp\islandbox@
 \ifdim\pagetotal=\z@\else
  \ifdim\lastskip<\abovebotfigskip\advance\dimen@\abovebotfigskip
  \advance\dimen@-\lastskip\fi
 \fi
 \advance\dimen@\pagetotal
 \ifdim\dimen@>\pagegoal
  \Aplace@
 \else
  \nointerlineskip
  \ifdim\lastskip<\abovebotfigskip\removelastskip\vskip\abovebotfigskip\fi
  \setbox\z@\vbox{\unvbox\islandbox@
   \breakisland@}%
  \printisland@
  \ifnum\topinscount@=\z@
   \setbox\z@\vbox{\vbox to-\belowtopfigskip{}}%
   \dimen@-\skip\topins
   \ht\z@\dimen@
   \storedim@
   \advancedimtopins@
   \insert\topins{\boxz@}%
   \global\advance\topinscount@\@ne
   \xdef\AAlist@{\AAlist@0}%
  \fi
  \ifMX@
   \ifnum\topinscount@=\@ne
    \setbox\z@\vbox{\vbox to-\abovebotfigskip{}}%
    \ht\z@\z@
    \dimen@\z@
    \storedim@
    \advancedimtopins@
    \insert\topins{\boxz@}%
    \global\advance\topinscount@\@ne
    \xdef\AAlist@{\AAlist@0}%
   \fi
  \fi
  \nointerlineskip
  \vskip\belowtopfigskip
 \fi}
\expandafter\newbox\csname Parbox1\endcsname
\expandafter\newbox\csname Parbox2\endcsname
\expandafter\newbox\csname Parbox3\endcsname
\expandafter\newbox\csname Parbox4\endcsname
\expandafter\newbox\csname Parbox5\endcsname
\def\endPar{\egroup
 \count@\@ne
 {\vbadness\@M\vfuzz\maxdimen\splitmaxdepth\maxdimen\splittopskip\ht\strutbox
 \setbox\z@\vsplit\Parbox@ to\ht\Parbox@
 \loop
  \ifnum\count@<\Parcount@
  \expandafter\expandafter\expandafter\global\expandafter\setbox
   \csname Parbox\number\count@\endcsname\vsplit\Parbox@ to\ht\Parbox@
  \count@@\count@\advance\count@@\@ne
  \global\csname Parprev\number\count@@\endcsname
   \dp\csname Parbox\number\count@\endcsname
  \advance\count@\@ne
  \repeat}%
 \vskip\parskip
 \count@\@ne
 \def\nextv@##1##2\nextv@{\DN@{##1}\gdef\Parlist@{##2}}%
 \loop
  \ifnum\count@<\Parcount@
   \dimen@\csname Parprev\number\count@\endcsname
   \advance\dimen@\ht\strutbox
   \ifdim\dimen@<\baselineskip
    \advance\dimen@-\baselineskip\vskip-\dimen@
   \else
    \vskip\lineskip
   \fi
   \unvbox\csname Parbox\number\count@\endcsname
   \global\setbox\islandbox@\box\csname Parfigbox\number\count@\endcsname
   \expandafter\nextv@\Parlist@\nextv@
   \if a\next@\Aplace@\else
   \if A\next@\AAplace@\else
   \if b\next@\Bplace@\else
   \if c\next@\Cplace@\else
   \if m\next@\Mplace@\else
   \if M\next@\MXplace@\fi\fi\fi\fi\fi\fi
  \advance\count@\@ne
  \repeat
 \global\Par@false
 \ifvoid\Parbox@
  \prevdepth\csname Parprev\number\count@\endcsname
 \else
  \dimen@\csname Parprev\number\count@\endcsname\advance\dimen@\ht\strutbox
  \ifdim\dimen@<\baselineskip
    \advance\dimen@-\baselineskip\vskip-\dimen@
  \else
    \vskip\lineskip
  \fi
  \dimen@\dp\Parbox@
  \unvbox\Parbox@
  \prevdepth\dimen@
 \fi}
\def\folio{{\page@F\page@S{\page@P\page@N{\number\page@C}\page@Q}}}
\def\advancepageno{\global\advance\pageno\@ne}
\newif\ifspecialsplit@
\newbox\outbox@
\let\shipout@\shipout
\def\plainoutput{\specialsplit@false\ifvoid\topins\else\ifdim\ht\topins=\z@
 \specialsplit@true\advance\minpagesize-\skip\topins\fi\fi
 \fliptopins@
 \setbox\outbox@\vbox{\makeheadline\pagebody\makefootline}%
 {\noexpands@\let\style\relax
 \shipout@\box\outbox@}%
 \advancepageno
 \resetdimtopins@
 \ifvoid\@cclv\else\unvbox\@cclv\penalty\outputpenalty\fi
 \ifnum\outputpenalty>-\@MM\else\dosupereject\fi}
\def\pagebody{\vbox to\vsize{\boxmaxdepth\maxdepth
 \ifvoid\margin@\else
 \rlap{\kern\hsize\vbox to\z@{\kern4\p@\box\margin@\vss}}\fi
 \pagecontents}}
\newif\ifonlytop@
\def\pagecontents{%
 \onlytop@false
 \ifdim\ht\@cclv<\minpagesize\ifnum\flipcount@<\tw@\ifvoid\footins
  \onlytop@true\fi\fi\fi
 \test@false
 \ifC@
  \ifnum\flipcount@=\@ne
   \global\multiply\Cdim@\tw@
   \ifdim\Cdim@>\ht\@cclv
    \test@true
   \fi
  \fi
 \fi
 \global\C@false
 \iftest@
  \dimen@\ht\@cclv
  \advance\dimen@\skip\topins
  {\vfuzz\maxdimen\vbadness\@M
  \splitmaxdepth\maxdepth\splittopskip\topskip
  \setbox\z@\vsplit\@cclv to\dimen@
  \unvbox\z@}%
  \global\setbox\topins\vbox{\unvbox\topins
   \global\setbox\@ne\lastbox}%
  \setbox\z@\vbox{\unvbox\@ne
   \breakisland@}%
  \nointerlineskip
  \vskip\abovebotfigskip
  \printisland@
 \else
  \ifnum\flipcount@>\z@
   \global\setbox\topins\vbox{\unvbox\topins\global\setbox\@ne\lastbox}%
   \setbox\z@\vbox{\unvbox\@ne
    \breakisland@}%
   \printisland@
   \ifonlytop@\kern-\prevdepth\vfill\else\vskip\belowtopfigskip\fi
  \fi
 \fi
 \ifdim\ht\@cclv<\minpagesize
  \ifonlytop@\else\vfill\fi
 \else
  \ifspecialsplit@
   {\vfuzz\maxdimen\vbadness\@M
   \splitmaxdepth\maxdepth\splittopskip\topskip
   \dimen@ii\ht\@cclv \advance\dimen@ii\skip\topins
   \setbox\z@\vsplit\@cclv to\dimen@ii
   \unvbox\z@}%
  \else
   \unvbox\@cclv
  \fi
 \fi
 \bottomfigs@
 \ifvoid\footins\else\vskip\skip\footins\footnoterule\unvbox\footins\fi}
\newread\readdata@
\def\readthedata@#1{\expandafter
 \ifx\csname#1@D\endcsname\relax
  \immediate\openin\readdata@=#1.dat
  \ifeof\readdata@
   \Err@{No file #1.dat}%
  \else
   {\endlinechar\m@ne\gdef\Next@{}%
   \DNii@##1 ##2 ##3pt{\global\data@ht##1\global\data@dp##2%
    \global\data@wd##3pt}%
   \loop
    \ifeof\readdata@
    \else
    \read\readdata@ to\next@
    \ifx\next@\empty\else
     \edef\next@{\expandafter\nextii@\next@}%
     \expandafter\rightadd@\next@\to\Next@
    \fi
    \repeat}%
   \immediate\closein\readdata@
   \expandafter\expandafter\expandafter\global\expandafter
    \let\csname#1@D\endcsname\Next@\global\let\Next@\relax
  \fi
 \fi}
\newdimen\data@ht
\newdimen\data@dp
\newdimen\data@wd
\newif\ifgetdata@
\def\getdata@#1#2{\global\getdata@true\count@#2\relax
 {\let\\\or\xdef\Next@{\ifcase\number\count@#1\else
 \global\noexpand\getdata@false\fi}}\Next@}
\def\paste#1#2{\readthedata@{#1}%
 \getdata@{\csname#1@D\endcsname}{#2}%
 \ifgetdata@
 \dimen@\data@ht \advance\dimen@\data@dp
  \hbox{\special{dvipaste: #1 #2}%
   \lower\data@dp\vbox to\dimen@{\hbox to\data@wd{}\vfil}}%
 \else
  {\lccode`\Z=`\#\lccode`\N=`\N\lccode`\F=`\F%
   \lowercase{\Err@{No data for File [#1], Z#2}}}%
 \fi}
\newdimen\httable
\newdimen\dptable
\newdimen\wdtable
\def\measuretable#1#2{\readthedata@{#1}%
 \getdata@{\csname#1@D\endcsname}{#2}%
 \ifgetdata@
  \httable\data@ht \dptable\data@dp \wdtable\data@wd
 \else
  {\lccode`\Z=`\#\lccode`\N=`\N\lccode`\F=`\F%
  \lowercase{\Err@{No data for File [#1], Z#2}}}%
 \fi}
\def\East#1#2{\setboxz@h{$\m@th\ssize\;{#1}\;\;$}%
 \setbox\tw@\hbox{$\m@th\ssize\;{#2}\;\;$}\setbox4=\hbox{$\m@th#2$}%
 \dimen@\minaw@
 \ifdim\wdz@>\dimen@\dimen@\wdz@\fi\ifdim\wd\tw@>\dimen@\dimen@\wd\tw@\fi
 \ifdim\wd4 >\z@
  \mathrel{\mathop{\hbox to\dimen@{\rightarrowfill}}\limits^{#1}_{#2}}%
 \else
  \mathrel{\mathop{\hbox to\dimen@{\rightarrowfill}}\limits^{#1}}%
 \fi}
\def\West#1#2{\setboxz@h{$\m@th\ssize\;\;{#1}\;$}%
 \setbox\tw@\hbox{$\m@th\ssize\;\;{#2}\;$}\setbox4=\hbox{$\m@th#2$}%
 \dimen@\minaw@
 \ifdim\wdz@>\dimen@\dimen@\wdz@\fi\ifdim\wd\tw@>\dimen@\dimen@\wd\tw@\fi
 \ifdim\wd4 >\z@
  \mathrel{\mathop{\hbox to\dimen@{\leftarrowfill}}\limits^{#1}_{#2}}%
 \else
  \mathrel{\mathop{\hbox to\dimen@{\leftarrowfill}}\limits^{#1}}%
 \fi}
\font\arrow@i=lams1
\font\arrow@ii=lams2
\font\arrow@iii=lams3
\font\arrow@iv=lams4
\font\arrow@v=lams5
\newdimen\standardcgap
\standardcgap40\p@
\newdimen\hunit
\hunit\tw@\p@
\newdimen\standardrgap
\standardrgap32\p@
\newdimen\vunit
\vunit1.6\p@
\def\Cgaps#1{\RIfM@
 \standardcgap#1\standardcgap\relax\hunit#1\hunit\relax
 \else\nonmatherr@\Cgaps\fi}
\def\Rgaps#1{\RIfM@
 \standardrgap#1\standardrgap\relax\vunit#1\vunit\relax
 \else\nonmatherr@\Rgaps\fi}
\newdimen\getdim@
\def\getcgap@#1{\ifcase#1\or\getdim@\z@\else\getdim@\standardcgap\fi}
\def\getrgap@#1{\ifcase#1\getdim@\z@\else\getdim@\standardrgap\fi}
\def\cgaps{\RIfM@\expandafter\cgaps@\else\expandafter\nonmatherr@
 \expandafter\cgaps\fi}
\def\cgaps@{\ifnum\catcode`\;=\active\expandafter\cgapsA@\else
 \expandafter\cgapsO@\fi}
\def\cgapsO@#1{\toks@{\ifcase\i@\or\getdim@=\z@}%
 \gaps@@\standardcgap#1;\gaps@@\gaps@@
 \edef\next@{\the\toks@\noexpand\else\noexpand\getdim@\noexpand\standardcgap
  \noexpand\fi}%
 \toks@=\expandafter{\next@}%
 \edef\getcgap@##1{\i@##1\relax\the\toks@}\toks@{}}
{\catcode`\;=\active
 \gdef\cgapsA@#1{\toks@{\ifcase\i@\or\getdim@=\z@}%
 \gaps@@\standardcgap#1;\gaps@@\gaps@@
 \edef\next@{\the\toks@\noexpand\else\noexpand\getdim@\noexpand\standardcgap
  \noexpand\fi}%
 \toks@=\expandafter{\next@}%
 \edef\getcgap@##1{\i@##1\relax\the\toks@}\toks@{}}
}
\def\Gaps@@{\gaps@@}
\def\gaps@@#1#2;#3{\mgaps@#1#2\mgaps@
 \edef\next@{\the\toks@\noexpand\or\noexpand\getdim@
  \noexpand#1\the\mgapstoks@@}%
 \toks@\expandafter{\next@}%
 \DN@{#3}%
 \ifx\next@\Gaps@@\def\next@##1\gaps@@{}\else
  \def\next@{\gaps@@#1#3}\fi\next@}
{\catcode`\;=\active
 \gdef\rgaps#1{\RIfM@{\ifnum\catcode`\;=\active\def;{\string;}\fi
   \xdef\Next@{\noexpand\rgaps@{#1}}}%
  \Next@\edef\getrgap@##1{\i@##1\relax\the\toks@}\toks@{}\else
  \nonmatherr@\rgaps\fi}
}
\def\rgaps@#1{\toks@{\ifcase\i@\getdim@=\z@}%
 \gaps@@\standardrgap#1;\gaps@@\gaps@@
 \edef\next@{\the\toks@\noexpand\else\noexpand\getdim@\noexpand\standardrgap
  \noexpand\fi}%
 \toks@=\expandafter{\next@}}
\newbox\ZER@
\def\mgaps@#1{\let\mgapsnext@#1\FNSS@\mgaps@@}
\def\mgaps@@{\ifx\next\w\expandafter\mgaps@@@\else
 \expandafter\mgaps@@@@\fi}
\newtoks\mgapstoks@@
\def\mgaps@@@@#1\mgaps@{\getdim@\mgapsnext@\getdim@#1\getdim@
 \edef\next@{\noexpand\getdim@\the\getdim@}%
 \mgapstoks@@\expandafter{\next@}}
\def\mgaps@@@\w#1#2\mgaps@{\mgaps@@@@#2\mgaps@
 \setbox\ZER@\hbox{$\m@th\hskip15\p@\tsize@#1$}%
 \dimen@\wd\ZER@
 \ifdim\dimen@>\getdim@\getdim@\dimen@\fi
 \edef\next@{\noexpand\getdim@\the\getdim@}%
 \mgapstoks@@\expandafter{\next@}}
\def\changewidth#1#2{\setbox\ZER@{$\m@th#2}%
 \hbox to\wd\ZER@{\hss$\m@th#1$\hss}}
\atdef@({\FN@\ARROW@}
\def\ARROW@{\ifx\next)\let\next@\OPTIONS@\else
 \DN@{\csname\string @(\endcsname}\fi\next@}
\newif\ifoptions@
\def\OPTIONS@){\ifoptions@\let\next@\relax\else
 \DN@{\global\options@true\begingroup\optioncodes@}\fi\next@}
\newif\ifN@
\newif\ifE@
\newif\ifNESW@
\newif\ifH@
\newif\ifV@
\newif\ifHshort@
\expandafter\def\csname\string @(\endcsname #1,#2){%
 \ifoptions@\expandafter\endgroup\fi
 \N@false\E@false\H@false\V@false\Hshort@false
 \ifnum#1>\z@\E@true\fi
 \ifnum#1=\z@\V@true\global\tX@false\global\tY@false\global\a@false\fi
 \ifnum#2>\z@\N@true\fi
 \ifnum#2=\z@\H@true\global\tX@false\global\tY@false\global\a@false
  \ifshort@\Hshort@true\fi\fi
 \NESW@false
 \ifN@\ifE@\NESW@true\fi\else\ifE@\else\NESW@true\fi\fi
 \arrow@{#1}{#2}%
 \global\options@false
 \global\scount@\z@\global\tcount@\z@\global\arrcount@\z@
 \global\s@false\global\sxdimen@\z@\global\sydimen@\z@
 \global\tX@false\global\tXdimen@i\z@\global\tXdimen@ii\z@
 \global\tY@false\global\tYdimen@i\z@\global\tYdimen@ii\z@
 \global\a@false\global\exacount@\z@
 \global\x@false\global\xdimen@\z@
 \global\X@false\global\Xdimen@\z@
 \global\y@false\global\ydimen@\z@
 \global\Y@false\global\Ydimen@\z@
 \global\p@false\global\pdimen@\z@
 \global\label@ifalse\global\label@iifalse
 \global\dl@ifalse\global\ldimen@i\z@
 \global\dl@iifalse\global\ldimen@ii\z@
 \global\short@false\global\unshort@false}
\newif\iflabel@i
\newif\iflabel@ii
\newcount\scount@
\newcount\tcount@
\newcount\arrcount@
\newif\ifs@
\newdimen\sxdimen@
\newdimen\sydimen@
\newif\iftX@
\newdimen\tXdimen@i
\newdimen\tXdimen@ii
\newif\iftY@
\newdimen\tYdimen@i
\newdimen\tYdimen@ii
\newif\ifa@
\newcount\exacount@
\newif\ifx@
\newdimen\xdimen@
\newif\ifX@
\newdimen\Xdimen@
\newif\ify@
\newdimen\ydimen@
\newif\ifY@
\newdimen\Ydimen@
\newif\ifp@
\newdimen\pdimen@
\newif\ifdl@i
\newif\ifdl@ii
\newdimen\ldimen@i
\newdimen\ldimen@ii
\newif\ifshort@
\newif\ifunshort@
\def\zero@#1{\ifnum\scount@=\z@
 \if#1e\global\scount@\m@ne\else
 \if#1t\global\scount@\tw@\else
 \if#1h\global\scount@\thr@@\else
 \if#1'\global\scount@6 \else
 \if#1`\global\scount@7 \else
 \if#1(\global\scount@8 \else
 \if#1)\global\scount@9 \else
 \if#1s\global\scount@12 \else
 \if#1H\global\scount@13 \else
 \Err@{\Invalid@@ option \string\0}\fi\fi\fi\fi\fi\fi\fi\fi\fi
 \fi}
\def\one@#1{\ifnum\tcount@=\z@
 \if#1e\global\tcount@\m@ne\else
 \if#1h\global\tcount@\tw@\else
 \if#1t\global\tcount@\thr@@\else
 \if#1'\global\tcount@4 \else
 \if#1`\global\tcount@5 \else
 \if#1(\global\tcount@\ten@ \else
 \if#1)\global\tcount@11 \else
 \if#1s\global\tcount@12 \else
 \if#1H\global\tcount@13 \else
 \Err@{\Invalid@@ option \string\1}\fi\fi\fi\fi\fi\fi\fi\fi\fi
 \fi}
\def\a@#1{\ifnum\arrcount@=\z@
 \if#10\global\arrcount@\m@ne\else
 \if#1+\global\arrcount@\@ne\else
 \if#1-\global\arrcount@\tw@\else
 \if#1=\global\arrcount@\thr@@\else
 \Err@{\Invalid@@ option \string\a}\fi\fi\fi\fi
 \fi}
\def\ds@{\ifnum\catcode`\;=\active\expandafter\dsA@\else
 \expandafter\dsO@\fi}
\def\dsO@(#1;#2){\ds@@{#1}{#2}}
\def\ds@@#1#2{\ifs@\else
 \global\s@true
 \global\sxdimen@\hunit\global\sxdimen@#1\sxdimen@\relax
 \global\sydimen@\vunit\global\sydimen@#2\sydimen@\relax
 \fi}
\def\dtX@{\ifnum\catcode`\;=\active\expandafter\dtXA@\else
 \expandafter\dtXO@\fi}
\def\dtXO@(#1;#2){\dtX@@{#1}{#2}}
\def\dtX@@#1#2{\iftX@\else
 \global\tX@true
 \global\tXdimen@i\hunit\global\tXdimen@i#1\tXdimen@i\relax
 \global\tXdimen@ii\vunit\global\tXdimen@ii#2\tXdimen@ii\relax
 \fi}
\def\dtY@{\ifnum\catcode`\;=\active\expandafter\dtYA@\else
 \expandafter\dtYO@\fi}
\def\dtYO@(#1;#2){\dtY@@{#1}{#2}}
\def\dtY@@#1#2{\iftY@\else
 \global\tY@true
 \global\tYdimen@i\hunit\global\tYdimen@i#1\tYdimen@i\relax
 \global\tYdimen@ii\vunit\global\tYdimen@ii#2\tYdimen@ii\relax
 \fi}
{\catcode`\;=\active
 \gdef\dsA@(#1;#2){\ds@@{#1}{#2}}
 \gdef\dtXA@(#1;#2){\dtX@@{#1}{#2}}
 \gdef\dtYA@(#1;#2){\dtY@@{#1}{#2}}
}
\def\da@#1{\ifa@\else\global\a@true\global\exacount@#1\relax\fi}
\def\dx@#1{\ifx@\else
 \global\x@true
 \global\xdimen@\hunit\global\xdimen@#1\xdimen@\relax
 \fi}
\def\dX@#1{\ifX@\else
 \global\X@true
 \global\Xdimen@\hunit\global\Xdimen@#1\Xdimen@\relax
 \fi}
\def\dy@#1{\ify@\else
 \global\y@true
 \global\ydimen@\vunit\global\ydimen@#1\ydimen@\relax
 \fi}
\def\dY@#1{\ifY@\else
 \global\Y@true
 \global\Ydimen@\vunit\global\Ydimen@#1\Ydimen@\relax
 \fi}
\def\p@@#1{\ifp@\else
 \global\p@true
 \global\pdimen@\hunit\global\divide\pdimen@\tw@
 \global\pdimen@#1\pdimen@\relax
 \fi}
\def\L@#1{\iflabel@i\else
 \global\label@itrue\gdef\label@i{#1}%
 \fi}
\def\l@#1{\iflabel@ii\else
 \global\label@iitrue\gdef\label@ii{#1}%
 \fi}
\def\dL@#1{\ifdl@i\else
 \global\dl@itrue\global\ldimen@i\hunit\global\ldimen@i#1\ldimen@i\relax
 \fi}
\def\dl@#1{\ifdl@ii\else
 \global\dl@iitrue\global\ldimen@ii\hunit\global\ldimen@ii#1\ldimen@ii\relax
 \fi}
\def\s@{\ifunshort@\else\global\short@true\fi}
\def\uns@{\ifshort@\else\global\unshort@true\global\short@false\fi}
\def\optioncodes@{\let\0\zero@\let\1\one@\let\a\a@\let\ds\ds@\let\dtX\dtX@
 \let\dtY\dtY@\let\da\da@\let\dx\dx@\let\dX\dX@\let\dY\dY@\let\dy\dy@
 \let\p\p@@\let\L\L@\let\l\l@\let\dL\dL@\let\dl\dl@\let\s\s@\let\uns\uns@}
\def\slopes@{\\161\\152\\143\\134\\255\\126\\357\\238\\349\\45{10}\\56{11}%
 \\11{12}\\65{13}\\54{14}\\43{15}\\32{16}\\53{17}\\21{18}\\52{19}\\31{20}%
 \\41{21}\\51{22}\\61{23}}
\newcount\tan@i
\newcount\tan@ip
\newcount\tan@ii
\newcount\tan@iip
\newdimen\slope@i
\newdimen\slope@ip
\newdimen\slope@ii
\newdimen\slope@iip
\newcount\angcount@
\newcount\extracount@
\def\slope@{{\slope@i\secondy@\advance\slope@i-\firsty@
 \ifN@\else\multiply\slope@i\m@ne\fi
 \slope@ii\secondx@\advance\slope@ii-\firstx@
 \ifE@\else\multiply\slope@ii\m@ne\fi
 \ifdim\slope@ii<\z@
  \global\tan@i6 \global\tan@ii\@ne\global\angcount@23
 \else
  \dimen@\slope@i\multiply\dimen@6
  \ifdim\dimen@<\slope@ii
   \global\tan@i\@ne\global\tan@ii6 \global\angcount@\@ne
  \else
   \dimen@\slope@ii\multiply\dimen@6
   \ifdim\dimen@<\slope@i
    \global\tan@i6 \global\tan@ii\@ne\global\angcount@23
   \else
    \global\tan@ip\z@\global\tan@iip\@ne
    \def\\##1##2##3{\global\angcount@##3\relax
     \slope@ip\slope@i\slope@iip\slope@ii
     \multiply\slope@iip##1\relax\multiply\slope@ip##2\relax
     \ifdim\slope@iip<\slope@ip
      \global\tan@ip##1\relax\global\tan@iip##2\relax
     \else
      \global\tan@i##1\relax\global\tan@ii##2\relax
      \def\\####1####2####3{}%
     \fi}%
    \slopes@
    \slope@i\secondy@\advance\slope@i-\firsty@
    \ifN@\else\multiply\slope@i\m@ne\fi
    \multiply\slope@i\tan@ii\multiply\slope@i\tan@iip\multiply\slope@i\tw@
    \count@\tan@i\multiply\count@\tan@iip
    \extracount@\tan@ip\multiply\extracount@\tan@ii
    \advance\count@\extracount@
    \slope@ii\secondx@\advance\slope@ii-\firstx@
    \ifE@\else\multiply\slope@ii\m@ne\fi
    \multiply\slope@ii\count@
    \ifdim\slope@i<\slope@ii
     \global\tan@i\tan@ip\global\tan@ii\tan@iip
     \global\advance\angcount@\m@ne
    \fi
   \fi
  \fi
 \fi}%
}
\def\slope@a#1{{\def\\##1##2##3{\ifnum##3=#1\global\tan@i##1\relax
 \global\tan@ii##2\relax\fi}\slopes@}}
\newcount\i@
\newcount\j@
\newcount\colcount@
\newcount\Colcount@
\newcount\tcolcount@
\newdimen\rowht@
\newdimen\rowdp@
\newcount\rowcount@
\newcount\Rowcount@
\newcount\maxcolrow@
\newtoks\colwidthtoks@
\newtoks\Rowheighttoks@
\newtoks\Rowdepthtoks@
\newtoks\widthtoks@
\newtoks\Widthtoks@
\newtoks\heighttoks@
\newtoks\Heighttoks@
\newtoks\depthtoks@
\newtoks\Depthtoks@
\newif\iffirstCDcr@
\def\dotoks@i{%
 \global\widthtoks@\expandafter{\the\widthtoks@\else\getdim@\z@\fi}%
 \global\heighttoks@\expandafter{\the\heighttoks@\else\getdim@\z@\fi}%
 \global\depthtoks@\expandafter{\the\depthtoks@\else\getdim@\z@\fi}}
\def\dotoks@ii{%
 \global\widthtoks@{\ifcase\j@}%
 \global\heighttoks@{\ifcase\j@}%
 \global\depthtoks@{\ifcase\j@}}
\def\preCD@#1\endCD{\setbox\ZER@
 \vbox{%
  \def\arrow@##1##2{{}}%
  \global\rowcount@\m@ne\global\colcount@\z@\global\Colcount@\z@
  \global\firstCDcr@true\toks@{}%
  \global\widthtoks@{\ifcase\j@}%
  \global\Widthtoks@{\ifcase\i@}%
  \global\heighttoks@{\ifcase\j@}%
  \global\Heighttoks@{\ifcase\i@}%
  \global\depthtoks@{\ifcase\j@}%
  \global\Depthtoks@{\ifcase\i@}%
  \global\Rowheighttoks@{\ifcase\i@}%
  \global\Rowdepthtoks@{\ifcase\i@}%
  \Let@
  \everycr{%
   \noalign{%
    \global\advance\rowcount@\@ne
    \ifnum\colcount@<\Colcount@
    \else
     \global\Colcount@\colcount@\global\maxcolrow@\rowcount@
    \fi
    \global\colcount@\z@
    \iffirstCDcr@
     \global\firstCDcr@false
    \else
     \edef\next@{\the\Rowheighttoks@\noexpand\or\noexpand\getdim@\the\rowht@}%
      \global\Rowheighttoks@\expandafter{\next@}%
     \edef\next@{\the\Rowdepthtoks@\noexpand\or\noexpand\getdim@\the\rowdp@}%
      \global\Rowdepthtoks@\expandafter{\next@}%
     \global\rowht@\z@\global\rowdp@\z@
     \dotoks@i
     \edef\next@{\the\Widthtoks@\noexpand\or\the\widthtoks@}%
      \global\Widthtoks@\expandafter{\next@}%
     \edef\next@{\the\Heighttoks@\noexpand\or\the\heighttoks@}%
      \global\Heighttoks@\expandafter{\next@}%
     \edef\next@{\the\Depthtoks@\noexpand\or\the\depthtoks@}%
      \global\Depthtoks@\expandafter{\next@}%
     \dotoks@ii
    \fi}}%
  \tabskip\z@
  \halign{&\setbox\ZER@\hbox{\vrule\height\ten@\p@\width\z@\depth\z@     
   $\m@th\displaystyle{##}$}\copy\ZER@
   \ifdim\ht\ZER@>\rowht@\global\rowht@\ht\ZER@\fi
   \ifdim\dp\ZER@>\rowdp@\global\rowdp@\dp\ZER@\fi
   \global\advance\colcount@\@ne
   \edef\next@{\the\widthtoks@\noexpand\or\noexpand\getdim@\the\wd\ZER@}%
    \global\widthtoks@\expandafter{\next@}%
   \edef\next@{\the\heighttoks@\noexpand\or\noexpand\getdim@\the\ht\ZER@}%
    \global\heighttoks@\expandafter{\next@}%
   \edef\next@{\the\depthtoks@\noexpand\or\noexpand\getdim@\the\dp\ZER@}%
    \global\depthtoks@\expandafter{\next@}%
   \cr#1\crcr}}%
 \Rowcount@\rowcount@
 \global\Widthtoks@\expandafter{\the\Widthtoks@\fi\relax}%
 \edef\Width@##1##2{\i@##1\relax\j@##2\relax\the\Widthtoks@}%
 \global\Heighttoks@\expandafter{\the\Heighttoks@\fi\relax}%
 \edef\Height@##1##2{\i@##1\relax\j@##2\relax\the\Heighttoks@}%
 \global\Depthtoks@\expandafter{\the\Depthtoks@\fi\relax}%
 \edef\Depth@##1##2{\i@##1\relax\j@##2\relax\the\Depthtoks@}%
 \edef\next@{\the\Rowheighttoks@\noexpand\fi\relax}%
 \global\Rowheighttoks@\expandafter{\next@}%
 \edef\Rowheight@##1{\i@##1\relax\the\Rowheighttoks@}%
 \edef\next@{\the\Rowdepthtoks@\noexpand\fi\relax}%
 \global\Rowdepthtoks@\expandafter{\next@}%
 \edef\Rowdepth@##1{\i@##1\relax\the\Rowdepthtoks@}%
 \global\colwidthtoks@{\fi}%
 \setbox\ZER@\vbox{%
  \unvbox\ZER@
  \count@\rowcount@
  \loop
   \unskip\unpenalty
   \setbox\ZER@\lastbox
   \ifnum\count@>\maxcolrow@\advance\count@\m@ne
   \repeat
  \hbox{%
   \unhbox\ZER@
   \count@\z@
   \loop
    \unskip
    \setbox\ZER@\lastbox
    \edef\next@{\noexpand\or\noexpand\getdim@\the\wd\ZER@\the\colwidthtoks@}%
     \global\colwidthtoks@\expandafter{\next@}%
    \advance\count@\@ne
    \ifnum\count@<\Colcount@
    \repeat}}%
 \edef\next@{\noexpand\ifcase\noexpand\i@\the\colwidthtoks@}%
  \global\colwidthtoks@\expandafter{\next@}%
 \edef\Colwidth@##1{\i@##1\relax\the\colwidthtoks@}%
 \global\colwidthtoks@{}\global\Rowheighttoks@{}\global\Rowdepthtoks@{}%
 \global\widthtoks@{}\global\Widthtoks@{}\global\heighttoks@{}%
 \global\Heighttoks@{}\global\depthtoks@{}\global\Depthtoks@{}%
}
\newcount\xoff@
\newcount\yoff@
\newcount\endcount@
\newcount\rcount@
\newdimen\firstx@
\newdimen\firsty@
\newdimen\secondx@
\newdimen\secondy@
\newdimen\tocenter@
\newdimen\charht@
\newdimen\charwd@
\def\outside@{\Err@{This arrow points outside the \string\CD}}
\newif\ifsvertex@
\newif\iftvertex@
\def\arrow@#1#2{\global\xoff@#1\relax\global\yoff@#2\relax
 \count@\rowcount@\advance\count@-\yoff@
 \ifnum\count@<\@ne\outside@\else\ifnum\count@>\Rowcount@\outside@\fi\fi
 \count@\colcount@\advance\count@\xoff@
 \ifnum\count@<\@ne\outside@\else\ifnum\count@>\Colcount@\outside@\fi\fi
 \tcolcount@\colcount@\advance\tcolcount@\xoff@
 \Width@\rowcount@\colcount@\divide\getdim@\tw@\tocenter@-\getdim@
 \ifdim\getdim@=\z@
  \firstx@\z@\firsty@\mathaxis@\svertex@true
 \else
  \svertex@false
  \ifHshort@
   \Colwidth@\colcount@\divide\getdim@\tw@
   \ifE@ \firstx@\getdim@ \else \firstx@-\getdim@ \fi
  \else
   \ifE@ \firstx@\getdim@ \else \firstx@-\getdim@ \fi
  \fi
  \ifE@
   \ifH@ \advance\firstx@\thr@@\p@ \else \advance\firstx@-\thr@@\p@ \fi  
  \else
   \ifH@ \advance\firstx@-\thr@@\p@ \else \advance\firstx@\thr@@\p@ \fi  
  \fi
  \ifN@
   \Height@\rowcount@\colcount@ \firsty@\getdim@                         
   \ifV@ \advance\firsty@\thr@@\p@ \fi                                   
  \else
   \ifV@
    \Depth@\rowcount@\colcount@ \firsty@-\getdim@                        
    \advance\firsty@-\thr@@\p@                                           
   \else
    \firsty@\z@                                                          
   \fi
  \fi
 \fi
 \ifV@
 \else
  \Colwidth@\colcount@\divide\getdim@\tw@
  \ifE@\secondx@\getdim@\else\secondx@-\getdim@\fi
  \ifE@\else\getcgap@\colcount@\advance\secondx@-\getdim@\fi
  \endcount@\colcount@\advance\endcount@\xoff@
  \count@\colcount@
  \ifE@
   \advance\count@\@ne
   \loop
    \ifnum\count@<\endcount@
    \Colwidth@\count@\advance\secondx@\getdim@
    \getcgap@\count@\advance\secondx@\getdim@
    \advance\count@\@ne
    \repeat
  \else
   \advance\count@\m@ne
   \loop
    \ifnum\count@>\endcount@
    \Colwidth@\count@\advance\secondx@-\getdim@
    \getcgap@\count@\advance\secondx@-\getdim@
    \advance\count@\m@ne
    \repeat
  \fi
  \Colwidth@\count@\divide\getdim@\tw@
  \ifHshort@
  \else
   \ifE@\advance\secondx@\getdim@\else\advance\secondx@-\getdim@\fi
  \fi
  \ifE@\getcgap@\count@\advance\secondx@\getdim@\fi
  \rcount@\rowcount@\advance\rcount@-\yoff@
  \Width@\rcount@\count@\divide\getdim@\tw@
  \tvertex@false
  \ifH@\ifdim\getdim@=\z@\tvertex@true\Hshort@false\fi\fi
  \ifHshort@
  \else
   \ifE@\advance\secondx@-\getdim@\else\advance\secondx@\getdim@\fi
  \fi
  \iftvertex@
   \advance\secondx@.4\p@
  \else
   \ifE@\advance\secondx@-\thr@@\p@\else\advance\secondx@\thr@@\p@\fi    
  \fi
 \fi
 \ifH@
 \else
  \ifN@
   \Rowheight@\rowcount@\secondy@\getdim@
  \else
   \Rowdepth@\rowcount@\secondy@-\getdim@
   \getrgap@\rowcount@\advance\secondy@-\getdim@
  \fi
  \endcount@\rowcount@\advance\endcount@-\yoff@
  \count@\rowcount@
  \ifN@
   \advance\count@\m@ne
   \loop
    \ifnum\count@>\endcount@
    \Rowheight@\count@\advance\secondy@\getdim@
    \Rowdepth@\count@\advance\secondy@\getdim@
    \getrgap@\count@\advance\secondy@\getdim@
    \advance\count@\m@ne
    \repeat
  \else
   \advance\count@\@ne
   \loop
    \ifnum\count@<\endcount@
    \Rowheight@\count@\advance\secondy@-\getdim@
    \Rowdepth@\count@\advance\secondy@-\getdim@
    \getrgap@\count@\advance\secondy@-\getdim@
    \advance\count@\@ne
    \repeat
  \fi
  \tvertex@false
  \ifV@\Width@\count@\colcount@\ifdim\getdim@=\z@\tvertex@true\fi\fi
  \ifN@
   \getrgap@\count@\advance\secondy@\getdim@
   \Rowdepth@\count@\advance\secondy@\getdim@
   \iftvertex@
    \advance\secondy@\mathaxis@
   \else
    \Depth@\count@\tcolcount@\advance\secondy@-\getdim@
    \advance\secondy@-\thr@@\p@                                          
   \fi
  \else
   \Rowheight@\count@\advance\secondy@-\getdim@
   \iftvertex@
    \advance\secondy@\mathaxis@
   \else
    \Height@\count@\tcolcount@\advance\secondy@\getdim@
    \advance\secondy@\thr@@\p@                                           
   \fi
  \fi
 \fi
 \ifV@\else\advance\firstx@\sxdimen@\fi
 \ifH@\else\advance\firsty@\sydimen@\fi
 \iftX@
  \advance\secondy@\tXdimen@ii
  \advance\secondx@\tXdimen@i
  \slope@
 \else
  \iftY@
   \advance\secondy@\tYdimen@ii
   \advance\secondx@\tYdimen@i
   \slope@
   \secondy@\secondx@\advance\secondy@-\firstx@
   \ifNESW@\else\multiply\secondy@\m@ne\fi
   \multiply\secondy@\tan@i\divide\secondy@\tan@ii\advance\secondy@\firsty@
  \else
   \ifa@
    \slope@
    \ifNESW@\global\advance\angcount@\exacount@\else
     \global\advance\angcount@-\exacount@\fi
    \ifnum\angcount@>23 \global\angcount@23 \fi
    \ifnum\angcount@<\@ne\global\angcount@\@ne\fi
    \slope@a\angcount@
    \ifY@
     \advance\secondy@\Ydimen@
    \else
     \ifX@
      \advance\secondx@\Xdimen@
      \dimen@\secondx@\advance\dimen@-\firstx@
      \ifNESW@\else\multiply\dimen@\m@ne\fi
      \multiply\dimen@\tan@i\divide\dimen@\tan@ii
      \advance\dimen@\firsty@\secondy@\dimen@
     \fi
    \fi
   \else
    \ifH@\else\ifV@\else\slope@\fi\fi
   \fi
  \fi
 \fi
 \ifH@\else\ifV@\else\ifsvertex@\else
  \dimen@6\p@\multiply\dimen@\tan@ii
  \count@\tan@i\advance\count@\tan@ii\divide\dimen@\count@
  \ifE@\advance\firstx@\dimen@\else\advance\firstx@-\dimen@\fi
  \multiply\dimen@\tan@i\divide\dimen@\tan@ii
  \ifN@\advance\firsty@\dimen@\else\advance\firsty@-\dimen@\fi
 \fi\fi\fi
 \ifp@
  \ifH@\else\ifV@\else
   \getcos@\pdimen@\advance\firsty@\dimen@\advance\secondy@\dimen@
   \ifNESW@\advance\firstx@-\dimen@ii\else\advance\firstx@\dimen@ii\fi
  \fi\fi
 \fi
 \ifH@\else\ifV@\else
  \ifnum\tan@i>\tan@ii
   \charht@\ten@\p@\charwd@\ten@\p@
   \multiply\charwd@\tan@ii\divide\charwd@\tan@i
  \else
   \charwd@\ten@\p@\charht@\ten@\p@
   \divide\charht@\tan@ii\multiply\charht@\tan@i
  \fi
  \ifnum\tcount@=\thr@@
   \ifN@\advance\secondy@-.3\charht@\else\advance\secondy@.3\charht@\fi
  \fi
  \ifnum\scount@=\tw@
   \ifE@\advance\firstx@.3\charht@\else\advance\firstx@-.3\charht@\fi
  \fi
  \ifnum\tcount@=12
   \ifN@\advance\secondy@-\charht@\else\advance\secondy@\charht@\fi
  \fi
  \iftY@
  \else
   \ifa@
    \ifX@
    \else
     \secondx@\secondy@\advance\secondx@-\firsty@
     \ifNESW@\else\multiply\secondx@\m@ne\fi
     \multiply\secondx@\tan@ii\divide\secondx@\tan@i
     \advance\secondx@\firstx@
    \fi
   \fi
  \fi
 \fi\fi
 \ifH@\harrow@\else\ifV@\varrow@\else\arrow@@\fi\fi}
\newdimen\mathaxis@
\mathaxis@90\p@\divide\mathaxis@36
\def\harrow@b{\ifE@\hskip\tocenter@\hskip\firstx@\fi}
\def\harrow@bb{\ifE@\hskip\xdimen@\else\hskip\Xdimen@\fi}
\def\harrow@e{\ifE@\else\hskip-\firstx@\hskip-\tocenter@\fi}
\def\harrow@ee{\ifE@\hskip-\Xdimen@\else\hskip-\xdimen@\fi}
\def\harrow@{\dimen@\secondx@\advance\dimen@-\firstx@
 \ifE@\let\next@\rlap\else\multiply\dimen@\m@ne\let\next@\llap\fi
 \next@{%
  \harrow@b
  \smash{\raise\pdimen@\hbox to\dimen@
   {\harrow@bb\arrow@ii
    \ifnum\arrcount@=\m@ne\else\ifnum\arrcount@=\thr@@\else
     \ifE@
      \ifnum\scount@=\m@ne
      \else
       \ifcase\scount@\or\or\char118 \or\char117 \or\or\or\char119 \or
       \char120 \or\char121 \or\char122 \or\or\or\arrow@i\char125 \or
       \char117 \hskip\thr@@\p@\char117 \hskip-\thr@@\p@\fi
      \fi
     \else
      \ifnum\tcount@=\m@ne
      \else
       \ifcase\tcount@\char117 \or\or\char117 \or\char118 \or\char119 \or
       \char120 \or\or\or\or\or\char121 \or\char122 \or\arrow@i\char125
       \or\char117 \hskip\thr@@\p@\char117 \hskip-\thr@@\p@\fi
      \fi
     \fi
    \fi\fi
    \dimen@\mathaxis@\advance\dimen@.2\p@
    \dimen@ii\mathaxis@\advance\dimen@ii-.2\p@
    \ifnum\arrcount@=\m@ne
     \let\leads@\null
    \else
     \ifcase\arrcount@
      \def\leads@{\hrule\height\dimen@\depth-\dimen@ii}\or
      \def\leads@{\hrule\height\dimen@\depth-\dimen@ii}\or
      \def\leads@{\hbox to\ten@\p@{%
       \leaders\hrule\height\dimen@\depth-\dimen@ii\hfil
       \hfil
      \leaders\hrule\height\dimen@\depth-\dimen@ii\hskip\z@ plus2fil\relax
       \hfil
       \leaders\hrule\height\dimen@\depth-\dimen@ii\hfil}}\or
     \def\leads@{\hbox{\hbox to\ten@\p@{\dimen@\mathaxis@\advance\dimen@1.2\p@
       \dimen@ii\dimen@\advance\dimen@ii-.4\p@
       \leaders\hrule\height\dimen@\depth-\dimen@ii\hfil}%
       \kern-\ten@\p@
       \hbox to\ten@\p@{\dimen@\mathaxis@\advance\dimen@-1.2\p@
       \dimen@ii\dimen@\advance\dimen@ii-.4\p@
       \leaders\hrule\height\dimen@\depth-\dimen@ii\hfil}}}\fi
    \fi
    \cleaders\leads@\hfil
    \ifnum\arrcount@=\m@ne\else\ifnum\arrcount@=\thr@@\else
     \arrow@i
     \ifE@
      \ifnum\tcount@=\m@ne
      \else
       \ifcase\tcount@\char119 \or\or\char119 \or\char120 \or\char121 \or
       \char122 \or \or\or\or\or\char123 \or\char124 \or
       \char125 \or\char119 \hskip-\thr@@\p@\char119 \hskip\thr@@\p@\fi
      \fi
     \else
      \ifcase\scount@\or\or\char120 \or\char119 \or\or\or\char121 \or\char122
      \or\char123 \or\char124 \or\or\or\char125 \or
      \char119 \hskip-\thr@@\p@\char119 \hskip\thr@@\p@\fi
     \fi
    \fi\fi
    \harrow@ee}}%
  \harrow@e}%
 \iflabel@i
  \dimen@ii\z@\setbox\ZER@\hbox{$\m@th\tsize@@\label@i$}%
  \ifnum\arrcount@=\m@ne
  \else
   \advance\dimen@ii\mathaxis@
   \advance\dimen@ii\dp\ZER@\advance\dimen@ii\tw@\p@
   \ifnum\arrcount@=\thr@@\advance\dimen@ii\tw@\p@\fi
  \fi
  \advance\dimen@ii\pdimen@
  \next@{\harrow@b\smash{\raise\dimen@ii\hbox to\dimen@
   {\harrow@bb\hskip\tw@\ldimen@i\hfil\box\ZER@\hfil\harrow@ee}}\harrow@e}%
 \fi
 \iflabel@ii
  \ifnum\arrcount@=\m@ne
  \else
   \setbox\ZER@\hbox{$\m@th\tsize@\label@ii$}%
   \dimen@ii-\ht\ZER@\advance\dimen@ii-\tw@\p@
   \ifnum\arrcount@=\thr@@\advance\dimen@ii-\tw@\p@\fi
   \advance\dimen@ii\mathaxis@\advance\dimen@ii\pdimen@
   \next@{\harrow@b\smash{\raise\dimen@ii\hbox to\dimen@
    {\harrow@bb\hskip\tw@\ldimen@ii\hfil\box\ZER@\hfil\harrow@ee}}\harrow@e}%
  \fi
 \fi}
\let\tsize@\tsize
\def\tsizeCDlabels{\let\tsize@\tsize}
\def\ssizeCDlabels{\let\tsize@\ssize}
\def\tsize@@{\ifnum\arrcount@=\m@ne\else\tsize@\fi}
\def\varrow@{\dimen@\secondy@\advance\dimen@-\firsty@
 \ifN@\else\multiply\dimen@\m@ne\fi
 \setbox\ZER@\vbox to\dimen@
  {\ifN@\vskip-\Ydimen@\else\vskip\ydimen@\fi
   \ifnum\arrcount@=\m@ne\else\ifnum\arrcount@=\thr@@\else
    \hbox{\arrow@iii
     \ifN@
      \ifnum\tcount@=\m@ne
      \else
       \ifcase\tcount@\char117 \or\or\char117 \or\char118 \or\char119 \or
       \char120 \or\or\or\or\or\char121 \or\char122 \or\char123 \or
       \vbox{\hbox{\char117}\nointerlineskip\vskip\thr@@\p@
       \hbox{\char117}\vskip-\thr@@\p@}\fi
      \fi
     \else
      \ifcase\scount@\or\or\char118 \or\char117 \or\or\or\char119 \or
      \char120 \or\char121 \or\char122 \or\or\or\char123 \or
      \vbox{\hbox{\char117}\nointerlineskip\vskip\thr@@\p@
      \hbox{\char117}\vskip-\thr@@\p@}\fi
     \fi}%
    \nointerlineskip
   \fi\fi
   \ifnum\arrcount@=\m@ne
    \let\leads@\null
   \else
    \ifcase\arrcount@\let\leads@\vrule\or\let\leads@\vrule\or
    \def\leads@{\vbox to\ten@\p@{%
     \hrule\height1.67\p@\depth\z@\width.4\p@
     \vfil
     \hrule\height3.33\p@\depth\z@\width.4\p@
     \vfil
     \hrule\height1.67\p@\depth\z@\width.4\p@}}\or
    \def\leads@{\hbox{\vrule\height\p@\hskip\tw@\p@\vrule}}\fi
   \fi
  \cleaders\leads@\vfill\nointerlineskip
   \ifnum\arrcount@=\m@ne\else\ifnum\arrcount@=\thr@@\else
    \hbox{\arrow@iv
     \ifN@
      \ifcase\scount@\or\or\char118 \or\char117 \or\or\or\char119 \or
      \char120 \or\char121 \or\char122 \or\or\or\arrow@iii\char123 \or
      \vbox{\hbox{\char117}\nointerlineskip\vskip-\thr@@\p@
      \hbox{\char117}\vskip\thr@@\p@}\fi
     \else
      \ifnum\tcount@=\m@ne
      \else
       \ifcase\tcount@\char117 \or\or\char117 \or\char118 \or\char119 \or
       \char120 \or\or\or\or\or\char121 \or\char122 \or\arrow@iii\char123 \or
       \vbox{\hbox{\char117}\nointerlineskip\vskip-\thr@@\p@
       \hbox{\char117}\vskip\thr@@\p@}\fi
      \fi
     \fi}%
   \fi\fi
   \ifN@\vskip\ydimen@\else\vskip-\Ydimen@\fi}%
 \ifN@
  \dimen@ii\firsty@
 \else
  \dimen@ii-\firsty@\advance\dimen@ii\ht\ZER@\multiply\dimen@ii\m@ne
 \fi
 \rlap{\smash{\hskip\tocenter@\hskip\pdimen@\raise\dimen@ii\box\ZER@}}%
 \iflabel@i
  \setbox\ZER@\vbox to\dimen@{\vfil
   \hbox{$\m@th\tsize@@\label@i$}\vskip\tw@\ldimen@i\vfil}%
  \rlap{\smash{\hskip\tocenter@\hskip\pdimen@
  \ifnum\arrcount@=\m@ne\let\next@\relax\else\let\next@\llap\fi
  \next@{\raise\dimen@ii\hbox{\ifnum\arrcount@=\m@ne\hskip-.5\wd\ZER@\fi
   \box\ZER@\ifnum\arrcount@=\m@ne\else\hskip\tw@\p@\fi}}}}%
 \fi
 \iflabel@ii
  \ifnum\arrcount@=\m@ne
  \else
   \setbox\ZER@\vbox to\dimen@{\vfil
    \hbox{$\m@th\tsize@\label@ii$}\vskip\tw@\ldimen@ii\vfil}%
   \rlap{\smash{\hskip\tocenter@\hskip\pdimen@
   \rlap{\raise\dimen@ii\hbox{\ifnum\arrcount@=\thr@@\hskip4.5\p@\else
    \hskip2.5\p@\fi\box\ZER@}}}}%
  \fi
 \fi
}
\newdimen\goal@
\newdimen\shifted@
\newcount\Tcount@
\newcount\Scount@
\newbox\shaft@
\newcount\slcount@
\def\getcos@#1{%
 \ifnum\tan@i<\tan@ii
  \dimen@#1%
  \ifnum\slcount@<8 \count@9 \else \ifnum\slcount@<12 \count@8 \else
   \count@7 \fi\fi
  \multiply\dimen@\count@\divide\dimen@\ten@
  \dimen@ii\dimen@\multiply\dimen@ii\tan@i\divide\dimen@ii\tan@ii
 \else
  \dimen@ii#1%
  \count@-\slcount@\advance\count@24
  \ifnum\count@<8 \count@9 \else \ifnum\count@<12 \count@8
   \else\count@7 \fi\fi
  \multiply\dimen@ii\count@\divide\dimen@ii\ten@
  \dimen@\dimen@ii\multiply\dimen@\tan@ii\divide\dimen@\tan@i
 \fi}
\newdimen\adjust@
\def\Nnext@{\ifN@\let\next@\raise\else\let\next@\lower\fi}
\def\arrow@@{\slcount@\angcount@
 \ifNESW@
  \ifnum\angcount@<\ten@
   \let\arrowfont@\arrow@i\global\advance\angcount@\m@ne
   \global\multiply\angcount@13
  \else
   \ifnum\angcount@<19
    \let\arrowfont@\arrow@ii\global\advance\angcount@-\ten@
    \global\multiply\angcount@13
   \else
    \let\arrowfont@\arrow@iii\global\advance\angcount@-19
    \global\multiply\angcount@13
  \fi\fi
  \Tcount@\angcount@
 \else
  \ifnum\angcount@<5
   \let\arrowfont@\arrow@iii\global\advance\angcount@\m@ne
   \global\multiply\angcount@13 \global\advance\angcount@65
  \else
   \ifnum\angcount@<14
    \let\arrowfont@\arrow@iv\global\advance\angcount@-5
    \global\multiply\angcount@13
   \else
    \ifnum\angcount@<23
     \let\arrowfont@\arrow@v\global\advance\angcount@-14
     \global\multiply\angcount@13
    \else
     \let\arrowfont@\arrow@i\global\angcount@117
  \fi\fi\fi
  \ifnum\angcount@=117 \Tcount@115 \else\Tcount@\angcount@\fi
 \fi
 \Scount@\Tcount@
 \ifE@
  \ifnum\tcount@=\z@\advance\Tcount@\tw@\else\ifnum\tcount@=13
   \advance\Tcount@\tw@\else\advance\Tcount@\tcount@\fi\fi
  \ifnum\scount@=\z@\else\ifnum\scount@=13 \advance\Scount@\thr@@\else
   \advance\Scount@\scount@\fi\fi
 \else
  \ifcase\tcount@\advance\Tcount@\thr@@\or\or\advance\Tcount@\thr@@\or
  \advance\Tcount@\tw@\or\advance\Tcount@6 \or\advance\Tcount@7
  \or\or\or\or\or\advance\Tcount@8 \or\advance\Tcount@9 \or
  \advance\Tcount@12 \or\advance\Tcount@\thr@@\fi
  \ifcase\scount@\or\or\advance\Scount@\thr@@\or\advance\Scount@\tw@\or
  \or\or\advance\Scount@4 \or\advance\Scount@5 \or\advance\Scount@\ten@
  \or\advance\Scount@11 \or\or\or\advance\Scount@12 \or\advance
  \Scount@\tw@\fi
 \fi
 \ifcase\arrcount@\or\or\global\advance\angcount@\@ne\else\fi
 \ifN@\shifted@\firsty@\else\shifted@-\firsty@\fi
 \ifE@\else\advance\shifted@\charht@\fi
 \goal@\secondy@\advance\goal@-\firsty@
 \ifN@\else\multiply\goal@\m@ne\fi
 \setbox\shaft@\hbox{\arrowfont@\char\angcount@}%
 \ifnum\arrcount@=\thr@@
  \getcos@{1.5\p@}%
  \setbox\shaft@\hbox to\wd\shaft@{\arrowfont@
   \rlap{\hskip\dimen@ii
    \smash{\ifNESW@\let\next@\lower\else\let\next@\raise\fi
     \next@\dimen@\hbox{\arrowfont@\char\angcount@}}}%
   \rlap{\hskip-\dimen@ii
    \smash{\ifNESW@\let\next@\raise\else\let\next@\lower\fi
      \next@\dimen@\hbox{\arrowfont@\char\angcount@}}}\hfil}%
 \fi
 \rlap{\smash{\hskip\tocenter@\hskip\firstx@
  \ifnum\arrcount@=\m@ne
  \else
   \ifnum\arrcount@=\thr@@
   \else
    \ifnum\scount@=\m@ne
    \else
     \ifnum\scount@=\z@
     \else
      \setbox\ZER@\hbox{\ifnum\angcount@=117 \arrow@v\else\arrowfont@\fi
       \char\Scount@}%
      \ifNESW@
       \ifnum\scount@=\tw@
        \dimen@\shifted@\advance\dimen@-\charht@
        \ifN@\hskip-\wd\ZER@\fi
        \Nnext@
        \next@\dimen@\copy\ZER@
        \ifN@\else\hskip-\wd\ZER@\fi
       \else
        \Nnext@
        \ifN@\else\hskip-\wd\ZER@\fi
        \next@\shifted@\copy\ZER@
        \ifN@\hskip-\wd\ZER@\fi
       \fi
       \ifnum\scount@=12
        \advance\shifted@\charht@\advance\goal@-\charht@
        \ifN@\hskip\wd\ZER@\else\hskip-\wd\ZER@\fi
       \fi
       \ifnum\scount@=13
        \getcos@{\thr@@\p@}%
        \ifN@\hskip\dimen@\else\hskip-\wd\ZER@\hskip-\dimen@\fi
        \adjust@\shifted@\advance\adjust@\dimen@ii
        \Nnext@
        \next@\adjust@\copy\ZER@
        \ifN@\hskip-\dimen@\hskip-\wd\ZER@\else\hskip\dimen@\fi
       \fi
      \else
       \ifN@\hskip-\wd\ZER@\fi
       \ifnum\scount@=\tw@
        \ifN@\hskip\wd\ZER@\else\hskip-\wd\ZER@\fi
        \dimen@\shifted@\advance\dimen@-\charht@
        \Nnext@
        \next@\dimen@\copy\ZER@
        \ifN@\hskip-\wd\ZER@\fi
       \else
        \Nnext@
        \next@\shifted@\copy\ZER@
        \ifN@\else\hskip-\wd\ZER@\fi
       \fi
       \ifnum\scount@=12
        \advance\shifted@\charht@\advance\goal@-\charht@
        \ifN@\hskip-\wd\ZER@\else\hskip\wd\ZER@\fi
       \fi
       \ifnum\scount@=13
        \getcos@{\thr@@\p@}%
        \ifN@\hskip-\wd\ZER@\hskip-\dimen@\else\hskip\dimen@\fi
        \adjust@\shifted@\advance\adjust@\dimen@ii
        \Nnext@
        \next@\adjust@\copy\ZER@
        \ifN@\hskip\dimen@\else\hskip-\dimen@\hskip-\wd\ZER@\fi
       \fi
      \fi
  \fi\fi\fi\fi
  \ifnum\arrcount@=\m@ne
  \else
   \loop
    \ifdim\goal@>\charht@
    \ifE@\else\hskip-\charwd@\fi
    \Nnext@
    \next@\shifted@\copy\shaft@
    \ifE@\else\hskip-\charwd@\fi
    \advance\shifted@\charht@\advance\goal@-\charht@
    \repeat
   \ifdim\goal@>\z@
    \dimen@\charht@\advance\dimen@-\goal@
    \divide\dimen@\tan@i\multiply\dimen@\tan@ii
    \ifE@\hskip-\dimen@\else\hskip-\charwd@\hskip\dimen@\fi
    \adjust@\shifted@\advance\adjust@-\charht@\advance\adjust@\goal@
    \Nnext@
    \next@\adjust@\copy\shaft@
    \ifE@\else\hskip-\charwd@\fi
   \else
    \adjust@\shifted@\advance\adjust@-\charht@
   \fi
  \fi
  \ifnum\arrcount@=\m@ne
  \else
   \ifnum\arrcount@=\thr@@
   \else
    \ifnum\tcount@=\m@ne
    \else
     \setbox\ZER@
      \hbox{\ifnum\angcount@=117 \arrow@v\else\arrowfont@\fi\char\Tcount@}%
     \ifnum\tcount@=\thr@@
      \advance\adjust@\charht@
      \ifE@\else\ifN@\hskip-\charwd@\else\hskip-\wd\ZER@\fi\fi
     \else
      \ifnum\tcount@=12
       \advance\adjust@\charht@
       \ifE@\else\ifN@\hskip-\charwd@\else\hskip-\wd\ZER@\fi\fi
      \else
       \ifE@\hskip-\wd\ZER@\fi
     \fi\fi
     \Nnext@
     \next@\adjust@\copy\ZER@
     \ifnum\tcount@=13
      \hskip-\wd\ZER@
      \getcos@{\thr@@\p@}%
      \ifE@\hskip-\dimen@\else\hskip\dimen@\fi
      \advance\adjust@-\dimen@ii
      \Nnext@
      \next@\adjust@\box\ZER@
     \fi
  \fi\fi\fi}}%
 \iflabel@i
  \rlap{\hskip\tocenter@
  \dimen@\firstx@\advance\dimen@\secondx@\divide\dimen@\tw@
  \advance\dimen@\ldimen@i
  \dimen@ii\firsty@\advance\dimen@ii\secondy@\divide\dimen@ii\tw@
  \global\multiply\ldimen@i\tan@i\global\divide\ldimen@i\tan@ii
  \ifNESW@\advance\dimen@ii\ldimen@i\else\advance\dimen@ii-\ldimen@i\fi
  \setbox\ZER@\hbox{\ifNESW@\else\ifnum\arrcount@=\thr@@\hskip4\p@\else
   \hskip\tw@\p@\fi\fi
   $\m@th\tsize@@\label@i$\ifNESW@\ifnum\arrcount@=\thr@@\hskip4\p@\else
   \hskip\tw@\p@\fi\fi}%
  \ifnum\arrcount@=\m@ne
   \ifNESW@\advance\dimen@.5\wd\ZER@\advance\dimen@\p@\else
    \advance\dimen@-.5\wd\ZER@\advance\dimen@-\p@\fi
   \advance\dimen@ii-.5\ht\ZER@
  \else
   \advance\dimen@ii\dp\ZER@
   \ifnum\slcount@<6 \advance\dimen@ii\tw@\p@\fi
  \fi
  \hskip\dimen@
  \ifNESW@\let\next@\llap\else\let\next@\rlap\fi
  \next@{\smash{\raise\dimen@ii\box\ZER@}}}%
 \fi
 \iflabel@ii
  \ifnum\arrcount@=\m@ne
  \else
   \rlap{\hskip\tocenter@
   \dimen@\firstx@\advance\dimen@\secondx@\divide\dimen@\tw@
   \ifNESW@\advance\dimen@\ldimen@ii\else\advance\dimen@-\ldimen@ii\fi
   \dimen@ii\firsty@\advance\dimen@ii\secondy@\divide\dimen@ii\tw@
   \global\multiply\ldimen@ii\tan@i\global\divide\ldimen@ii\tan@ii
   \advance\dimen@ii\ldimen@ii
   \setbox\ZER@\hbox{\ifNESW@\ifnum\arrcount@=\thr@@\hskip4\p@\else
    \hskip\tw@\p@\fi\fi
    $\m@th\tsize@\label@ii$\ifNESW@\else\ifnum\arrcount@=\thr@@\hskip4\p@
    \else\hskip\tw@\p@\fi\fi}%
   \advance\dimen@ii-\ht\ZER@
   \ifnum\slcount@<9 \advance\dimen@ii-\thr@@\p@\fi
   \ifNESW@\let\next@\rlap\else\let\next@\llap\fi
   \hskip\dimen@\next@{\smash{\raise\dimen@ii\box\ZER@}}}%
  \fi
 \fi
}
\def\outCD@#1{\def#1{\Err@{\noexpand#1must not be used within \string\CD}}}
\newskip\preCDskip@
\newskip\postCDskip@
\preCDskip@\z@
\postCDskip@\z@
\def\preCDspace#1{\RIfMIfI@
 \onlydmatherr@\preCDspace\else\advance\preCDskip@#1\relax\fi\else
 \onlydmatherr@\preCDspace\fi}
\def\postCDspace#1{\RIfMIfI@
 \onlydmatherr@\postCDspace\else\advance\postCDskip@#1\relax\fi\else
 \onlydmatherr@\postCDspace\fi}
\def\predisplayspace#1{\RIfMIfI@
 \onlydmatherr@\predisplayspace\else
 \advance\abovedisplayskip#1\relax
 \advance\abovedisplayshortskip#1\relax\fi
 \else\onlydmatherr@\preCDspace\fi}
\def\postdisplayspace#1{\RIfMIfI@
 \onlydmatherr@\postdisplayspace\else
 \advance\belowdisplayskip#1\relax
 \advance\belowdisplayshortskip#1\relax\fi
 \else\onlydmatherr@\postdisplayspace\fi}
\def\PreCDSpace#1{\global\preCDskip@#1\relax}
\def\PostCDSpace#1{\global\postCDskip@#1\relax}
\def\CD#1\endCD{%
 \outCD@\cgaps\outCD@\rgaps\outCD@\Cgaps\outCD@\Rgaps
 \preCD@#1\endCD
 \advance\abovedisplayskip\preCDskip@
 \advance\abovedisplayshortskip\preCDskip@
 \advance\belowdisplayskip\postCDskip@
 \advance\belowdisplayshortskip\postCDskip@
 \vcenter{\offinterlineskip
  \vskip\preCDskip@\Let@\global\colcount@\@ne\global\rowcount@\z@
  \everycr{%
   \noalign{%
    \ifnum\rowcount@=\Rowcount@
    \else
     \getrgap@\rowcount@\vskip\getdim@
     \global\advance\rowcount@\@ne\global\colcount@\@ne
    \fi}}%
  \tabskip\z@
  \halign{&\global\xoff@\z@\global\yoff@\z@
   \getcgap@\colcount@\hskip\getdim@
   \hfil\vrule\height\ten@\p@\width\z@\depth\z@
   $\m@th\displaystyle{##}$\hfil
   \global\advance\colcount@\@ne\cr
   #1\crcr}\vskip\postCDskip@}%
 \preCDskip@\z@\postCDskip@\z@
 \def\getcgap@##1{\ifcase##1\or\getdim@\z@\else\getdim@\standardcgap\fi}%
 \def\getrgap@##1{\ifcase##1\getdim@\z@\else\getdim@\standardrgap\fi}%
 \let\Width@\relax\let\Height@\relax\let\Depth@\relax\let\Rowheight@\relax
 \let\Rowdepth@\relax\let\Colwidth@\relax
}
\let\enddocument\bye
\def\alloc@#1#2#3#4#5{\global\advance\count1#1by\@ne
  \ch@ck#1#4#2%
  \allocationnumber=\count1#1%
  \global#3#5=\allocationnumber
  \wlog{\string#5=\string#2\the\allocationnumber}}
\catcode`\@=\active

\catcode`\"=12
\font\blackital=cmmib10  \skewchar\blackital='177
\font\sblackital=cmmib7  \skewchar\sblackital='177
\font\ssblackital=cmmib5  \skewchar\ssblackital='177
\font\sanss=cmss10
\font\ssanss=cmss8 scaled 900
\font\sssanss=cmss8 scaled 600
\font\blackboard=msbm10
\font\sblackboard=msbm7
\font\ssblackboard=msbm5
\font\caligr=eusm10
\font\scaligr=eusm7
\font\sscaligr=eusm5

\def\tx#1{{\fam0\relax#1}}

\newfam\bifam
\textfont\bifam=\blackital
\scriptfont\bifam=\sblackital
\scriptscriptfont\bifam=\ssblackital
\def\bi#1{{\fam\bifam\relax#1}}

\newfam\bbfam
\textfont\bbfam=\blackboard
\scriptfont\bbfam=\sblackboard
\scriptscriptfont\bbfam=\ssblackboard
\def\bb#1{{\fam\bbfam\relax#1}}

\newfam\ssfam
\textfont\ssfam=\sanss
\scriptfont\ssfam=\ssanss
\scriptscriptfont\ssfam=\sssanss
\def\ss#1{{\fam\ssfam\relax#1}}

\newfam\clfam
\textfont\clfam=\caligr
\scriptfont\clfam=\scaligr
\scriptscriptfont\clfam=\sscaligr
\def\cl#1{{\fam\clfam\relax#1}}

\mathchardef\za="710B  
\mathchardef\zb="710C  
\mathchardef\zg="710D  
\mathchardef\zh="7111  
\mathchardef\zl="7115  
\mathchardef\zm="7116  
\mathchardef\zn="7117  
\mathchardef\zx="7118  
\mathchardef\zr="711A  
\mathchardef\zs="711B  
\mathchardef\zt="711C  
\mathchardef\zvf="711E 
\mathchardef\zc="7120  
\mathchardef\zf="7127  
\mathchardef\zL="7003  
\mathchardef\zW="700A  

\catcode`\"=\active

\def\*{{\textstyle *}}
\def\ast{{\scriptstyle *}}

\def\ra{\rightarrow }
\def\tT{\widetilde{\sT}}

\def\proof{\demo{Proof}}
\def\endproof{\hfill \vrule height4pt width6pt depth2pt \enddemo}

\def\N{{\bb N}}
\def\R{{\bb R}}

\def\Z{{\bb Z}}

\def\bV{{\bi V}}
\def\bT{{\bi T}}
\def\bZ{{\bi Z}}

\def\sP{{\ss P}}
\def\sT{{\ss T}}
\def\sV{{\ss V}}
\def\sv{{\ss v}}

\def\xd{\tx{d}}
\def\xid{\tx{id}}
\def\bT{\overline{\sT}}
\def\A{\cA}
\def\Li{\cL}

\def\cA{{\cl A}}
\def\cL{{\cl L}}
\def\cX{{\cl X}}

\input paper.st\relax
\hsize=37pc
\hoffset=-10pt
\vsize=53pc
\voffset=6pt
\TagsOnRight
\document
\RefWarnings

        \title
                LIE BRACKETS ON AFFINE BUNDLES
        \endtitle

        \author

      Katarzyna  Grabowska \\
        Division of Mathematical Methods in Physics \\
                University of Warsaw \\
                Ho\D{z}a 74, 00-682 Warszawa, Poland \\
                 {\tt konieczn\@fuw.edu.pl} \\
                        \\
            Janusz Grabowski \\
        Institute of Mathematics \\
                Polish Academy of Sciences \\
                ul. \'Sniadeckich 8, P.O.Box 21, 00-956 Warszawa 10, Poland\\
                {\tt jagrab\@mimuw.edu.pl} \\
                        \\
              Pawe\l\ Urba\'nski \\
                Division of Mathematical Methods in Physics \\
                University of Warsaw \\
                Ho\D{z}a 74, 00-682 Warszawa, Poland \\
                {\tt urbanski\@fuw.edu.pl}
        \endauthor

                \thanks{Supported by KBN, grant No 2 PO3A 041 18}

\abstract
Natural affine analogs of Lie brackets on affine bundles are studied.
In particular, a close relation to Lie algebroids and a duality with
certain affine analog of Poisson structure is established as well as
affine versions of complete lifts and Cartan exterior calculi.
\endabstract

        \maketitle
        \hl1{Introduction}
        It is known that the framework based on vector bundles is not
satisfactory for classical mechanics. For example, the phase space for the
charged particle is not the cotangent bundle of the configuration
manifold, but an affine bundle over it ([\Ref{We1}], [\Ref{MT}]). Also in
the frame-independent formulation of Newtonian analytical mechanics the
phase space is an affine bundle ([\Ref{Be}], [\Ref{Tu}], [\Ref{Ko}],
[\Ref{Pi}]). The third example is the frame-independent formulation of
time-dependent mechanics ([\Ref{U1}]).  In this case the space of
infinitesimal configurations (the domain of a Lagrangian) is no longer a
vector bundle (the tangent bundle), but the affine bundle of first jets.

        On the other hand,  the fundamental structure for the Lagrangian, first
order mechanics is the Lie algebroid structure of the tangent bundle.
In [\Ref{We2}], [\Ref{Li}], [\Ref{Mar}] one can find an attempt to apply more general Lie
algebroids in Lagrangian mechanics.
        It is natural question to combine affine and Lie-algebroidal aspects of
mechanics. An attempt to do this has been done in [\Ref{Sa},\Ref{Sa1},\Ref{Sa2}] in
the context of Euler-Lagrange equations.

        In this paper we discuss purely geometrical and algebraic aspects
motivated by natural examples of objects which carry a structure of an affine
bundle with a canonical bracket.

These examples are presented in Section 2, together with basic concepts
of differential geometry on affine bundles.

An abstract generalization of properties of the algebraic structures
studied in Section 2 -- Lie affgebroid -- is developed in Section 3.
It basically coincides with the concepts developed in [\Ref{Sa},\Ref{Sa2}].
Affine versions of Poisson and Jacobi structures are introduced.

In Section 4 we introduce some cohomology of Lie algebroids which helps
to classify Lie affgebroids up to an isomorphism.

An affine-linear duality studied in Section 5 is used then, in Section 6,
to show that Lie affgebroids are particular substructures of Lie algebroids.
It makes possible to define complete lifts and Cartan exterior calculus
(Section 7).

Natural affine bundle structures associated with epimorphisms of vector
bundles are used in Section 8, where an affine analog of the well-known
correspondence: Lie algebroid structure on $E$ $\longleftrightarrow$
linear Poisson structure on $E^{\textstyle *}$, is found.

\hl1{Affine bundles: notation and preliminaries}

\noindent
    Let $\zt_i\colon A_i\rightarrow M$ be an affine bundle modelled on a
vector bundle $\sv(\zt_i)\colon \sV(A_i)\rightarrow M$, $i=1,2,3$.
Note that the space $\cA_i$ of sections of $\zt_i$ is an affine
space modelled on the space $\sV(\cA_i)$ of sections of
$\sv(\zt_i)$. Moreover, $\sV(\cA_i)$ is a $C^\infty(M)$ module.
    For an affine bundle morphism $\zvf\colon A_1\rightarrow A_2$ we
denote by $\zvf_\sv \colon \sV({A_1})\rightarrow \sV({A_2})$ its
linear part, i.e.
        $$\zvf_\sv(v) = \zvf(a+v) - \zvf(a) \ \ \text{for}\ \ a\in A,\ v\in
\sV(A), \ \ \zt_1(a) = \sV(\zt_1)(v).
                                                         \tag \label{Fb1}$$
\noindent
        Throughout this paper we shall consider  only  bundle  morphisms
over the identity on the base. We will denote by $\ss{Aff}_M (A_1,A_2)$ (resp.
$\ss{Hom}_M(V_1,V_2)$) the set of such morphisms in the affine
(resp. vector) case. We shall also write $\ss{Aff}(A,\R)$ instead of
$\ss{Aff}_M(A, M\times \R)$ and $\ss{Lin}(V,\R)$ instead of
$\ss{Hom}_M(V,M\times \R)$.
   For a bi-affine mapping
            $$ F\colon A_1\times _M A_2 \rightarrow A_3
                                                            $$
we denote by $F^\sv$ and ${}^\sv F$, respectively, the mappings
                $$\align
    F^\sv&\colon A_1\times_ M \sV(A_2) \rightarrow \sV(A_3) \\
    &\colon (a_1,v_2)  \mapsto (F(a_1,\cdot ))_\sv(v_2)
                                          \tag \label{Fb2}\endalign$$
and
            $$\align
    {}^\sv F&\colon \sV(A_1)\times_ M A_2 \rightarrow \sV(A_3) \\
    &\colon (v_1,a_2)  \mapsto(F(\cdot,a_2 ))_\sv(v_1).
                                           \tag \label{Fb3}\endalign$$
        These mappings are, respectively, affine-linear and linear-affine in
the obvious sense.
By $F_\sv$ we denote the bilinear part of $F$, i.e.
        $$\align
    F_\sv&\colon \sV(A_1) \times _M \sV(A_2) \rightarrow \sV (A_3) \\
    &\colon  (v_1,v_2) \mapsto(F^\sv(\cdot ,v_2))_\sv(v_1) =
({}^\sv F(v_1,\cdot))_\sv(v_2)\\
    &= F(a_1 +v_1, a_2 +v_2) - F(a_1 + v_2, a_2) + F(a_1,a_2) -
F(a_1, a_2 +v_2).
                                           \tag \label{Fb4}\endalign$$

\noindent
There are obvious generalizations of $F^\sv$ and $F_\sv$ for the case of
$n$-affine mappings.

\bigskip
    By a {\it special vector bundle} we mean a vector bundle with a
distinguished non-vanishing section, and by a {\it special affine
bundle} we mean an affine bundle modelled on a special vector
bundle.
    Let $\bV =(V,\zf)$ and $\bV' = (V',\zf')$ be special vector bundles. A
morphism  $F\colon V\rightarrow V'$ of vector bundles is called a {\it
morphism of special vector bundles} if
$\zf$ and $\zf'$ are $F$-related.
    A {\it morphism of special affine bundles} is a morphism of affine bundles
such that its linear part is a morphism of special vector bundles.
        We call an $n$-dimensional affine bundle $\zt\colon A\rightarrow
M$ {\it trivial} if $A= M\times \R^n$ as an affine bundle, and {\it
vectorially trivial} if its model vector bundle is trivial: $\sV(A) =
M\times \R^n$.

\bigskip\noindent
{\bf Remark.} Any section $\zs\in \ss{Sec}(\zt)$ of $A$ induces an obvious
isomorphism $I_\zs\in \ss{Aff}_M(A,\sV(A))$ of affine bundles:
                $$ A_m\ni a\mapsto a-\zs(m) \in \sV(A_m)
                                         \tag \label{Fb5}$$
If $A$ is vectorially trivial, this gives a trivialization
                $$I_\zs\colon A\rightarrow M\times \R^n .
                                         \tag \label{Fb6}$$
Thus, vectorially trivial affine bundles are trivializable but not trivial
if no canonical section is given.

\bigskip
        An affine morphism
        $$D\colon \ss{Sec}(A_1)\rightarrow \ss{Sec}(A_2)$$
we call an {\it affine differential operator of order $n$} if
                        $$ D_\sv\colon \ss{Sec}(\sV(A_1))\rightarrow
\ss{Sec}(\sV(A_2))
                                         $$
is a differential operator of order $n$. Similarly,  $D$ is an {\it affine
quasi-derivation} (or, following Mackenzie [\Ref{Ma}], a {\it covariant
differential operator}) if $D_\sv$ is a quasi-derivation, i.e. if there is
a derivation in $C^\infty(M)$ (represented by a vector
field on $M$) $\widehat{D}_\sv$, called the {\it anchor} of $D$, such that
                $$ D_\sv (fX)=fD_\sv(X) + \widehat{D}_\sv(f)X.
                                         $$
        If $\sV(A_1) = \sV(A_2) = M\times \R$, i.e. $A_1, A_2$ are
one-dimensional and vectorially
trivial, then $D$ will be called an {\it affine derivation} if $D_\sv
\colon C^\infty(M)\rightarrow C^\infty (M)$ is a derivation in the algebra
$C^\infty(M)$ of smooth functions ($D_\sv$ is represented by a vector
field $\widehat{D}_\sv$ on $M$).

        We denote by $\ss{AffDiff}^n(A_1,A_2)$ the space of affine
differential
operators of order $n$, by $\ss{AffQder}(A_1,A_2)$ the space of affine quasi
derivations, and by $\ss{AffDer}(A_1,A_2)$ the space of affine derivations
(we will write simply $\ss{AffDer}(A,\R)$ instead of $\ss{AffDer}
(A,M\times\R)$, etc.).
        These spaces are affine spaces modelled on the corresponding vector
spaces $\ss{AffDiff}^n(A_1,\sV(A_2))$, $\ss{AffQder}(A_1,\sV(A_2))$, and
$\ss{AffDer}(A_1,\sV(A_2))$, respectively.

                \claim \c{e}{Example}{}{\rm
\label{Cb1}
        For a one-dimensional vectorially trivial affine bundle $Z$, the
space $\ss{AffDer}(Z,\R)$ is not only a vector space, but also a
$C^\infty(M)$-module. Let $D\colon \ss{Sec}(Z)\rightarrow C^\infty(M)$ be
an affine derivation with the linear part
$\widehat{D}$. It is easy to  see that the formula
      $$ [D_1,D_2] = \widehat{D}_1\circ D_2 -\widehat{D}_2 \circ D_1
                                         \tag \label{Fb7}$$
defines a Lie bracket on $\ss{AffDer}(Z,\R)$. Moreover,
       $$ [D_1,fD_2] = f [D_1,D_2] + \widehat{D}_1(f)D_2,
                                         \tag \label{Fb8}$$
      i.e. we have a Lie pseudoalgebra structure on the
$C^\infty(M)$-module  $\ss{AffDer}(Z,\R)$ (cf. [\Ref{Ma}]).

        Let $\zs$ be a section of $Z$. With the trivialization
$I_\zs \colon Z\rightarrow M\times\R$ (see \Ref{Fb6})we can identify
sections of $Z$ with functions on $M$. An affine derivation $D\in
\ss{AffDer}(Z,\R)$ induces an affine derivation on $C^\infty(M)$:
                $$ f\mapsto \widehat{D}(f) + D(\zs),
                                         \tag \label{Fb9}$$
         i.e, $\ss{AffDer}(Z,\R)$ can be identified with $\cX(M)\times
C^\infty(M)$. Here $\cX(M)$ denotes the Lie algebra of vector fields on
$M$. With this identification, the bracket~\Ref{Fb7} takes the form
                $$ [(\widehat{X},g),(\widehat{Y},h)] =
([\widehat{X},\widehat{Y}],
\widehat{X}(h) - \widehat{Y}(g)).
                                         \tag \label{Fb10}$$
           A basic fact in differential geometry is that derivations on
$C^\infty(M)$ can be identified with vector fields, i.e. sections of the
Lie algebroid $\sT M$. Now, we show that we have a similar interpretation
for affine derivations.

        The canonical $\R$-action $\zc\colon (t,z_p)\mapsto (z_p+t)$ on $Z$
induces an $\R$-action $\zc_{\ast}$ on $\sT Z$. The space of orbits of
this action we will denote $\tT Z$ and we will call it the {\it reduced
tangent bundle} [\Ref{TU}]. This  is  a  vector  bundle  over  $M$   and
sections of $\tT Z$
can be identified with $\zc_{\ast}$-invariant vector fields on $Z$.
The projection $\zt \colon Z\rightarrow M$ induces a projection
                $$ \widetilde{\zt} = \tT(\zt) \colon \tT Z\rightarrow \sT M.
                                         \tag \label{Fb11}$$
      With this projection, $\tT Z$ becomes a one-dimensional affine bundle
modelled on $\sT M \times \R$ (for details see Section~\Ref{S1}).

    A $\zc_\ast$-invariant vector field on $Z$ (a section of
$\widetilde{\sT}Z$) is a derivation in the algebra of functions on $Z$
which sends affine functions to functions which are constant on fibers.
In local affine coordinates $(x^a,s)$ on $Z$, providing an  isomorphism  of
$Z$ with $\sV(Z)=M\times\R$, an invariant vector field reads
        $$ X=g^a(x)\frac{\partial}{\partial x^a}
-\za(x)\frac{\partial}{\partial s}= X_M-\za(x)\frac{\partial}{\partial
s}.                                              \tag \label{Fb12}$$
The action $\zc(t,x^a,s)=(x^a,s+t)$  gives  rise  to  the  corresponding
fundamental vector field $X_0=-\frac{\partial}{\partial s}$.
We can identify a section $\zs$ of $Z$ with an
affine function $f_\zs$ on $Z$ with the directional coefficient $-1$  (i.e.
$X_0(f_\zs)=1$) and such that $f_\zs\circ \zs =0$. Consequently, we can
interpret a section $X$ of $\widetilde{\sT}Z$ as a mapping
        $$ D_X\colon \ss{Sec}(Z) \rightarrow C^\infty (M).
                                             \tag \label{Fb13}$$
      It is easy to verify that $D_X$ is an affine derivation, i.e. its
linear part is a derivation in $C^\infty(M)$. The Lie bracket of
invariant vector fields is an invariant vector field. Thus, we
have also the Lie bracket of function-valued affine derivations on
$\ss{Sec}(Z)$. It can be trivially verified that this bracket coincides with
the bracket~\Ref{Fb7}.

        The Lie bracket on sections of $\tT Z$ defines a Lie algebroid
structure with the anchor $\widetilde{\zt}$. It is the principal bundle
algebroid of $Z$ interpreted as a principal fibre bundle with the group
$(\R, +)$.
In local coordinates a section
$x=(x^a)\mapsto(x,\zs(x))$ can be represented by the affine function
$f_\sigma(x,  s)=\zs(x)-s$.  For  an  invariant  vector  field  $X$   as
in~\Ref{Fb12},
        $$ Xf_\sigma=g^a\frac{\partial\zs}{\partial x^a}+\za.
                                               \tag \label{Fb14}$$
    Thus $D_X(\sigma)=g^a\dfrac{\partial\zs}{\partial x^a}+\za$ is an
affine derivation on $Z$ with the linear part represented by the vector
field $X_M=g^a(x) \dfrac{\partial}{\partial x^a} =
\widetilde{\zt}\circ X$ on $M$. In fact, $\tT Z$ is a Lie algebroid with
the Lie bracket on sections induced from $\sT Z$. In local coordinates
        $$\align
    X&=X_M-\za\frac{\partial}{\partial s},  \\
            Y&=Y_M-\zb\frac{\partial}{\partial s},  \\
    \left[X,Y\right]& =\left[X_M,Y_M\right] -(X_M(\zb)-Y_M(\za))
\frac{\partial}{\partial s}   .
                                      \tag \label{Fb15}\endalign$$
        Moreover, $\tT Z$ has the distinguished section $X_0$ which
corresponds to the fundamental vector field
$X_0=-\frac{\partial}{\partial s}$ of the $\R$-action on $Z$. It
follows that $\tT Z$ is a special  vector  bundle.  We  summarize  these
observations in the following.
                \claim \c{t}{Theorem}{}
         \label{Cb2}
        There is a canonical isomorphism between $\ss{AffDer}(Z,\R)$
with the bracket~\Ref{Fb7}
and sections of the reduced tangent bundle $\tT Z$ with the canonical Lie
algebroid bracket.
        The fundamental vector field $X_0$ corresponds to the affine
derivation
                $$   \ss{Sec}(Z)\ni \zs \mapsto 1_M
                                         $$
        with vanishing linear part.
                \endclaim

\noindent
It is clear now  that  the  space  of  sections  of  $\bigwedge^n\tT  Z$
corresponds to the space $\bigwedge^n\ss{AffDer}(Z,\R)$ of  affine
multiderivations  on  $Z$   with   values   in    $C^\infty(M)$,    i.e.
multi-affine
skew-symmetric mappings
 $$D\colon \ss{Sec}(Z)\times\dots\times\ss{Sec}(Z)\rightarrow C^\infty(M)
     \tag\label{ma}$$
such that
$$D(a_1,\dots,a_{n-1},\cdot)\colon\ss{Sec}(Z)\rightarrow C^\infty(M)
      \tag\label{ma1}$$
is an affine derivation for any  $a_1,\dots,a_{n-1}\in\ss{Sec}(Z)$.  The
graded space
$$\bigwedge\ss{AffDer}(Z,\R)=\oplus_{n\in\Z}{\bigwedge}^{\! n}\ss{AffDer}
(Z,\R),     \tag\label{ma2}$$
being isomorphic with $\ss{Sec}(\bigwedge\tT Z)=
\oplus_{n\in\Z}\ss{Sec}(\bigwedge^n\tT Z)$ is therefore  a  Gerstenhaber
algebra (see [\Ref{KS}]) with the  wedge  product  and  the  Lie
algebroid Schouten-Nijenhuis
bracket induced by the Lie  algebroid  bracket  on  $\tT  Z$.  With  the
identification  $\ss{Sec}(\tT  Z)=\cX(M)\times  C^\infty(M)$, the   Lie
algebroid bracket on  $\ss{Sec}(\tT  Z)$  is  identified  with  the  Lie
algebroid bracket of first-order differential operators on $M$, so  that
the Gerstenhaber algebra $\bigwedge\ss{AffDer}(Z,\R)$ is identified with
the  Gerstenhaber  algebra  of  skew-symmetric  multilinear  first-order
differential operators. Thus we get the following.

                \claim \c{t}{Theorem}{}
         \label{Cb2a}
There is a canonical Gerstenhaber algebra structure on the graded space
$\bigwedge\ss{AffDer}(Z,\R)$ of multi-affine  skew-symmetric  derivations
on $Z$ with values in $C^\infty(M)$. Using a section  $\zs$  of  $Z$  to
identify     $Z$ with $M\times\R$  we   can      identify
$\bigwedge^n\ss{AffDer}(Z,\R)$ with $\ss{Sec}(\bigwedge^n\sT M)\times
\ss{Sec}(\bigwedge^{n-1}\sT M)$ by
$$D_{(X_n,X_{n-1})}(\zs+f_1,\dots,\zs+f_n)=\langle   X_n,\xd   f_1\wedge
\dots\wedge\xd f_n\rangle+\sum_i(-1)^i\langle   X_{n-1},\xd   f_1\wedge
\overset i\atop\vee \to \dots\wedge\xd f_n\rangle
                        \tag\label{ma4}$$
for $(X_n,X_{n-1})\in\ss{Sec}(\bigwedge^n\sT M)\times
\ss{Sec}(\bigwedge^{n-1}\sT M)$.
With this identification, the wedge product is
$$(X_n,X_{n-1})\wedge(Y_k,Y_{k-1})=(X_n\wedge  Y_k,  X_{n-1}\wedge  Y_k+
(-1)^nX_n\wedge Y_{k-1})
                       \tag\label{ma5}$$
and the Schouten-Nijenhuis bracket reads
$$[(X_n,X_{n-1}),(Y_k,Y_{k-1})]^{SN}=([X_n,Y_k]^{SN}, [X_{n-1},Y_k]^{SN}+
(-1)^{n-1}[X_n,Y_{k-1}]^{SN}),
                       \tag\label{ma6}$$
where the Schouten-Nijenhuis bracket  on  the  right-hand  side  is  the
classical Schouten-Nijenhuis bracket of multivector fields.
\endclaim

                \claim \c{e}{Example}{}{\rm              \label{Cb3}
        Let $Z$ be as above. Consider now the affine space $\cA =
\ss{AffDer}(Z,Z)$ modelled on $\sV(\cA) = \ss{AffDer}(Z,\R) $. It
is easy to see that the commutator
                $$ [D_1,D_2] = D_1\circ D_2 - D_2\circ D_1
                                         \tag \label{Fb16}$$
        defines a bracket on $\cA$ with values in $\sV(\cA)$. This bracket is
an affine analog of the Lie bracket of vector fields. It is bi-affine and
it has the  following properties:
        \list
        \item it is skew-symmetric: $ [D_1,D_2] = - [D_2,D_1]$,
        \item it satisfies the Jacobi identity:
                $$ [D_1,[D_2,D_3]]^\sv + [D_2,[D_3,D_1]]^\sv +
[D_3,[D_1,D_2.]]^\sv
=0 ,
                                         $$
        \item for any $D\in\cA$ the map $\ss{ad}_D = [D,\cdot ]$  is  an
        affine quasi-derivation, i.e.
                $$ [D,fX]^\sv = f[D,X]^\sv + \widehat{D}(f)X,
                                         $$
        for $f\in C^\infty(M)$, $X\in \sV(\cA)$, and some
        $\widehat{D}\in\cX(M)$,
        where $\cX(M)$ denotes the Lie algebra of vector fields on $M$,
        \endlist
        We write here $\widehat{D}$ instead of ${D}_\sv$ to denote the
vector field representing the linear part of $D$. The affine-linear part
$[\, ,\,]^\sv$ of $[\,,\,]$ is given by
                $$ [D,X]^\sv = \widehat{D}\circ X - X\circ D,
                                         $$
$D\in \ss{AffDer}(Z,Z)$, $X\in \ss{AffDer}(Z,\R)$.  The linear part
$[\,,\,]_\sv$ of $[\,,\,]$ coincides with the Lie algebra bracket on
$\ss{AffDer}(Z,\R)$ from the previous example.
        Given a  global section $\zs$ of $Z$, we get~\Ref{Fb6} an
isomorphism
$I_\zs \colon Z\rightarrow M\times \R$  and, consequently,  we can
identify $\ss{AffDer}(Z,Z)$ with $\cX(M)\times  C^\infty(M)$ and sections
of $Z$ with functions on $M$. With these identifications,
                $$\align
        (X,g)(f)&= X(f) + g\, ,  \\
        [(X,g),(Y,h)]&=([X,Y], X(h) - Y(g) + g-h)\, .  \\
                 \tag \label{Fb17}\endalign$$
        Let us recall that $I_\zs$ provides also the identification of
$\ss{AffDer}(Z,\R)$ with $\cX(M)\times C^\infty(M)$,
but the induced bracket is  $[(X,g),(Y,h)] = ([X,Y],
X(h)-Y(g)) $.
        In the previous example we identified an element of
$\ss{AffDer}(Z,\R)$ with a section of the reduced tangent bundle
$\tT Z$. Now, we give a similar construction for $\ss{AffDer}(Z,Z)$.

      There is a canonical diffeomorphism $\sT Z = \tT Z\times _MZ$ and a
$\R^2$-group action on $\sT Z$:
                $$((t,r),(v,z))\mapsto (v + r X_0, z+t),
                                         \tag \label{Fb18}$$
         where $X_0$ is the distinguished section $X_0$ of $\tT Z$ which
corresponds to the fundamental vector field of the $\R$-action on
$Z$ (Example~\Ref{Cb1}). In local coordinates this action reads
                $$ (x^a,s,{\dot x}^b, \dot s)\mapsto  (x^a, t+s, {\dot x}^b,
\dot s -r).
                                         \tag \label{Fb19}$$
        We reduce $\sT Z$ with respect to the group homomorphism
                $$ \R^2\rightarrow \R \colon (r,t) \mapsto r+t
                                         \tag \label{Fb20}$$
i.e., we consider a manifold $\bT Z$ of orbits of the kernel group $(r,-r)$
with the induced group action
                $$ (t, [(v,z)]) \mapsto [(v,z+t)] = [(v+tX_0, z)].
                                         \tag \label{Fb21}$$
        $\bT Z$ is an affine bundle over $M$ and its model bundle is
canonically isomorphic to $\tT Z$. The pair $(\bT Z, X_0)$ is a special
affine bundle over $M$ denoted by $\bT \bZ$.
        Sections of $\bT Z$ correspond to vector fields on $\sT Z$,
invariant
with respect to the $\R$-action
                $$ (r,(v,z))\mapsto (v+rX_0, z-r).
                                         \tag \label{Fb22}$$
         Such a vector field has the local form
                $$ X = g^a(x)\frac{\partial}{\partial
x^a}+(s-g(x))\frac{\partial}{\partial s}=
X_M+(s-g(x))\frac{\partial}{\partial s}.
                                         \tag \label{Fb23}$$
        The  vector  field  $X$  preserves,   as   the   derivation   in
$C^\infty(Z)$, the
affine subspace of functions with the directional coefficient $-1$, i.e.
such that $X_0(f) =1$. Since such affine functions represent sections of
$Z$, we get a mapping
                $$ D_X \colon \ss{Sec}(Z) \rightarrow \ss{Sec}(Z).
                                         \tag \label{Fb24}$$
Let $\zs$ be a section of $Z$ and $f_\zs$ the corresponding function on $Z$.
In local coordinates
                $$ f_\zs(x,s) =\zs(x) -s
                                         \tag \label{Fb25}$$
and
                $$X(f_\zs) = X_M(\zs) + g(x) -s.
                                         \tag \label{Fb26}$$
        Consequently,
                $$ D_X(\zs) = X_M(\zs) + g.
                                         \tag \label{Fb27}$$
        The linear part of $D_X$ is a derivation in the algebra $C^\infty(M)$
given by the vector field $X_M$.
We conclude that sections of $\bT Z$ can be viewed as the affine derivations
$D \colon \ss{Sec}(Z) \ra\ss{Sec}(Z)$, i.e. we have the canonical isomorphism
of the space of the affine derivations $\ss{AffDer}(Z,Z)$ and the
sections of $\bT Z$.
Moreover,
      $$[X_M+(s-g(x))\frac{\partial}{\partial s},
Y_M+(s-h(x))\frac{\partial}{\partial s}]=[X_M,Y_M]-(X_M(h)-Y_M(g)+g-h)
\frac{\partial}{\partial s},
                            \tag\label{zz} $$
so  that  $[D_X,D_Y]$  corresponds   exactly   to   the   element   from
$\ss{AffDer}(Z,\R)$ associated with the invariant vector field
$[X,Y]$. Thus we get a theorem analogous to Theorem 1.
                \claim \c{t}{Theorem}{}
         \label{Cb4}
        There is a canonical isomorphism between $\ss{AffDer}(Z,Z)$
with the bracket~\Ref{Fb16} and sections of the affine bundle $\bT Z$
with the canonical bracket.
        \endclaim
        \hfill $\spadesuit$     }\endclaim


\hl1{Lie brackets on affine bundles}

\noindent
        The Example~\Ref{Cb3} justifies the following definition.

      \claim \c{d}{Definition}{}{\rm                \label{Cb5}
An {\it affine Lie bracket} on  an  affine  space  $\A$  is a bi-affine map
   $$\left[\cdot,\cdot\right]\colon  \A\times\A\ra \sV(\A)$$
which
   \list
   \item is skew-symmetric:
$\left[\za_1,\za_2\right]=-\left[\za_2,\za_1\right]$ and
   \item satisfies the Jacobi identity:
        $$\left[\za_1,\left[\za_2, \za_3\right]\right]^\sv+\left[\za_2,
\left[\za_3,
\za_1\right] \right]^\sv+ \left[\za_3,\left[\za_1,
\za_2\right]\right]^\sv=0.$$
   \endlist
   An affine space equipped with an affine Lie bracket we shall call a
{\it Lie affgebra}.
Note that the term {\it affine Lie algebra} has been  already  used  for
certain types of Kac-Moody algebras.
      }\endclaim

      \claim \c{d}{Definition}{}{\rm                \label{Cb6}
   If $A$ is an affine bundle over $M$ modelled on $\sV(A)$ then a {\it Lie
affgebroid structure} on $A$ is an affine Lie bracket on sections of $A$
and a morphism $\zg\colon A\ra\sT M$ of affine bundles (over $\xid $
on $M$) such that $[\za,\cdot]^\sv$ is a quasi-derivation with the anchor
$\zg(\za)$, i.e.
      $$ \left[\za,fX\right]^\sv=f\left[\za,X\right]^\sv+\zg(\za)(f)X
                                             \tag \label{Fb28}$$
   for all $\za\in\A$, $X\in\sV(\A)$, $f\in C^\infty(M)$.
      }\endclaim

\bigskip\noindent{\bf Remark.}
The above definition goes back to the one used in [\Ref{Sa}] but
without the additional assumption that the base manifold  $M$  is
fibred  over $\R$  and   that   $\zg(\za)$   are   vector   fields
projectable   onto $\frac{\partial}{\partial t}$. This requirement
was motivated by the particular aim of [\Ref{Sa}] to get
time-dependent Euler-Lagrange equations, but it makes impossible
to get all Lie affgebras. We will  see  later
(Theorem~\Ref{Cba0}) that the particular  case  of [\Ref{Sa}]
means just vanishing  of  certain cohomology.

\bigskip\noindent
The following  fact  has  been  proved,  in  principle,  in  [\Ref{Sa}],
Proposition 1
.
      \claim \c{t}{Theorem}{}                     \label{Cb7}
   For every Lie affgebroid structure
$\left[\cdot,\cdot\right]\colon \A\times\A\ra\sV(\A)$ its vector part
     $$\left[\cdot,\cdot\right]_\sv\colon \sV(\A)\times\sV(\A)\ra\sV(\A)$$
     is a  Lie  algebroid  structure  with  the  anchor  map  $\zg_\sv$.
Moreover,
$\left[\za, \cdot\right]^\sv$ is a derivation of $\left[\cdot,
\cdot\right]_\sv$ for every $\za\in\A$. Conversely, if we have an
affine-linear map $\left[\cdot,\cdot\right]^0 \colon
\A\times\sV(\A)\ra\sV(\A)$ which satisfies \Ref{Fb28} for certain affine
morphism $\zg\colon A\ra \sT M$ and
such that $\left[\cdot,\cdot\right]_0={\left[\cdot,\cdot\right]}^0_\sv$ is
skew symmetric and $\left[\za,\cdot\right]^0$ is a derivation of
$\left[\cdot,\cdot\right]_0$, for any $\za\in\A$, then there is a unique Lie
affgebroid bracket $\left[\cdot,\cdot\right]$ on $\A$ such that
$\left[\cdot,\cdot\right]^0= \left[\cdot,\cdot\right]^\sv$.
      \endclaim

       \claim \c{d}{Definition}{}{\rm           \label{Cb8}
        Let $\zt\colon Z\rightarrow M$ be a one-dimensional affine bundle
with $\sV(Z) = M\times \R$. An affine Lie bracket on $\ss{Sec}(Z)$
        $$      \{\,,\,\} \colon \ss{Sec}(Z)\times \ss{Sec}(Z) \rightarrow
C^\infty(M)$$
        is called an {\it aff-Poisson} (resp.{\it aff-Jacobi}) {\it
bracket} if
        $$ \{\za,\cdot \} \colon \ss{Sec}(Z) \rightarrow C^\infty(M)
                                         \tag \label{Fb29}$$
        is an affine derivation (resp. an affine first order differential
operator) for every $\za\in \ss{Sec}(Z)$.       }\endclaim

\noindent
We use  the  term {\it aff-Poisson}, since  {\it  affine  Poisson
structure}  has already  a different meaning in the literature.

                \claim \c{t}{Theorem}{}          \label{Cb9}
        For every aff-Poisson (resp. aff-Jacobi) bracket
        $$\{\,,\,\} \colon \ss{Sec}(Z)\times \ss{Sec}(Z) \rightarrow
C^\infty(M)$$
        its vector part
        $$\{\,,\,\}_\sv \colon C^\infty(M)\times C^\infty(M) \rightarrow
C^\infty(M)$$
        is a Poisson (resp. Jacobi) bracket.  Moreover,
         $$ \{\za,\cdot \}^\sv \colon C^\infty(M) \rightarrow C^\infty(M)
                                         $$
        is a derivation (resp. first-order differential operator) for every
section $\za\in \ss{Sec}(Z)$, which is simultaneously a derivation of the
bracket  $\{\,,\,\}_\sv$.  Conversely,  if  we  have  a  Poisson  (resp.
Jacobi)
bracket $\{\,,\,\}_0$  on  $C^\infty(M)$  and  a  derivation  (resp.  a
first-order differential operator)
          $$ D \colon C^\infty(M) \rightarrow C^\infty(M)
                                         $$
which is simultaneously a derivation of the bracket $\{\,,\,\}_0$, then
there is a unique aff-Poisson (resp.  aff-Jacobi)  bracket   $\{\,,\,\}$
on  $\ss{Sec}(Z)$ such that $\{\,,\,\}_0=\{\,,\,\}_\sv$ and
$D=\{\za,\cdot\}^\sv$  for  a chosen section $\za\in \ss{Sec}(Z)$.
          \endclaim
          \proof
        The above theorem is a direct consequence of Theorem~\Ref{Cb7}.
                \endproof

\noindent
{\bf Remark.} Using a section $\zs$ to identify $\ss{Sec}(Z)$ with
$C^\infty(M)$, we get the aff-Poisson (resp. aff-Jacobi) bracket on
$\ss{Sec}(Z)$ in the form
                $$ \{f,g\} = D(g-f) + \{f,g\}_\sv,
                                         \tag \label{Fb30}$$
        where $D$ is a vector field (resp first-order differential operator)
which  is  a  derivation  of  the   Poisson   (resp. Jacobi)   bracket
$\{\,,\,\}_\sv$. We have clearly $\{f,\cdot \}^\sv = D + \{f,\cdot
\}_\sv$. In particular, for the Poisson case, the aff-Poisson bracket is
a bi-affine skew-symmetric derivation, so it is represented by the pair
$(\zL,D)$, where $\zL$ is the Poisson tensor  of $\{\ ,\ \}_\sv$, in
accordance  with Theorem~\Ref{Cb2a}.  The  Jacobi  identity means  that
$[(\zL,D),(\zL,D)]^{SN}=0$,  i.e.  $[\zL,\zL]^{SN}=0$   and
$[D,\zL]^{SN}=0$.


\hl1{Lie affgebroid brackets and $Qder$-cohomology}

\noindent
        Every Lie affgebroid structure on $A$ determines a Lie algebroid
structure on $\sV(A)$. Moreover, by Theorem~\Ref{Cb7}, every section $\za$
of $A$ determines a derivation $[\za,\cdot ]^\sv$ of the Lie algebroid
bracket on $\ss{Sec}(\sV(A))$ which at the  same  time  is  a
quasi-derivation  on sections of $\sV(A)$.
        Passing to another section we change this derivation by an inner
derivation. As in the case of a Lie algebra, this determines certain
cohomology class of the Lie algebroid $\sV(A)$.
        To be more precise, consider a Lie algebroid structure on a vector
bundle $E$ over $M$.  Let $C^n_{Qder}(E,E)$ be the space of $n$-cochains
with coefficients in the adjoint representation which are quasi-derivations
with respect to each argument. This means that elements of
$C^n_{Qder}(E,E)$ are skew-n-linear maps
                $$ \zm\colon \ss{Sec}(E)\times \cdots\times \ss{Sec}(E)
\rightarrow \ss{Sec}(E),                                         $$
        such that for $X_1,\dots,X_n \in \ss{Sec}(E)$, $f\in C^\infty(M)$,
                $$ \zm(fX_1, X_2,\dots, X_n) =f\zm(X_1,\dots,X_n) +
\widehat{\zm}_{(X_2,\dots,X_n)}(f)X_1
                                       \tag \label{Fb31}$$
        for certain vector field $\widehat{\zm}_{(X_2,\dots,X_n)}$ on $M$.

                \claim \c{t}{Theorem}{}         \label{Cb10}
        The standard Chevalley coboundary operator
            $$ \partial\zm(X_0,X_1,\dots,X_n) = \sum_i (-1)^i [X_i,
\zm(X_0,\overset i\atop\vee \to \dots, X_n )] + \sum_{i<j} (-1)^{i+j}
\zm([X_i,X_j],X_1,\overset i\atop\vee \to \dots\overset j\atop\vee \to
\dots, X_n)
         \tag \label{Fb32}$$
        can be restricted to the subcomplex of multi-quasi-derivations, i.e.
          $$\partial C^n_{Qder}(E,E)\subset C^{n+1}_{Qder}(E,E).
                                         \tag \label{Fb33}$$
        \endclaim
        \proof
        The proof consists of standard calculations in order
        to show that, for $\zm\in C^n_{Qder}(E,E)$,
          $$\partial\zm(fX_0,X_1,\dots,X_n) =
            f\partial\zm(X_0,X_1,\dots,X_n) +\widehat{D}(f)X_0
                                         \tag \label{Fb34}$$
        for certain vector field $\widehat{D}$ on $M$ depending on
$X_1,\dots,X_n$ and $\zm$.
        \endproof

\medskip\noindent
        The corresponding cohomology spaces will be denoted by
$H^n_{Qder}(E,E)$. Let us notice that a Qder-1-cocycle is a
quasi-derivation $D\colon \ss{Sec}(E) \rightarrow \ss{Sec}(E)$ such that
     $$\partial D(X_0,X_1) = [X_0,D(X_1)] + [D(X_0),X_1]-D([X_0,X_1])=0,
                                         \tag \label{Fb49}$$
        i.e. $D$ is a derivation of the Lie bracket. Any Qder-1-coboundary is
of the form
                $$ (\partial X_0)(X) =[X,X_0],
                                         \tag \label{Fb50}$$
         so that 1-coboundaries are just inner derivations. Hence,
                $$ H^1_{Qder}(E,E) = \ss{Der}_{Qder}(E)/\ss{InnDer}(E).
                                         \tag \label{Fb51}$$

                \claim \c{t}{Theorem}{}
         \label{Cb14}
For an affine bundle $A$ over $M$, every Lie affgebroid bracket  on  $A$
has the form
                $$ [\za_0 +X, \za_0 +Y] = D(Y) -D(X) + [X,Y]_\sv,
                                         \tag \label{laf1}$$
where          $\za_0\in\ss{Sec}(A)$,          $X,Y\in\ss{Sec}(\sV(A))$,
$D=[{\za_0},\cdot]^\sv$.
     The isomorphism class of a Lie affgebroid structure  on  an  affine
bundle $\zt\colon  A\rightarrow  M$  with  the  prescribed  vector  part
$[\,,\,]_\sv$ on $\sV(A)$ is determined by the class $[D]\in
H^1_{Qder}(\sV(A),     \sV(A))$     for      the      Lie      algebroid
$(\sV(A),[\,,\,]_\sv)$.
                \endclaim
                \proof
        We have already seen that a Lie affgebroid structure on $A$ determines
a Lie algebroid structure on $\sV(A)$ and a class from $H^1_{Qder}(\sV(A),
\sV(A))$   represented   by   $D=[\za_0,\cdot]^\sv$   for   a    section
$\za_0\in\ss{Sec}(A)$. Let now  be given a  Lie
algebroid bracket $[\,,\,]_\sv$ on
$\ss{Sec}(\sV(A))$ and $[D]\in H^1_{Qder}(\sV(A), \sV(A))$. We choose a
section $\za_0\in \ss{Sec}(A)$ and define a bracket $[\,,\,]\colon
\ss{Sec}(A) \times \ss{Sec}(A)\rightarrow\ss{Sec}(\sV(A))$ by
                $$ [\za_0 +X, \za_0 +Y] = D(Y) -D(X) + [X,Y]_\sv.
                                         \tag \label{Fb52}$$
         It is easy to see that this formula defines a Lie affgebroid bracket
on $A$ with the linear part $[\,,\,]_\sv$. The Jacobi  identity  follows
from the fact that $D$ is a derivation of $[\,,\,]_\sv$, and the existence
of the anchor follows from the fact that $D$ is a quasi-derivation.
If we choose another representant $D'$ of $[D]$, then $D' = D +
\ss{ad}_{X_0}$ for certain $X_0 \in \ss{Sec}(\sV(A))$, and we get another
bracket
     $$ [\za_0 +X, \za_0+Y]'= (D+\ss{ad}_{X_0})(Y-X)+[X,Y]_\sv.
                                         \tag \label{Fb53}$$
         We have clearly an isomorphism between $[\,,\,]$ and $[\,,\,]'$
induced by the isomorphism $A\ni\za \mapsto\za+ X_0(\zt(\za))$.
                \endproof

\bigskip\noindent
The next theorem shows  that  every  Lie  affgebroid  structure  can  be
extended to a standard Lie algebroid.

                \claim \c{t}{Theorem}{}
         \label{Cb15}
        Let $[\,,\,]_0$ be a Lie algebroid bracket on sections of a vector
bundle $E$ over $M$ and let $[D] \in  H^1_{Qder}(E,E)$, i.e.  $D$  is  a
$\ss{Qder}$-1-cocycle. Then, the bracket
defined on sections of $E\oplus\R$ by
                $$ [(X,f),(Y,g)] = \left([X,Y]_0 +fD(Y) - gD(X),
(\widehat{X}+f\widehat{D})(g)-(\widehat{Y}+g\widehat{D})(f)\right)
                                         \tag\label{laf}$$
          is  a  Lie  algebroid  bracket  on  $E\oplus   \R$   extending
$[\,,\,]_0$ on sections of $E=E\oplus\{ 0\}\subset E\oplus \R$ and
such that $[\ss{Sec}(A),\ss{Sec}(A)]\subset\ss{Sec}(E)$,  where
$A$  is  the affine  subbundle  $A=E\oplus\{1\}$  of   $E\oplus
\R$,   modelled   on $\sV(A)=E=E\oplus\{ 0\}$.

Conversely, every Lie algebroid extension of $[\,,\,]_0$ with  the above
properties has the form \Ref{laf} and
the isomorphism class of such a Lie algebroid bracket on $E\oplus \R$
is determined by the class $[D]$. The Lie algebroid bracket \Ref{laf}
reduced  to  the  sections  of  $A$  gives  a  Lie  affgebroid   bracket
\Ref{laf1} on $A$ with $\za_0=(0,1)\in\ss{Sec}(E\oplus \R)$.
                \endclaim
                \proof
    The  $C^\infty(M)$-module  $\ss{Sec}(E\oplus  \R)$  is generated  by
the sections of $E$ and the section $\za_0=(0,1)\in\ss{Sec}(E\oplus \R)$.
There is a unique extension of the bracket $[\,,\,]_0$  on  sections  of
$E$ to a Lie algebroid  bracket  on  $\ss{Sec}(E\oplus  \R)$  for  which
$\ss{ad}_{\za_0}=D$, namely the bracket \Ref{laf}. The Jacobi  identity
follows from the fact that $D$ is a derivation of $[\,,\,]_0$, and  the
existence of the anchor follows from the fact that $D$ is a quasi-derivation.
If we choose another  representant  $D'\in  [D]$, $D'=D+\ss{ad}_Y$, the
corresponding Lie algebroid bracket $[\,,\,]'$ on $\ss{Sec}(E\oplus \R)$
goes  to \Ref{laf} with respect to the isomorphisms $(X_p,t)\rightarrow
(X_p+tY(p),t)$ of $\ss{Sec}(E\oplus \R)$. Indeed, with respect  to  this
isomorphism $(X,f)$ goes to $(X+fY,f)$, i.e. sections  $X$  of  $E$  are
preserved      and  $\za_0$      goes      to   $\za_0+Y$.      Since
$\ss{ad}_{\za_0+Y}=D+\ss{ad}_Y=D'$,  the   bracket   $[\,,\,]$   goes   to
$[\,,\,]'$.

Conversely, if $[\,,\,]$ is a Lie algebroid bracket on $\ss{Sec}(E\oplus
\R)$ extending $[\,,\,]_0$ on the sections of $E$  and such  that
$[\ss{Sec}(A),\ss{Sec}(A)]\subset\ss{Sec}(E)$,   then   the bracket has
clearly the form~\Ref{laf} for $D=\ss{ad}_{\za_0}$. Since
     $$[(X,1),(Y,1)]=[\za_0+X,\za_0+Y]=([X,Y]_0+D(Y-X),0),
                       $$
the bracket \Ref{laf}, reduced to $A$, gives \Ref{laf1}.
                \endproof

\medskip\noindent
        {\bf Remark.} The mentioned isomorphisms are, in general,
non-canonical and depend on the choice of a section or a representant of
the class in $H^1_{Qder}(\sV(A),\sV(A))$. A canonical way to interpret a
Lie affgebroid as an affine subbundle of a Lie algebroid will be
described in the next section.

        \hl1{Affine - special linear duality}

\noindent
        For  an  affine  bundle  $\zt  \colon   A\rightarrow   M$,   let
$\zt^\dag\colon  A^\dag \rightarrow M$ be the vector bundle of
affine functions on fibers of $A$, i.e. the fiber $A^\dag _p$
at $p\in M$ is the space $\ss{Aff}(A_p,\R)$ of affine functions on $A_p$.
It is easy to see that $A^\dag $ carries a canonical smooth
vector bundle structure.
    Locally, for coordinates $(x^a)$ on $M$ and for a local basis $e_0,
e_1, \dots, e_n$ of sections of $A$ (i.e. $e_1 -e_0,e_2 -e_0, \dots , e_n
-e_0 $, form a local  basis  of  sections  of  $\sV(A)$),  we  have  the
corresponding coordinates $(x^a,s^i)$ on $A$ such that for $\za\in A_p$
            $$\za = e_0(p) +\sum_i s^i(\za)(e_i(p) - e_0(p)).
                                                            $$
    The adapted local coordinates on $A^\dag $ are $(x^a,a_i,
b)$ and for $f\in A^\dag _p$, we have
        $$ f(\za) =  a_i(f)s^i(\za) +b(f) .
                                                            $$
    Note that $A^\dag $ is a special vector bundle with the
distinguished section $1_A$ corresponding to the constant function with
the value $1$. The special vector bundle $(A^\dag ,1_A)$ will
be called the {\it vector dual} to the affine bundle $A$.  Conversely,
with a special vector bundle $\bV = (V,\zf_0)$ we associate an affine
bundle $\bV^\ddag$ which
is identified with the affine subbundle of the
dual vector bundle $V^{\textstyle *} $, defined by the equation
        $$ (V,\zf_0)^\ddag = \{ f\in V^{\textstyle *} \colon  \langle f,
\zf_0\rangle  = 1\} .
                                             \tag \label{Fb43}$$
For the detailed discussion on affine structures involved  see
Section~\Ref{S1}
Since for any vector  bundle  $V$  we  have the obvious  identification
  $$\ss{Sec}(V^{\textstyle *})\ni X\mapsto \imath_X\in\ss{Lin}(V,\R)
                 $$
of  sections  of  $V^{\textstyle *}$  with  linear  functions  on
the   dual   bundle $V$, sections of  $(V,\zf_0)^\ddag$
can  be  viewed  as  those  linear functions on $V$ which take the value 1
on $\zf_0(M)\subset V$.

\noindent
       The evaluation mappings
        $$\align
    \langle \cdot ,\cdot \rangle_A &\colon A\times _MA^\dag
\rightarrow  M\times \R \\
    &\colon (a,\zf) \mapsto \zf(a)
                                       \tag \label{Fb44}\endalign$$
        and
        $$\align
    \langle \cdot ,\cdot \rangle_{\zf_0} &\colon  (V,\zf_0)^\ddag\times _M
V \rightarrow  M\times \R \\
    &\colon (f,v) \mapsto f(v)
                                       \tag \label{Fb444}\endalign$$
are affine-linear and such that $\langle f ,\zf_0 \rangle_{\zf_0} =1$,
$\langle a, 1_A\rangle_A =1$, i.e. for each $a\in A_p$, $f\in
(V,\zf_0)^\ddag_p$, the mappings
        $$ \langle f ,\cdot \rangle_{\zf_0} \colon V_p
\rightarrow \R, \qquad  \langle a, \cdot \rangle _A \rightarrow \R
                                                \tag \label{Fb45}$$
are morphisms of special vector spaces.

By a {\it pairing} between an affine bundle $A$ and a special vector
bundle $\bV=(V,\zf_0)$ over $M$ we mean an affine-linear mapping
        $$ \langle \cdot ,\cdot \rangle \colon A\times _M V\longrightarrow
M\times \R
                                               \tag \label{Fb46}$$
with the following properties:
    \list
        \item"(i)" $\langle a,\zf_0\rangle =1$ for each $a\in A$,
        \item"(ii)" it is nondegenerate, i.e. the induced morphisms
$V\rightarrow A^\dag$ and $A\rightarrow V^{\textstyle *} $
are injective.
    \endlist
    Note, that the image of the induced mapping $A\rightarrow
V^{\textstyle *} $ is contained in $\bV^\ddag$.

        \claim \c{t}{Theorem}{}                         \label{Cb13}
    Let
        $$ \langle \cdot ,\cdot \rangle \colon A\times _M V\rightarrow
M\times \R
                                \tag \label{Fb47}$$
be a pairing between an affine bundle $A$ and a special vector bundle
$\bV= (V,\zf_0)$. Then the induced morphisms $V\rightarrow A^\dag $ and
$A\rightarrow \bV^\ddag $ are isomorphisms.
            \endclaim
        \claim \c{c}{Corollary}{}                       \label{Cbx14}
We have natural isomorphisms
\list
    \item"(a)" $A\simeq (A^\dag, 1_A )^\ddag$,
    \item"(b)" $(V,\zf_0)\simeq (((V,\zf_0)^\ddag)^\dag,
1_{(V,\zf_0)^\ddag})$. \endlist
        \endclaim
        \demo{Proof of the Theorem}
    Let $\zvf_r\colon V\rightarrow A^\dag$ and $\zvf_l\colon A\rightarrow
\bV^\ddag $ be the mappings  induced by the pairing $\langle \cdot ,\cdot
\rangle$. Since they are injective, it is enough to prove that they are
surjective. Since $1\in\zvf_r(V_p)$ and the intersection of the kernels
$\bigcap_{v\in  V_p}  \ss{Ker}{\langle\cdot,v\rangle}_\sv$
is trivial, $\zvf_r$ is an epimorphism.
Moreover, $\zvf_r(\zf_0)=1_A$.

It remains to show that $A\simeq(V,\zf_0)^\ddag$. The map $A_p\ni
a\mapsto\zvf_l(a)=\langle a, \cdot\rangle\in V_p^{\textstyle *} $ induces
an affine injective map $\zvf_l\colon A\ra \bV^\ddag $. Moreover
$\zvf_l(A)=(V,\zf_0)^\ddag$, since if $\text{codim}\,\zvf_l(A_p)>1$ then
there are infinitely many $\zf\in V_p$ for which $\langle\cdot,\zf
\rangle=1$, in spite of the fact that $\zvf_r$ is injective.
        \endproof

\bigskip\noindent
        In this way we have defined duality between affine bundles on one
hand and special vector bundles on the other hand. The vector bundle
$(A^\dag)^{\textstyle *}$ we  denote by $\widehat A$ and call the {\it
vector hull} of  the  affine  bundle  $A$.  The  affine  bundle  $A$  is
canonically isomorphic to the affine subbundle of $\widehat A$ described
as a 1-level set of the linear function $\imath_{\zf_0}$ on $\widehat{A}$
defined by a non-vanishing section $\zf_0$ of $A^\dag$.

        \claim \c{e}{Example}{}{\rm                         \label{Cbx15}
    Let $A$ be an affine bundle modelled on a one-dimensional vector bundle
$\sV(A)$ over $M$. A fiber of $A^\dag$ is a vector space of
dimension 2 and
for each $a\in A_p$ the fiber $A_p^\dag$ can be split into a direct sum of
two 1-dimensional subspaces: one contains constant functions and the other
one consists of functions vanishing at $a$. The first subspace can be
parameterized by reals and the second by the linear parts, i.e. by
elements of $(\sV(A))^{\textstyle *}_p $. We have then a surjective
mapping
        $$\align
    \zh&\colon A\times _M\sV(A)^{\textstyle *} \times \R\rightarrow
A^\dag \\
    &\colon  (a,\zf,t)\mapsto (b\mapsto \zf(b-a) +t)
                                          \tag \label{Fb48}\endalign$$
which is a linear isomorphism for each $a\in A$.
    The equality $\zh(a,\zf,t) = \zh(a',\zf',t')$, i.e. $$ \zf(b-a) +t =
\zf'(b-a') + t' = \zf'(b-a) + \zf'(a-a') +t', $$ is equivalent to $\zf
=\zf'$ and $t= t' + \zf'(a-a')$. It follows that the level sets of $\zh$
are orbits of the fiber-wise action of $\sV(A)$ on $A\times
_M\sV(A)^{\textstyle *} \times \R$:
        $$ (u,(a,\zf,t)) \mapsto  (a+u,\zf, t + \zf(u)).
                                                 \tag \label{Fbx49}$$

    Since  $\zh_a =\zh(a,\cdot ,\cdot )$ is an isomorphism for each $a$,
we have also the dual isomorphism for $a\in A_p$
        $$(\zh_a)^{\textstyle *} \colon (A_p^\dag)^{\textstyle *}
\rightarrow (\sV(A)_p^{\textstyle *} \times \R)^{\textstyle *}  =
\sV(A)_p\times\R .
                                            \tag \label{Fbx50}$$
   The level sets of the dual action
        $$ \zh^{\textstyle *} \colon A\times _M\sV(A)\times \R \rightarrow
(A^\dag)^{\textstyle *}  = \widehat{A},
                                              \tag \label{Fbx51}$$
$(\zh^{\textstyle *})_a=((\zh_a)^{\textstyle *})^{-1}$,
are the orbits of the  dual
action of $\sV(A)$ on $A\times _M\sV(A)\times \R$, given by the formula
        $$ (u,(a,v,t)) \mapsto  (a-u,v +su, t),
                                               \tag \label{Fbx52}$$
and $\widehat A$ can be viewed as the vector bundle of the orbits
$[(a,v,t)]$.
        In the representations given by $\zh$ and $\zh^{\textstyle *} $,
the canonical pairing between $A^\dag$ and $\widehat{A}$ takes the form
        $$ \langle [(a,\zf,t)],[(a',v,s)]\rangle =\zf(v) +
s\zf(a'-a) +st
                                                \tag \label{Fbx53}$$
        and the embedding $A\hookrightarrow \widehat{A}$ is given by
$a\mapsto [(a,0,1)]$. \hfill $\spadesuit$ }\endclaim

        \claim \c{e}{Example}{}{\rm                   \label{Cb16}
    Let $Z$ be a one-dimensional affine bundle over $M$ modelled on
the trivial bundle $M\times\R$. In this case
        $$ Z\times _M\sV(Z)^{\textstyle *} \times \R = Z\times \R^2 =
Z\times _M\sV_(Z)\times \R.
                                               \tag \label{Fb54}$$
    With these identifications the level sets of $\zh$ and
$\zh^{\textstyle *} $ are, respectively, the orbits of $\R$-actions:
        $$ (\za,(z,s,t)) \mapsto (z+\za,s,t +s\za)
                                                            $$
and
        $$ (\za,(z,s,t))\mapsto (z-\za, s+\za t, t).
                                                            $$
    We can identify $z\in Z_p$ with an affine function $Z_p\ni y\mapsto
z-y$ represented by the triple $(z,-1,0)$. Thus we have inclusions $Z
\hookrightarrow \widehat{Z}$ and $Z\hookrightarrow Z^\dag$ which define an
isomorphism  of $Z^\dag$ and $\widehat{Z}$, and the pairing
        $$\align
    \langle \cdot ,\cdot \rangle_Z&\colon Z^\dag\times _M
\widehat{Z} = Z^\dag\times _MZ^\dag  \rightarrow
\R  \\
    &\colon ([(z,s,t)],[(z',s',t')])  \mapsto st'-s's(z-z')-s't
=s(t'+s'(y-z')) -s'(t +s(y-z)).
                                    \tag \label{Fb55}\endalign$$
    Let $f,g$ be affine functions represented by $(z,s,t)$ and
$(z',s',t')$ respectively, i.e. $f(y) = t+s(y-z)$ and $g(y)= t'
+s'(y-z')$. Then the formula~\Ref{Fb55} reads
        $$ \langle f,g\rangle = X_0(g)f-X_0(f)g,
                                               \tag \label{Fb56}$$
where $X_0=-\frac{\partial}{\partial y}$ is the fundamental vector field
associated with the $\R$-action on $Z$.
        For $s\neq 0$ the equivalence classes can be parameterized by
$Z\times \R_{\textstyle *} $:
        $$ (z,s)\mapsto [(z,s,0)]
                                              \tag \label{Fb57}$$
     and for $s=0$ by $\R$:
        $$ t\rightarrow [(z,0,t)].
                                              \tag \label{Fb58}$$
It follows that $Z^\dag$ can be parameterized by the
set $Z^0=(Z\times\R_{\ast})\cup(M\times\R\times\{0\})$ with the
obvious projection on $M$. In this representation the vector
bundle structure of $Z^\dag$ is given by the formulae
        $$\align
    (s,a) + (s',a')&=\cases (s + \frac{a'}{a+a'}(s'-s),a+a' ) \ \
&\text{if} \ a+a' \neq 0 \\  (a(s'-s),0) &\text{if} \ \ a+a'=0 \endcases\\
    (t,0) + (t',0)      &= (t+t',0) \\
    (t,0)+(s,a)     &= (s-\frac{t}{a}, a )\\
        \zl(s,a) &= (s,\zl a)\\
        \zl(t,0) &= (\zl t,0).
                                       \tag \label{Fb59}\endalign$$
    The vector $(s,a)$ represents the affine function on $Z_p$ vanishing
at $s$ and with directional coefficient $a$, $a\neq 0$ and $(t,0)$
represents the constant affine function equal to $t$.
    \hfill $\spadesuit$ }\endclaim

       \claim \c{e}{Example}{}{\rm       \label{Cb17}
        Let $Z $ be a one-dimensional affine bundle which is vectorially
trivial. The vector bundle $\tT Z$ introduced in Example~\Ref{Cb1} has a
non-vanishing distinguished section $X_0$, so that it is a special vector
bundle.
        An element of the dual affine bundle $\tT^\ddag Z= (\tT Z)^\ddag$
is a hyperspace in $\tT_p Z$, transversal to $X_0(p)$. It follows that a
section of $\tT^\ddag Z$ can be identified with a $\zc$-invariant
horizontal distribution, i.e. a  connection on $Z$.  Standard
representation of such connection is a $\zc$-invariant 1-form $\za$ on $Z$
such that $\langle \za, X_0\rangle = 1$.
        In [\Ref{TU}] the affine bundle $\tT^\ddag Z$ has been denoted by
$\sP  Z$. This is  an  affine  version  of  the   cotangent   bundle
$T^{\textstyle *}M$.
      \hfill $\spadesuit$       }\endclaim

                \hl1{The Lie algebroid hull of a Lie affgebroid}

\noindent
        We begin this section with the following lemma.

      \claim \c{l}{Lemma}{}                         \label{Cb21}
   Let $A$ be an affine subspace of a finite  dimensional  vector  space
$V$, $A=\{v\in V \colon \zf(v)=1\}$ for certain $\zf\in V^{\textstyle
*} $, $\zf\ne 0$. Then,
   \list
   \item $A^{\otimes n}=\{a_1\otimes\ldots\otimes a_n\colon  a_i\in A\}$
spans $V\otimes\ldots\otimes V$ -- the n-fold tensor product of $V$;

   \item $A^{\wedge n}=\{a_1\wedge\ldots\wedge a_n\colon  a_i\in A\}$
spans $\bigwedge^nV$;

   \item $A\wedge\sV(A)^{\wedge n}=\{ a\wedge  v_1\wedge\ldots\wedge  v_n:
a\in A, v_i\in\sV(A)\}$ spans $\bigwedge^{n+1}V$.
   \endlist
      \endclaim
      \proof
        {(1)} Inductively, $A=a_0+\sV(A)$  for  any  $a_0\in A$
and $a_0\notin\sV(A)$, so
        $$\ss{span}A=\ss{span }\{a_0,\sV(A)\}=V.$$
Assume, that $A^{\otimes n}$ spans the
$n$-fold tensor product of $V$ and let $v_1,\ldots,v_{n+1}\in V$. Since
$v_1\otimes\ldots\otimes v_n$ is a linear combination $\sum  \zl_j  u_j$
for  a  basis  $u_j$  of  $V^{\otimes  n}$  from  $A^{\otimes  n}$   and
$v_{n+1}=\lambda a_0+v$, $v\in\sV(A)$,
we can write $v_{n+1}=a_0+v+(\lambda-1)a_0$, so that
      $$\multline
v_1\otimes\ldots \otimes v_{n+1}= v_1\otimes\ldots\otimes v_n\otimes
\left((a_0+v)+(\lambda-1)a_0\right)=
\sum \zl_i u_i\otimes\left( (a_0+v)+(\zl-1)a_0\right) \\=
\sum \zl_i u_i\otimes(a_0+v) + \sum (\zl-1)\zl_i u_i\otimes a_0,
            \endmultline                     \tag \label{Fb72}$$
i.e. $v_1\otimes\ldots\otimes v_{n+1}$ is a linear combination of elements
of $A^{\otimes (n+1)}$.

\noindent(2) follows easily from (1).

\noindent{(3)} It follows from (2) that $A^{\wedge (n+1)}$ spans
$\bigwedge^{(n+1)}V$. Since for $a_0\in A$, $v_i\in\sV(A)$,
      $$
(a_0+v_1)\wedge\ldots\wedge(a_0+v_{n+1})= \qquad\sum_{i\geq
2}^{n+1}(-1)^{i+1}a_0\wedge v_1\wedge\overset i\atop\vee\to\dots\wedge
v_{n+1}\,+\,(a_0+v_1)\wedge v_2\wedge\ldots\wedge v_{n+1},
                 $$
$A\wedge \sV(A)^{\wedge n}$ spans $A^{\wedge(n+1)}$, so $\bigwedge^{n+1}V$.
      \endproof

                \claim \c{t}{Theorem}{}
         \label{Cbx22}
        For every Lie affgebroid bracket $[\,,\,]$ on an affine bundle $A$
over
$M$ there is a unique Lie algebroid bracket $[\,,\,]^\wedge$  on
$\widehat{A}$ such that $[\,,\,]$ is the restriction of $[\,,\,]^\wedge$ to
sections of $A$.
                \endclaim
                \proof
For a given section $\za_0$ of $A$ there is a unique isomorphism of  $A$
with the affine subbundle $A'=\sV(A)\oplus\{ 1\}$ of  $\sV(A)\oplus  \R$
such that $\za_0+X$ corresponds  to  $(X,1)$  for  any  section  $X$  of
$\sV(A)$.
This isomorphism extends uniquely to an isomorphism  of  vector  bundles
$\widehat{A}$ and $\sV(A)\oplus \R$ which maps  sections  $f\za_0+X$  to
$(X,f)$.
        According to Theorem~\Ref{Cb15},
there is a Lie algebroid bracket $[\,,\,]$ on $\sV(A)\oplus
\R$   which,   restricted   to   $\ss{Sec}(A)\simeq\ss{Sec}(A')$, gives
a Lie affgebroid bracket, so the existence follows. This Lie affgebroid
bracket defines the bracket and anchors for a
generating set of sections of $\widehat{A}$ ($A$ spans $\widehat{A}$), so
the Lie algebroid bracket on $\widehat{A}$ extending the Lie  affgebroid
bracket on $A$ is unique.
        \endproof

\bigskip\noindent
Note that the above result has been already mentioned in [\Ref{Sa}].

        Now, let $\zt\colon A\rightarrow M$ be an affine bundle and let $\zf
=1_A$ be the section of $A^\dag$ such that $A$ is the level-1-set of $\zf$
in $\widehat{A} = (A^\dag)^{\textstyle *} $. Let $[\,,\,]$ be a Lie
algebroid bracket on sections of $\widehat{A}$, let $\zL$ be the
corresponding linear Poisson tensor on $A^\dag$ and let $\xd$ be the
corresponding exterior derivative. Let $\ss{Sec}(\widehat{A})\ni  X\mapsto
X^c\in\cX(\widehat{A})$ be the complete lift of sections of the Lie
algebroid $\widehat{A}$ to vector fields on $\widehat{A}$ (cf.
[\Ref{MX},\Ref{GU2},\Ref{GU3}]).

                        \claim \c{t}{Theorem}{}
                 \label{Cb23}
        The following are equivalent
        \list
        \item the restriction of $[\,,\,]$ to the sections of $A$ is a Lie
affgebroid bracket;
        \item $[X,Y] \in \ss{Sec}(\sV(A))$ for each $X,Y \in
\ss{Sec}({A})$;
        \item $\xd \zf$ =0;
        \item  $ \Li_{V_\sT(\zf)}\Lambda=0$, where $\Li$ is the Lie derivative
and $V_\sT(\zf)$ is the vertical lift of $\zf$;
        \item the complete Lie algebroid lift $X^c$ of any section $X$ of $A$
is a vector field on $\widehat{A}$ which is tangent to $A\subset \widehat{A}$.
        \endlist
                \endclaim
                \proof
(1)$\Rightarrow $(2) \quad is trivial.

(2)$\Rightarrow $(3) \quad Let $X,Y\in\ss{Sec}(A)$. Then
                $$ \xd\zf(X,Y) = \widehat{X} (\zf(Y)) - \widehat{Y}(\zf(X)) -
\zf([X,Y]) =0,$$
        since $\zf(X)= \zf(Y) =1$ and $\zf$ vanishes on $\sV(A)$. It
follows from Lemma~\Ref{Cb21} that $\xd \zf= 0$.

(3)$\Leftrightarrow $(4)\quad follows from the identity $[\zL,
V_\sT(\zm)]^{SN}= V_\sT(\xd\zm)$ ([\Ref{GU1}], Thm 15.d), where $\zm$ is a
multisection of $A^\dag$ and $[\,,\,]^{SN}$ is the Schouten-Nijenhuis
bracket.

(4)$\Rightarrow $(5)\quad For every $\zm \in \ss{Sec}(A^\dag)$,
$X^c(\imath_\zm)  =  \imath_{\Li_X\zm}$,  where   $\Li$   is   the   Lie
derivative and $\imath_\zm$ is  the  linear  function  on  $\widehat{A}$
corresponding to $\zm$ (see [\Ref{MX},\Ref{GU2},\Ref{GU3}]).
Since, for $X\in \ss{Sec}(A)$,
       $$ \Li_X\zf = \tx{i}_X\xd\zf + \xd\tx{i}_X\zf = \tx{i}_X
0+\xd(1)=0,
                                         $$
we get $X^c(\imath_\zf) = 0$. Since $A$ is a level-set of $\imath_\zf$, it
follows that $X^c$ is tangent to $A$.

(5)$\Rightarrow $(1)\quad  Let $X,Y\in \ss{Sec}(A)$. Since $X^c$ is tangent
to  $A$, the function $X^c(\imath_\zf)$  vanishes   on   $A$.   But   the
function $X^c(\imath_\zf)$ is
linear and $A$ spans $\widehat{A}$, so that $X^c(\imath_\zf) = 0$.
        Since $V_\sT(Y)(\imath_\zf)  =  V_\sT(\zf(Y))  =  1$,  and  (cf.
[\Ref{GU2},\Ref{GU3}])
$V_\sT([X,Y])(\imath_\zf) = [X^c, V_\sT(Y)](\imath_\zf)$, we get
$V_\sT([X,Y])(\imath_\zf) = 0$.
        But $V_\sT([X,Y])(\imath_\zf) = V_\sT(\zf([X,Y]))$, so $\zf([X,Y])=0$
and $[X,Y] \in \ss{Sec}(\sV(A))$. Now, the Jacobi identity and the  existence
of the anchor map for the bracket restricted to a map on
$\ss{Sec}(A)\times\ss{Sec}(A)$
follows directly from the corresponding properties of the Lie  algebroid
bracket $[\,,\,]$.
                \endproof

\bigskip\noindent
According to Theorem~\Ref{Cb23}  (5),  the  complete  lift  $X^c$of  any
section $X$ of $A$, so also for any section of $\sV(A)$, is  tangent  to
the submanifold $A$ of $\widehat{A}$ and we  can  define  the  {\it  Lie
affgebroid complete lift} $X^A$ as the restriction of the  vector  field
$X^c$ to $A$. Since $V_\sT(X)(\imath_\zf)=\langle X,\zf\rangle\circ\zt$,
also the vertical lift of any section $X$ of  $\sV(A)$  is  tangent  to
$A$. Its restrictions to $A$ we will denote $V_\sT^A(X)$ and call  the
{\it Lie affgebroid vertical lift} of $X$.

        \claim \c{t}{Theorem}{}         \label{Cba23}
For any Lie affgebroid bracket  $[\,  ,\,]$  on  $A$  and  all  sections
$X_1,X_2\in\ss{Sec}(A)$, $Y_1,Y_2\in\ss{Sec}(\sV(A))$,
  \list\item"(a)" $[X_1^A,X_2^A]=([X_1,X_2])^A$,
       \item"(b)" $[Y_1^A,Y_2^A]=([Y_1,Y_2])^A$,
       \item"(c)" $[X_1^A,V_\sT(Y_1)]=V_\sT^A([X_1,Y_1])$,
       \item"(d)" $[V_\sT^A(Y_1),V_\sT^A(Y_2)]=0$.\endlist
       Here the brackets on the left-hand sides are clearly the brackets
of vector fields on $A$.
\endclaim
\proof
The above identities follow immediately from  the  analogous  identities
for the complete and vertical lifts of sections of a Lie algebroid  (cf.
[\Ref{GU2},\Ref{GU3}]).
\endproof

\bigskip\noindent
Theorem~\Ref{Cb23} implies that a Lie affgebroid  structure  on  $A$  is
determined, in fact, by  a  Lie algebroid structure on a vector bundle
$E$ (which is $\widehat{A}$) and a non-vanishing 1-cocycle
$\zf\in\ss{Sec}(E^{\textstyle *})$. This is  a  particular  case  of  a
{\it generalized Lie algebroid} in the sense of Iglesias and Marrero
[\Ref{IM}] or a {\it Jacobi algebroid} in the sense  of  [\Ref{GM}].  In
the particular case, when  the 1-cocycle $\zf$ is trivial, we get the
affine  Lie  algebroid  in  the sense of [\Ref{Sa}] as shows the
following.

        \claim \c{t}{Theorem}{}         \label{Cba0}
In the notation  preceding Theorem \Ref{Cb23},  assume  that  we   have
a Lie affgebroid bracket on $A$ with an anchor $\zg$ such  that
$\zf=\xd  t$ for certain $t\in C^\infty(M)$. Then $t\colon
M\rightarrow\R$ is  a  fibration
over the open subset $t(M)$ of $\R$ and $\zg$ maps sections  of  $A$  to
vector    fields    on    $M$     which     are     projectable     onto
the canonical vector field $\frac{\partial}{\partial t}$ on $\R$.
\endclaim
\proof Since
$$\zg(X)(t)=\langle X,\xd t\rangle=\langle X,\zf\rangle=1
                     \tag\label{sa}$$
for any section $X$ of $A$, the function $t$ is regular, so that
$t:M\rightarrow\R$ is a fibration and $\zg(X)$ is projectable onto
$\frac{\partial}{\partial t}$.
\endproof

\medskip\noindent
{\bf Remark.} In the case when $\zf$ is not trivial  one  can  make  it
trivial for the prize of extending  the  base  manifold  by  $\R$.  This
construction for generalized Lie algebroids can be found in [\Ref{IM}]. The
sections  of  the  corresponding  Lie  algebroid  can   be   viewed   as
time-dependent sections of the original one.


                \hl1{Exterior calculus}

        \claim \c{d}{Definition}{}{\rm
\label{Cb26}
        An {\it affine $k$-form} on an affine bundle $A$ over $M$ is a
skew-symmetric $k$-affine map
                $$ \overline{\zm}\colon  \ss{Sec}(A)\times \ss{Sec}(A)\times
\cdots
\times \ss{Sec}(A) \rightarrow C^\infty(M)
                                         $$
which comes from an affine morphism
             $$A\times_M\dots\times_MA\rightarrow M\times\R,
                          $$
i.e.          $\overline{\zm}(a_1,\dots,a_{k-1},\cdot)_\sv$           is
$C^\infty(M)$-linear for any $a_1,\dots,a_{k-1}\in\ss{Sec}(A)$.
                }\endclaim

\medskip\noindent
        {\bf Remark.} This definition coincides with that one in [\Ref{Sa}].
The
space $\ss{Aff\zW}^k(A)$ of affine $k$-forms  is a $C^\infty(M)$-module in
the obvious way. We have also the standard wedge product on the graded space
$\ss{Aff\zW}(A)=\underset k\in \Z \to\oplus \ss{Aff\zW}^k(A)$ which makes
$\ss{Aff\zW}(A)$ into a graded associative commutative algebra. Here,
$\ss{Aff\zW}^k(A) =\{0\}$ for $k<0$ and $\ss{Aff\zW}^0(A) =C^\infty(M)$.

                \claim \c{t}{Theorem}{}
         \label{Cb27}
        There is a canonical isomorphism between the space $\ss{Aff\zW}^k(A)$
of affine $k$-forms on $A$ and the space $\zW^k(\widehat{A})
=\ss{Sec}(\wedge^k A^\dag)$ of sections of the $k$-th exterior power of
$A^\dag$, given by restrictions of $k$-forms on $\widehat{A}$ to $A$. This
isomorphism can be extended to an isomorphism of graded associative
algebras $\ss{Aff\zW}(A)$ and $\zW(\widehat{A})$.
        \endclaim
                \proof
        By    definition,     every     affine     1-form     $\zm\colon
\ss{Sec}(A)\rightarrow C^\infty(M)$ is represented by a unique section of
$A^\dag$ and vice-versa. Sections of $\wedge ^k A^\dag$ can be identified
with $k$-forms on $\widehat{A} = (A^\dag)^{\textstyle *} $ and the
restriction of any $k$-form
                $$ \zm\colon \ss{Sec}(\widehat{A})\times \cdots\times
\ss{Sec}(\widehat{A}) \rightarrow C^\infty(M)
                                         $$
        to $\ss{Sec}(A)\times \cdots\times \ss{Sec}(A)$ gives clearly an
affine
$k$-form $\overline{\zm}$ on $A$. Conversely, let $\overline{\zm}$ be an
affine $k$-form on $A$. Let us choose a section $\za_0$ of $A$.
        Let $v_1^o,\dots,v_n^o $ be a local basis of sections of $\sV(A)$.
Then $\za_0^o =\za_0$ and $\za_i^o = \za_0  + v^o_i$, $i=1,\dots,n$, is a
basis of sections of $\widehat{A}$. There is a unique $k$-form $\zm$ on
$\widehat{A}$ such that
                $$ \zm(\za^o_{i_1},\dots, \za^o_{i_k}) =
\overline{\zm}(\za^o_{i_1},
\dots, \za^o_{i_k})
                $$
for all $i_1,\dots,i_k\in\{ 0,\dots,n\}$.
        We will show that $\overline{\zm}$ is the restriction of $\zm$ to
sections of $A$. First, we observe that, since $\overline{\zm}$ is
skew-symmetric and affine,
                $$ \overline{\zm}(\za_0 +v_1,\dots, \za_0 +v_k) = \sum
(-1)^{i+1}
\overline{\zm}^\sv(\za_0,v_1, \overset i\atop\vee \to \dots, v_k) +
\overline{\zm}_\sv (v_1,\dots,v_k)
\tag\label{A}$$
for certain affine-multilinear and linear maps $\overline{\zm}^\sv$ and
$\overline{\zm}_\sv$, respectively, and
        for all $v_1,\dots,v_k \in \ss{Sec}(\sV(A))$.
A similar formula for $\zm$ reads
                $${\zm}(\za_0 +v_1,\dots, \za_0 +v_k) = \sum
(-1)^{i+1}
{\zm}(\za_0,v_1, \overset i\atop\vee \to \dots, v_k) +
{\zm} (v_1,\dots,v_k).\tag\label{Aa}$$
Hence,
                $$ \overline{\zm}^\sv(\za_0,v_2,\dots,v_k) =
\overline{\zm}(\za_0,\za_0
+v_2,\cdots\za_0 + v_k),
                                         $$
        so that
                $$ \overline{\zm}^\sv(\za_0^o,v^o_{i_2},\dots,v^o_{i_k}) =
\overline{\zm}(\za_0^o,\za^o_{i_2},\dots,\za^o_{i_k}) =
\zm(\za_0^o,\za^o_{i_2},\dots,\za^o_{i_k}) =
\zm(\za_0,v^o_{i_2},\dots,v^o_{i_k}) .
                                         \tag\label{B}$$
        This implies further, by \Ref{A}, that
                $$  \overline{\zm}_\sv(v^o_{i_1},\dots,v^o_{i_k}) =
\zm(v^o_{i_1},\dots,v^o_{i_k})
                                         \tag\label{C}$$
for all $i_1,\dots,i_k\in\{ 0,\dots,n\}$.
    Since both $\zm$ and $\overline{\zm}^\sv$ are $C^\infty(M)$-linear in
the last (k-1)-arguments, \Ref{B} implies that
                        $$  \overline{\zm}^\sv(\za_0^o,v_{2},\dots,v_{k}) =
\zm(\za_0,v_{2},\dots,v_{k})
                                         $$
        for any $v_2, \dots, v_k \in \ss{Sec}(\sV(A)) $.
        Similarly, \Ref{C} implies
                        $$  \overline{\zm}_\sv(v_1,v_{2},\dots,v_{k}) =
\zm(v_1,v_{2},\dots,v_{k})
                                         $$
        and finally, due to \Ref{A}, that
                $$\overline{\zm}(\za_0 +v_1,\dots, \za_0 +v_k) = {\zm}(\za_0
+v_1,\dots, \za_0 +v_k),
                                         $$
        i.e. that $\overline{\zm}$ is the restriction of $\zm$ to sections of
$A$.
                \endproof

\bigskip\noindent
        It is well known that a Lie algebroid structure determines an exterior
derivative on forms. So, the Lie algebroid structure on $\widehat{A}$,
extending the Lie affgebroid structure on $A$, gives rise to an exterior
derivative
                $$ \xd \colon \zW^n(\widehat{A}) \rightarrow
\zW^{n+1}(\widehat{A})
                                         $$
        which, in view of the previous theorem, can be seen as an exterior
derivative on $\ss{Aff\zW}^n (A)$:
                $$ \xd\colon \ss{Aff\zW}^n (A) \rightarrow \colon
\ss{Aff\zW}^{n+1}
(A).
        $$
The derivation property of $\xd$ means that
                $$ \xd(\zm\wedge \zn) = \xd \zm \wedge \zn + (-1)^{deg
\zm}\zm\wedge \xd\zn.
                                                $$
        Since we obtain affine $k$-forms on $A$ from $k$-forms on
$\widehat{A}$ by restriction, the Cartan-like formula for $\xd$ is clearly
           $$\align
 \xd \overline{\zm}(\za_0,\dots,\za_{n})& = \sum (-1)^{i+1}
\zg({\za}_i)({\zm} (\za_0,\overset i\atop\vee \to
\dots,\za_{n}))
+ \sum_{i<j} (-1)^{i+j}{\zm}([\za_i,\za_j],\za_0,
\overset i\atop\vee \to \dots \overset j\atop\vee \to \dots,
\za_{n})=\\
&\sum (-1)^{i+1}\zg({\za}_i)( \overline{\zm} (\za_0,\overset
i\atop\vee \to\dots,\za_{n}))
- \sum_{i<j} (-1)^{n+i+j}\overline{\zm}^\sv(\za_0,\overline{\za}_1,
\overset i\atop\vee \to \dots \overset j\atop\vee \to \dots,
\overline{\za}_{n},[\za_i,\za_j]),
\tag \label{Fb79}\endalign$$
        where $\za_i\in\ss{Sec}(A)$ and $\overline{\za}_i = \za_i-\za_0$
for $i>0$.
        The corresponding affine de~Rham complex on $A$  is  equivalent,
via the isomorphism described  in  Theorem~\Ref{Cb27},  to  the  de~Rham
complex of the Lie algebroid $\widehat{A}$.
Using the well-known fact that graded derivations of degree 1 on
$\zW(\widehat{A})$ defined by the Cartan formula
      $$ \xd\zm(X_0,X_1,\dots,X_n) = \sum_i (-1)^i\widehat{X}_i
(\zm(X_0,\overset i\atop\vee \to \dots, X_n )) + \sum_{i<j} (-1)^{i+j}
\zm([X_i,X_j],X_1,\overset i\atop\vee \to \dots\overset j\atop\vee \to
\dots, X_n)
         \tag \label{Car}$$
satisfy
$\xd^2=0$  if  and  only  if  the  bracket  $[\,,\,]$  on  sections   of
$\widehat{A}$ is  a  Lie algebroid bracket  with  the  anchor  $X\mapsto
\widehat{X}$, we  get immediately  the following [\Ref{Sa}]:

                \claim \c{t}{Theorem}{}             \label{Cb28}
        For a graded derivation on  $\ss{Aff\zW}(A)$ of degree 1
        $$\xd\colon \ss{Aff\zW}(A) \rightarrow  \ss{Aff\zW}(A),$$
        defined by \Ref{Fb79}, $\xd^2 =0$ if and only if $[\,,\,]$ is a Lie
        affgebroid bracket with the anchor $\zg$.
         \endclaim

\noindent
{\bf Remark.} The 1-cocycle $\zf$ determining $A$ inside $\widehat A$
(cf. Theorem~\Ref{Cb23}) represents the affine 1-form with the constant
value 1 on $\ss{Sec}(A)$. The ideal $\zW(\widehat{A})_\zf$ generated
by $\zf$ in the exterior algebra $\zW(\widehat{A})$ is a differential ideal
and, since $A^\dag/\langle\zf\rangle=\sV(A)^*$, the differential algebra
              $$\zW(\sV(A))=\zW(\widehat{A})/\zW(\widehat{A})_\zf
                           $$
represents just the exterior calculus for the Lie algebroid structure on
$\sV(A)$.

       \hl1{Affine bundles of vector bundle epimorphisms \label{S1}}
\noindent
        According to Corollary~\Ref{Cbx14}, sections of $A$ may be identified
with linear
functions on $A^\dag$ which take the value $1$ on the section $\zf =1_A$
and, with respect to this identification, the Lie affgebroid bracket on $A$
may be identified with a linear Poisson bracket on $A^\dag$ restricted to
this affine subspace of functions. We shall use also another identification
in which the Lie affgebroid bracket comes from an aff-Poisson bracket.
To do this we need some preparations.

        Let $\zt_i\colon E_i \rightarrow M$, i=1,2, be vector bundles  and let
$\zr\in \ss{Hom}_M(E_1,E_2)$ be a vector bundle epimorphism, i.e.
$\zr(E_1) =E_2$.  The kernel of $\zr$  is a vector
subbundle of $E_1$ which will be denoted by $K_\zr$.

                \claim \c{p}{Proposition}{}                     \label{Cbx21}
        The fibration $\zr\colon E_1 \rightarrow E_2$ is an affine fibration
modelled on the vector bundle $E_2\times_MK_\zr$ over $E_2$.
                \endclaim
                \proof
        Let $e_1,e_2\in (E_1)_m$ with $\zr(e_1)= \zr(e_2) \in (E_2)_m$. Then
$\zr(e_1 -e_2) =0$, i.e. $e_1-e_2 \in (K_\zr)_m$. It is easy to
see that this defines a smooth affine bundle structure on the
total space $E_1$ and base $E_2$ which is modelled on $E_2 \times
_M K_\zr$.
                \endproof

\noindent
        The canonical embedding $K_\zr \hookrightarrow E_1$ induces the
projection
                $$ \zr^0\colon  E_1^{\textstyle *}  \rightarrow
K_\zr^{\textstyle *}.
                                         \tag \label{Fb62}$$
        The kernel of this projection can be identified with $E_2^{\textstyle
*} $ embedded by
                $$ \zr^{\textstyle *} \colon  E_2^{\textstyle *} \rightarrow
E_1^{\textstyle *}.
                                         \tag \label{Fb63}$$
         Let us denote:
        \list
        \item"(a)"  by $\ss{LinSec}(\zr)$ - the affine subspace of those
sections $\zs \in \ss{Sec}(\zr)$ which are morphisms of vector bundles,
i.e.
                $$ \ss{LinSec}(\zr) = \{\zs\in \ss{Hom}_M(E_2,E_1) \colon
\zr\circ
\zs = id_{E_2} \} .
                                         \tag \label{Fb64}$$
        Since for $\zs_1, \zs_2 \in \ss{LinSec}(\zr)$, we have $\zr\circ
(\zs_1-\zs_2)=0$, the affine space $\ss{LinSec}(\zr)$ is modelled
on $\ss{Hom}_M(E_2,K_\zr)$;
        \item"(b)" by $\ss{Proj}(E_1,K_\zr)$ - the affine space of those
$\zs\in \ss{Hom}_M(E_1, K_\zr)$ which are projections.
        For $\zs_1,\zs_2 \in \ss{Proj}(E_1,K_\zr)$,  we have $(\zs_1 -
\zs_2)|_{K_\zr} =0$, i.e. $(\zs_1 - \zs_2) \in
\ss{Hom}_M (E_1/K_\zr,K_\zr)\simeq\ss{Hom}_M (E_2,K_\zr)$ --
the model vector space for $\ss{Proj}(E_1,K_\zr)$.
        \item"(c)" by $\ss{ProjSec}(E_1^{\textstyle *} \otimes _MK_\zr)$ the
set of those sections $\zs\in \ss{Sec}(E_1^{\textstyle *} \otimes_M K_\zr)$
which represent projections onto $K_\zr$, i.e. $\langle X,\zs\rangle = X$
for any $X\in \ss{Sec}(K_\zr)$, where $\langle ,\rangle $ is the obvious
contraction. The model space is here clearly $\ss{Sec}(E_2^{\textstyle
*} \otimes_M K_\zr)$.
                \endlist

        \claim \c{t}{Theorem}{}         \label{Cb22}
        The following affine spaces are canonically isomorphic:
        \list
        \item"(a)" $\ss{LinSec}(\zr)$,
        \item"(b)" $\ss{Proj}(E_1,K_\zr)$,
        \item"(c)" $\ss{ProjSec}(E_1^{\textstyle *} \otimes_M K_\zr)$,
        \item"(d)" $\ss{LinSec}(\zr^0)$,
        \item"(e)" $\ss{Proj}(E_1^{\textstyle *} ,E_2^{\textstyle *} )$,
        \item"(f)" $\ss{ProjSec}(E_1^{\textstyle *} \otimes_M
E_2^{\textstyle *} )$.
        \endlist
                \endclaim
                \proof
        Any linear section $\zs\colon E_2\rightarrow E_1$ induces the vector
bundle decomposition
                $$ E_1 = K_\zr \oplus \zs(E_2),
                                         \tag\label{Fb65}$$
        so the projection of $E_1$ onto $K_\zr$, and vice versa. Any such
projection can be represented by the corresponding section of
$E_1^{\textstyle *} \otimes _MK_\zr$. A projection $p\colon E_1\rightarrow
K_\zr$ corresponds, by duality, to an embedding $p^{\textstyle *} \colon
K_\zr^{\textstyle *} \rightarrow E_1^{\textstyle *} $ which is a linear
section of $\zr^0$ and we can repeat all the identification on the dual
side for $\zr^0\colon E_1^{\textstyle *}  \rightarrow K_\zr^{\textstyle *}
$.
                \endproof
        \noindent
        {\bf Remark.} The linear part of the defined above isomorphism between
$\ss{LinSec}(\zr)$ and  $\ss{Proj}(E_2,K_\zr)$ is minus the identity on
$\ss{Hom}_M (E_1,K_\zr)$.  To see it better, let $E_2$ be a subspace of
$E_1$, i.e. $E_1 = E_2 \oplus K_\zr $.  Then a section  $\zs \in
\ss{LinSec}(\zr) $ is identified with a linear mapping
$\overline{\zs}\colon E_2 \rightarrow K_\zr  $ by $\zs(e) = (e,
\overline{\zs}(e))$. Similarly, a projection $\zl \in \ss{Proj}(E_1,K_\zr)$
is identified with  $\overline{\zl} \colon E_2\rightarrow K_\zr $ by
$\overline{\zl}(e) = \zl(e,0)$. In both cases the affine structure is the
canonical affine structure of the vector bundle $\ss{Hom}_M (E_1,K_\zr)$.
But the isomorphism of the theorem is given by $\zl\circ \zs = 0$, i.e.
$\zl(e, \overline{\zs}(e)) =\overline{\zl}(e) + \overline{\zs}(e)= 0$.

                \claim \c{t}{Theorem}{}
        \label{Cbx23}
Assume that $K_\zr$ is generated by a single section $\zf$ of $\zt_1\colon
E_1\rightarrow M$. Then we have the identification
            $$ \zs \longleftrightarrow f_\zs \longleftrightarrow a_\zs,
                                         $$
        where
                $$\align
        \zs&\in \ss{LinSec}(\zr),  \\
        f_\zs=\imath_{a_\zs}&\in \ss{Lin}_\zf(E_1,\R)=\{\ss{Lin}(E_1,\R)\ni f
\colon f\circ \zf =1\},   \\
                        a_\zs&\in \ss{Sec}((E_1,\zf)^\ddag),
                                                                        \tag
\label{Fb66}\endalign$$
        where $(E_1,\zf)^\ddag = \{E_1^{\textstyle *}\ni v\colon
\imath_\zf(v)=1 \}$  and the correspondence between $\zs$ and
$f_\zs$ is determined by the equation
                $$ f_\zs (e)\cdot \zf(\zt_1(e)) = e- \zs\circ \zr(e)
                                         \tag \label{Fb67}$$
        for all $e\in E_1$.
                \endclaim
                \proof
    The correspondence $a_\zs \longleftrightarrow f_\zs $ is obvious. Let
us see that~\Ref{Fb67} defines a unique $f_\zs\in \ss{Lin}_\zf(E_1,\R)$
for any $\zs \in \ss{LinSec}(\zr)$.  Since for any $e\in E_m$ the
right-hand side   of~\Ref{Fb67}    is    an    element    of
$(K_\zr)_m$,    the coefficient $f_\zs(e)$ is well
defined. Moreover, by linearity of $\zr$ and $\zs$,
                $$ f_\zs(e_1 +e_2)\zf(m) = (e_1- \zs\circ \zr(e_1)) + (e_2-
\zs\circ \zr(e_2)) = (f_\zs(e_1) +f_\zs(e_2))\zf(m)
                                         \tag \label{Fb68}$$
         for any $e_1, e_2 \in (E_1)_m$, so $f_\zs$ is linear. Further, by
definition,
        $$ f_\zs(\zf(m))\cdot \zf(m) = \zf(m) - \zs\circ \zr(\zf(m))
=\zf(m),
                                         \tag \label{Fb69}$$
so $f_\zs\circ\zf =1$.

        Conversely, let $f\in \ss{Lin}_\zf(E_1,\R)$. Let us notice that the
mapping $S\colon e\mapsto  e-f(e)\cdot \zf(\zt_1(e))$ is constant along
fibers of $\zr$. Indeed, for $e\in (E_1)_m$ and $v\in (K_\zr)_m$,
$v=a\zf(m)$, we have
        $$\align
        S(e+v)&= (e+v)-f(e+v)\zf(m)   \\
                        &= (e-f(e)\zf(m))+(v-f(v)\zf(m)) \\
                        &= S(e) +a\zf(m)-af(\zf(m))\zf(m)= S(e).
         \tag \label{Fb70}\endalign$$

        Here we have used the linearity of $f$ and the fact that $f\circ \zf
=1$. Since $S(e)$ depends only on $\zr(e)$, we can define $\zs\colon
E_2\rightarrow E_1$ by
                $$ \zs(\zr(e)) = S(e).
                                         \tag \label{Fb71}$$
         Since  $\zr(S(e))=\zr(e)$, the  mapping  $\zs$  is  a  section.
Moreover, the linearity of $f$ implies the linearity of $\zs$.
                \endproof

        Let us see how the above correspondences  look like in local
coordinates. Let us choose local coordinates $(x^a)$ on $M$, linear
coordinates $(x^a,y_i)$ on $E_2$, linear coordinates $(x^a, y_i,z)$ on
$E_1$. The dual coordinates on $E_1^{\textstyle *} $ will be denoted by
$(x^a,\zx^i,\zh)$. Then $K_\zr = \{e\in E_1\colon y_i(e) =0 \}$ is
generated by the section $\zf\colon (x^a) \mapsto (x^a,0,1)$ and
$(E_1,\zf)^\ddag = \{E_1^{\textstyle *} \ni e^{\textstyle *} \colon \zh
=1\}$. For the linear section $\zs(x^a,y_i) = (x^a,y_i, f^i(x)y_i)$
we have $a_\zs(x^a) = (x^a,-f^i(x),1)$ and $f_\zs =
\imath_{a_\zs}(x^a,y_i,z)= -f^i(x)y_i +z$.
This implies easily the following.

                \claim \c{t}{Theorem}{}
         \label{Cby23}
In   the   notation   of   Theorem~\Ref{Cbx23}, $f_\zs=\overline{\zr}^
{\textstyle *}(\zs)$ (i.e. $f_\zs$  is  the  pull-back  of  the  section
$\zs$) with respect to the morphism of the special affine bundles

    $$ \CD  E_1\times\overline{\R} @() \L{\overline{\zr}} @(1,0) @()
\L{pr_1}@(0,-1)& E_1 @()\L\zr @(0,-1) \\
   E_1 @() \L\zr @(1,0) & E_2=E_1/\langle\zf\rangle
         \endCD ,                 \tag \label{F99a}$$
where $\overline{\R}$ is the special vector space $(\R,-1)$ and
$\overline{\zr}(e_m,t)=e_m-t\zf(m)$.
\endclaim

  Let us remark that, completely analogously, one can consider the
affine bundles obtained from affine projections $\zr\colon A_1\rightarrow
A_2$ of affine bundles instead of vector bundles. We get results  similar
to the above ones for affine  sections and functions instead of linear ones.

                \claim \c{t}{Theorem}{}
         \label{Cb24}
   There is a one-to-one correspondence between Lie affgebroid brackets
$[\,,\,]$ on an affine bundle $A$ and linear aff-Poisson brackets
$\{\,,\,\}$ on the 1-dimensional affine bundle
                $$ \zr\colon A^\dag \rightarrow A^\dag/\langle 1_A\rangle
\simeq
\sV(A)^{\textstyle *}
                                         \tag \label{Fbx72}$$
          determined by
                $$ \{\zs_1,\zs_2\} = \imath_{[a_{\zs_1},a_{\zs_2}]}.
                                         \tag \label{Fb73}$$
                \endclaim

\medskip\noindent
{\bf Remark.} Here we call an aff-Poisson structure linear if the
bracket of linear sections is a linear function on $A^\dag /\langle
1_A\rangle $. Since $A^\dag /\langle
1_A\rangle $ can be canonically identified with $\sV(A)^{\textstyle *} $, the
linear aff-Poisson structure restricted to $\ss{LinSec}(\zr)$ gives a
bracket
                $$ \{\,,\,\} \colon \ss{LinSec}(\zr)\times \ss{LinSec}(\zr)
\rightarrow  \ss{Lin}(\sV(A)^{\textstyle *} ,\R).
                                         \tag \label{Fb74}$$
Since, due to Theorem~\Ref{Cbx23},
we have the canonical identification $\ss{LinSec}(\zr) \simeq
\ss{Sec}(A)$ and $\ss{Lin}(V(A)^{\textstyle *} ,\R) \simeq \sV(A)$, any Lie
affgebroid bracket on $A$ corresponds to a unique Lie affgebroid bracket on
$\ss{LinSec}(\zr)$ and vice-versa.
                        \claim \c{l}{Lemma}{}
                          \label{Cb25}
There is a one-to one correspondence between quasi-derivations  $D\colon
\ss{Sec}(E)\rightarrow \ss{Sec}(E)$ of  a vector bundle $E$ over $M$ and
linear vector fields $\overline{D}$ on $E^{\textstyle *} $ such that
                        $$  \imath_{D(X)}= \overline{D}(\imath_X).
                 \tag \label{Fb75}$$
                \endclaim
                \proof
        In local linear coordinates $(x^a,y_i)$ on $E^{\textstyle *} $
associated with a basis $(e_i)$ of sections of $E$ and local coordinates
$(x^a)$ on $M$, we put
 $$ \overline{D} = f_i^j(x) y_j\partial _{y_i} + g^a(x) \partial_{x^a}
                                         \tag \label{Fb76}$$
         where $D(e_i)= f_i^j(x) e_j$ and $\widehat{D} = g_a(x) \partial
_{x^a}$.
        It is a matter of simple calculations to check that the
formula~\Ref{Fb75} is satisfied. On the other hand, \Ref{Fb75} defines
$\overline{D}$ by $D$ (and vice versa) uniquely.
                \endproof
        \demo{Proof of the Theorem}
    Let us choose a linear section $\zs$ of $\zr$ and let $\{\,,\,\}_0$ be
the linear Poisson bracket on $\sV(A)^{\textstyle *} $ corresponding to the
Lie algebroid bracket $[\,,\,]_\sv$ on $\sV(A)$. Let $a_\zs$ be the section
of $A$ corresponding to $\zs$ and let $D=[a_\zs, \cdot ]^\sv$ be the quasi
derivation on $\sV(A)$ corresponding to $a_\zs$, and let $\overline{D}$ be
the linear vector field on $\sV(A)^{\textstyle *} $ corresponding to $D$,
according to Lemma~\Ref{Cb25}.

        Define a bracket on sections of $\zr$ by
                $$ \{\zs + f, \zs +g\} = \overline{D}(g-f) + \{f,g\}_0
                                         \tag \label{Fb77}$$
        for all $f,g \in C^\infty(\sV(A)^{\textstyle *} )$. It is obvious that this is an affine
derivative with respect to each argument separately. Moreover, if $f$ and
$g$ are linear functions, i.e. $f=\imath_{X_f}$, $g=\imath_{X_g}$ for
certain $X_f,X_g \in\ss{Sec}(\sV(A))$, we have
                $$\align
        \{\zs   +\imath_{X_f},\zs   +\imath_{X_g}   \}&=    \overline{D}
(\imath_{X_f}-\imath_{X_g}) +\{\imath_{X_f}, \imath_{X_g}\}_0 \\
                        &=\imath_{D(X_f -X_g)} +\imath_{[X_f,X_g]_\sv}  \\
                        &= \imath_{[a_\zs +X_f, a_\zs + X_g]} \\
                &= \imath_{[a_{\zs +f}, a_{\zs +g}]}.
         \tag \label{Fb78}\endalign$$
        It is clear that~\Ref{Fb73} determines $\{\,,\,\}$ by $[\,,\,]$
uniquely
and vice-versa.
                \endproof
\medskip
    \noindent {\bf Remark.} The Poisson structure on $A^\dag/\langle
1_A\rangle$ can be regarded as a reduction of the Poisson structure on
$A^\dag$ with respect to the canonical projection $\zr\colon A^\dag
\rightarrow A^\dag/\langle 1_A\rangle$ (cf. Theorem 11 (4)). This  means
that
      $$ \{\zr^{\textstyle *} (f), \zr^{\textstyle *} (g)\}^0 =
\zr^{\textstyle *}\{f,g\}_\sv,
                                 \tag\label{F97}$$
   where $\{\,,\,\}^0$ is the Poisson bracket on $A^\dag$ corresponding to
the Lie algebroid structure on $\widehat{A}$.
   Also the aff-Poisson bracket~\Ref{Fb73} can be obtained as a reduced
bracket, but this time with respect to the morphism of the special affine
bundles (cf. Theorem~\Ref{Cby23}):

      $$ \CD  A^\dag\times\overline{\R} @() \L{\overline{\zr}} @(1,0) @()
\L{pr_1}@(0,-1)& A^\dag @()\L\zr @(0,-1) \\
   A^\dag @() \L\zr @(1,0) & \sV(A)^{\textstyle *}
         \endCD .                 \tag \label{F99}$$
   This morphism defines the pull-back  $\overline{\zr}^{\textstyle *}\zs
$ of a section of $\zr$ to a function on $A^\dag$. Then
      $$ \zr^{\textstyle *} \{\zs_1,\zs_2\} =
\{\overline{\zr}^{\textstyle *} \zs_1, \overline{\zr}^{\textstyle *} \zs_2
\}^0.
                                             \tag \label{F98}$$
   We can say that $\{\,,\,\}$ is obtained from $\{\,,\,\}^0$ by an {\it
affine reduction}, in contrast to the standard reduction~\Ref{F97}.

      \claim \c{e}{Example}{}   {\rm   \label{C99}  ([\Ref{U1}])
   Let us consider the case important for the classical mechanics  when
the base  manifold  $M  $  is  the  space-time  fibred  over  the
time represented by the real line  $t\colon M\rightarrow \R$. We
parameterize first-jets of motions by the affine subbundle $A$ of
$\sT M$: $A= \{v\in \sT M\colon\sT t(v) = \frac{\partial
}{\partial t} \}$. We have then $\widehat{A} = \sT M$, $A^\dag =
\sT^{\textstyle *} M$, and $\zf = \xd t$. The Lie affgebroid
structure on $A$ is given by the canonical Lie algebroid structure
on $\sT M$. The quotient manifold $A^\dag/\langle \zf\rangle
=\sV^{\textstyle *} M$ is the dual to the bundle of vertical
vectors on $M$ and represents the phase manifold for a
non-relativistic particle. It carries a canonical Poisson
structure obtained by the linear reduction~\Ref{F97}. Hamiltonian
vector fields related to this structure are $t$-vertical and
cannot represent the dynamics of a particle.

   Similarly, we can associate a vector field $X_\zs$ to a section of
$\zr\colon \sT^{\textstyle *} M \rightarrow \sV^{\textstyle *} M$ by the
formula

        $$\{\zs,f\}^\sv = X_\zs(f), \quad \text{for} \ f\in
C^\infty(\sV^{\textstyle *} M),
                                    \tag \label{F96}$$
where  $\{\,,\,\}$ is the aff-Poisson bracket~\Ref{Fb73}.
   It is easy to check that $X_\zs$ projects to $\frac{\partial
}{\partial t}  $ and can represent the dynamics. The image of the
generating section is the energy surface of the system.
      }\endclaim

                \hl1{Concluding remarks}
\noindent
        In this paper  we have discussed some basic aspects of Lie-like
structures on affine bundles.  Such structures arise in a natural
way in differential geometry and classical mechanics. The next step would
be a more extended analysis  as one can do in the case of Lie algebroids
(see  e.g.  [\Ref{GU3}]).  One  needs  also  a  similar  theory  for  Lie
algebroids on special vector bundles with the affine dual structures as
well a theory for Lie affgebroids on special affine bundles. This would
include the duality between the affgebroid structure on $\bT Z$ and the
contact structure on the first jets bundle ${\ss J}^1(Z)$.
Another open question is  a  groupoid  version  of  Lie affgebroids. We
postpone these problems to a separate paper.

\newcounter\bibno
\def\bb{\bib\key \bibno}
        \makebib
        \bb \label{Be}
\by S. Benenti, S.
\paper Fibr\'es affines canoniques et m\'ecanique newtonienne
\jour S\'eminaire Sud-Rhodanien de G\'eometrie, Journ\'ees S.M.F
\paperinfo Lyon, 26-30 mai, 1986
        \endbib

        \bb \label{GU1} \by Grabowski, J. and Urba\'nski, P. \pp  447--486
\paper
Tangent and cotangent lifts and graded Lie algebras associated with Lie
algebroids \yr 1997 \vol 15 \jour Ann. Global Anal. Geom.
                \endbib
      \bb \label{GU2} \by Grabowski, J. and Urba\'nski, P. \pp 195--208
\paper Lie algebroids and Poisson-Nijenhuis structures \yr 1997
\vol 40\jour Rep. Math. Phys.
                \endbib

        \bb \label{GU3} \by Grabowski, J. and Urba\'nski, P. \pp 111--141
\paper Algebroids -- general differential calculi on vector bundles \yr 1999
\vol 31\jour J. Geom. Phys.
                \endbib
        \bb \label{GM} \by Grabowski, J. and Marmo, G. \pp 10975--10990
\paper Jacobi structures revisited \yr 2001 \vol 34
\jour J. Phys. A
                \endbib
        \bb  \label{IM}  \by   Iglesias, D.   and   Marrero, J.~C.   \paper
Generalized Lie bialgebroids and Jacobi structures
\jour J. Geom. Phys.\vol 40\yr 2001\pp 176--199
                \endbib
        \bb \label{Ko} \by Konieczna, K. \paper Affine formulation of the
Newtonian analytical mechanics (in Polish)\paperinfo Thesis, University of
Warsaw 1995
        \endbib
        \bb\label{KS}\by Kosmann-Schwarzbach, Y.\pp 153--165
\paper Exact Gerstenhaber algebras and Lie bialgebroids\yr 1995
\vol 41\jour Acta Appl. Math.
              \endbib
        \bb \label{Li}\by Libermann, P. \paper Lie algebroids and mechanics
\jour Arch.
Math. (Brno) \vol 32  \yr 1996 \pp 147--162
  \endbib
        \bb \label{Ma} \by Mackenzie, K.~C.~H. \pp 97--147 \paper Lie
algebroids and Lie pseudoalgebras \yr 1995 \vol 27 \jour Bull. London
Math. Soc.
        \endbib
        \bb\label{MX}\by Mackenzie, K.~C.~H. and Xu, P.\pp 59--85
\paper Classical lifting processes and multiplicative vector fields
\yr 1998\vol 49 (2)\jour Quartely J. Math. Oxford
              \endbib
        \bb \label{Mar} \by  Mart\'\i nez, E. \paper Lagrangian mechanics on
Lie
algebroids \jour Acta Appl. Math. \vol 67  \yr 2001 \pp 295--320
                \endbib
       \bb \label{Sa2} \by Mart\'\i nez, E.,  Mestdag, T.,  and Sarlet,W.
\paper Lie algebroid structures and Lagrangian systems on affine
bundles \jour J.~Geom.~Phys. \vol 44 \yr 2002\pp 70--95
                \endbib
                \bb \label{MT}
         \by Menzio, M.R., Tulczyjew, W.M. \paper Infinitesimal symplectic
relations and generalized Ha\-miltonian dynamics  \jour Ann. Inst. H.
Poincare, \vol 29 \pp 349--367 \yr 1978
        \endbib
        \bb \label{Pi}
\by Pidello, G.
\paper Una formulazione intrinseca della meccanica newtoniana
\paperinfo Tesi di dottorato di Ricerca in Matematica,
Consorzio Interuniversitario Nord - Ovest,
1987/1988
        \endbib

        \bb \label{Sa} \by Sarlet,W.,  Mestdag, T.,  and Mart\'\i nez, E.
\paper
Lie algebroid structures on a class of affine bundles \paperinfo
math.DG/0201264
                \endbib
                \bb \label{Sa1} \by Sarlet, W., Mestdag, T., and Mart\'\i
nez, E. \paper Lagrangian equations on affine Lie algebroids \paperinfo
In: {\it Differential Geometry and its Applications}, Proc. 8th Int. Conf.
(Opava 2001), D.~Krupka et al. (Eds.) \toappear
                \endbib

        \bb \label{Tu}
\by Tulczyjew, W.M.
\paper Frame independence of analytical mechanics
\jour Atti Accad. Sci. Torino,
\vol 119
\yr 1985
        \endbib

        \bb \label{TU} \by Tulczyjew, W.~M. and Urba\'nski, P.  \pp 257--265
\paper
An affine framework for the dynamics of char\-ged particles \yr 1992 \vol
126 \jour Atti Accad. Sci. Torino Suppl. n. 2
        \endbib
        \bb \label{U1} \by Urba\'nski, P. \pp 123--129 \paper Affine Poisson
structure in analytical mechanics \inbook Quantization and
Infinite-Dimensional Systems \eds J-P. Antoine and others
\yr 1994 \publ Plenum Press \publaddr New York and London
        \endbib

        \bb \label{We1} \by Weinstein, A. \paper A universal phase space for
particles in Yang-Mills fields \jour Lett. Math. Phys. \yr 1978\pp 417--420
\vol 2
        \endbib
        \bb \label{We2} \by Weinstein, A. \paper Lagrangian mechanics and
groupoids
\jour Fields Inst. Comm. \yr 1996\pp 207--231 \vol 7
        \endbib

\enddocument